\documentclass[a4paper,11pt]{report}

\usepackage{amsmath}
\usepackage{amsfonts}
\usepackage{amssymb}
\usepackage{amsthm}
\usepackage{graphicx}
\usepackage[all]{xy}
\usepackage{setspace}
\usepackage{hyperref} 
\usepackage{changepage}

\newcommand{\Z}{{\mathbb{Z}}}
\newcommand{\C}{{\mathbb{C}}}
\newcommand{\R}{{\mathbb{R}}}
\newcommand{\Q}{{\mathbb{Q}}}
\newcommand{\N}{{\mathbb{N}}}
\newcommand{\p}{{\mathbb{P}}}
\newcommand{\G}{{\mathbb{G}}}

\newcommand{\T}{{\mathbb{T}}}
\newcommand{\D}{{i \partial \bar{\partial}}}
\newcommand{\Do}{{i \partial_0 \bar{\partial}_0}}

\newcommand{\Dt}{{i \partial_t \bar{\partial}_t}}

\theoremstyle{plain}% default
\newtheorem{thm}{Theorem}[section]
\newtheorem{lem}{Lemma}[section]
\newtheorem{prop}{Proposition}[section] %contatori differenti nelle sezioni
\newtheorem{cor}{Corollary}[section]
\newtheorem{conj}{Conjecture}[section]

\theoremstyle{definition} % il corpo dell'enunciato non viene italicizzato
\newtheorem{defi}{Definition}[section]

\newtheorem{quest}{Question}[section]

\newtheorem{rmk}{Remark}[section] %nome ambiente in corsivo, testo normale
\begin{document}
\onehalfspacing

%titlepage

%Preliminary material
\begin{titlepage}

\thispagestyle{empty}
  \begin{adjustwidth}{-2cm}{-2cm}

\begin{center}

%title
\textbf{\Huge{ Degenerations of K\"ahler-Einstein Fano Manifolds} }\\[5cm]  

%author&institutions
\LARGE{Cristiano Spotti}\\
\Large{Department of Mathematics}\\
\Large{Imperial College London}\\

\vfill

% Bottom of the page
\large{\emph{ Thesis presented for the degree of  Doctor of Philosophy}}

\end{center}

  \end{adjustwidth}
\end{titlepage}
\clearpage
\thispagestyle{empty} 

\begin{center}
\Large{\textbf{Declaration}}\\[1cm]
\end{center}

I declare that the material presented in this Thesis is my own work, except where otherwise indicated.

\begin{flushright}\emph{Cristiano Spotti}\end{flushright}

\clearpage

\thispagestyle{empty} 

\begin{center}
\Large{\textbf{Abstract}}\\[1cm]
\end{center}

In this Thesis I investigate how Fano manifolds equipped with a K\"ahler-Einstein metric can degenerate as metric spaces (in the Gromov-Hausdorff topology) and some of the relations of this question with Algebraic Geometry. 

A central topic in the Thesis is the study of the deformation theory for singular K\"ahler-Einstein metrics. In particular, it is shown that K\"ahler-Einstein Fano varieties of dimension two (Del Pezzo surfaces) with only nodes as singularities and discrete automorphism group, admit (partial) smoothings which also carry (orbifold)  K\"ahler-Einstein metrics.
The above result is then used to study the metric compactification in the Gromov-Hausdorff topology of the space of K\"ahler-Einstein Del Pezzo surfaces. In the case of cubic surfaces some evidence is provided that the metric compactification agrees with the classical algebraic compactification given by the set of Chow polystable cubics. 

Finally, I study some higher dimensional analogous of the results outlined above: for example, we briefly discuss the existence and deformation theory for K\"ahler-Einstein metrics on nodal Fano varieties and the  compactifications of the space of intersections of two quadrics in $\p^5$.

\clearpage

\thispagestyle{empty} 

\begin{center}
\Large{\textbf{Acknowledgements}}\\[1cm]
\end{center}

I wish to express my deepest gratitude to my supervisor, Professor Sir Simon K. Donaldson, for generously sharing his ideas, for his infinite patience in explaining them and for his inestimable suggestions.  I am also grateful to my PhD examiners, Dr. Ivan Cheltsov and Dr. Mark Haskins, for their useful comments.

I would like to sincerely thank Professor Claudio Arezzo and Professor Frank Pacard for many helpful conversations, precious teachings and constant encouragement.  My gratitude goes also to all PhD students, Post-Doc researchers and academics of the Geometry Group at Imperial College: the long discussions with people and the beautiful seminars I attended, have been of immense value for my mathematical knowledge. 
In particular, I am grateful to Dr. Alexander Kasprzyk for the help given in a computation based on his database of toric Fano varieties (GRDb).  

I would also like to acknowledge the Department of Mathematics at Imperial College London and EPSRC who provided my financial support during the past three and a half years.
  
Last but not least, I would like to thank  Mamma, Pap\`a and Laura: without their support, it would have been impossible to overcome difficult times and arrive to the end of the PhD.

\clearpage

\tableofcontents

%Chapters

\begin{chapter}*{Introduction}

\addcontentsline{toc}{chapter}{Introduction}

The main subject of this Thesis is the study of compactified  moduli spaces of  K\"ahler-Einstein (KE) metrics on Fano manifolds, with special emphasis on the two dimensional case, i.e., on Del Pezzo surfaces.  Recall that a KE Fano manifold is a K\"ahler manifold $(M,J,\omega)$ with (necessarily)  ample anticanonical line bundle, satisfying the Einstein equation:
$$\mbox{Ric}(\omega)=\omega.$$
For a given Fano manifold it is usually a very difficult problem to detect if a KE metric exists or not: in contrast with the Calabi-Yau or negative first Chern class cases (S-T. Yau \cite{Y78} and, independently,  T. Aubin \cite{A78} for $c_1(X)<0$), not all Fano manifolds can admit a KE metric.
The first obstructions to have been discovered were related to the automorphism group. In \cite{M57}, Y. Matsushima proved that the automorphism group of a KE Fano manifold is always reductive, being given by the complexification of the compact group of isometries. This condition prevents many Fano manifolds, e.g., the blowup of the projective plane in one point, to carry a compatible KE metric.  
Another important obstruction to the existence of such special metrics is given by the vanishing of a character of the Lie algebra of holomorphic vector fields, the so-called (classical) Futaki invariant \cite{F83}.

E. Calabi conjectured that Fano manifolds with discrete automorphism groups always admit a KE metric. In the seminal work \cite{T90}, G. Tian verified this conjecture for smooth Del Pezzo surfaces. A few years later, he discovered in collaboration with  W. Y. Ding \cite{DT92}, that the Calabi Conjecture fails to be true for singular Del Pezzo surfaces. Finally in \cite{T97} the author showed that there are smooth Fano 3-folds with discrete automorphism group (deformations of the Mukai-Umemura Fano 3-fold) which actually cannot admit a KE metric. Thus the discreteness of the automorphism group is not sufficient for the existence of KE metrics.

Following a suggestion of S-T. Yau, G. Tian conjectured that the existence of a KE metric should be equivalent to some algebraic GIT notion of stability. 
Around the same period S. Donaldson (also A. Fujiki) showed how to interpret the constant scalar curvature (cscK) equation in the K\"ahler case (which generalizes the Einstein equation) as a moment map \cite{D97} and thus exhibited a formal connection to some ``stability'' notion via an infinite dimensional analog of the Kempf-Ness Theorem. In a series of seminal papers \cite{D01}, \cite{D05}, S. Donaldson explained how to relate the existence of a cscK metric to the classical GIT notion of Chow stability: more precisely, using the so-called Bergman kernel asymptotics obtained  by G. Tian, S. Zelditch, Z. Lu and others, he proved that if $X$ admits a cscK metric and $Aut(X,L)$ is discrete then $X$ needs to be asymptotically Chow stable.  

In \cite{D02} the notion of (algebraic) K-stability  is introduced. This  notion of stability is based on a generalization of the Futaki invariant (the so-called Donaldson-Futaki invariant $DF$) and it extends  Tian's original K-stability. In order to test K-stability one needs to consider a-priori all possible one parameter subgroup degenerations of a given variety (test-configurations). In the case of Fano manifolds, C. Li and C. Xu recently showed that one can restrict to test-configurations where the central fiber is a normal Fano variety with log-terminal singularities \cite{LX11}. Nevertheless testing K-stability remains a difficult problem, e.g., an algebraic proof of the K-polystability of $\p^2$ is still not known!

A generally accepted way to phrase the relations between KE metric and ``stability'' is the following: 

\begin{conj}[Yau-Tian-Donaldson (YTD)]\label{YTD} Let $X$ be a Fano manifold. Then $X$ admits a KE metric if and only if $X$ is K-polystable.
\end{conj}

Recall that $X$ is called K-polystable if $DF(X)\geq 0$ for every (normal) test configuration with equality if and only if the test-configuration is the trivial one. 

The direction KE implies K-polystability is now settled by the work of S. Donaldson \cite{D05b}, J. Stoppa \cite{S09} and T. Mabuchi \cite{M08} \cite{M09}. Very recently R. Berman has posted a preprint on ArXiv \cite{B12}, where it is shown that the direction KE implies K-polystability holds also for some singular Fano varieties.
On the other hand we should remark that it may be possible that one needs to restrict the class of ``stable'' Fano varieties to make the conjecture true (i.e., also the correct notion of ``stability'' one needs to choose is part of the conjecture). A refined notion of stability has been recently considered by  S. Donaldson in \cite{D10} and G. Sz\'ekelyhidi in \cite{SZ12}. However a counterexample to the YTD conjecture, as above stated, is not known.

In the previous paragraphs, we have very briefly (and incompletely) summarized  the state-of-art of  the ``KE problem for Fano manifolds'', focusing on the relations with algebraic notions of stability. Now we turn to the main object of the Thesis:  moduli spaces of KE Fano manifolds and their compactifications.

On a given smooth manifold we can sometimes find continuous families of Fano complex structures. Some of these structures could admit KE metrics. It is known that the set of KE Fano metrics, equipped with the so-called Gromov-Hausdorff (GH) topology is pre-compact. Using metric geometry, one can consider the space of KE structures on a fixed manifold $M$ together with their GH degenerations: $\overline{\mathcal{M}_M}^{GH}$. This space is compact and Hausdorff. By the works of G. Tian \cite{T90} and M. Anderson \cite{A89}, it is known that GH limits of smooth KE Del Pezzo surfaces are KE Del Pezzo orbifolds with isolated singularities. Moreover, in the very recent work of S. Donaldson and  S. Sun \cite{DS12} it is shown that GH limits of KE Fano manifolds are always singular normal Fano varieties, with at most log-terminal singularities.
Thus $\overline{\mathcal{M}_M}^{GH}$ is a compact Hausdorff space ``parameterizing'' possibly singular (KE) Fano varieties.

 The YTD conjecture implies that the above ``metric'' discussion should have an analogous in the algebraic geometric world, i.e., the GH compactification $\overline{\mathcal{M}_M}^{GH}$ should parameterize ``stable'' algebraic varieties (K-polystable?). It follows that one expects to find a \emph{compact, projective, algebraic} variety $\overline{\mathcal{M}_M}^{ALG}$ (a ``coarse moduli space of stable smoothable Fano varieties'') which is \emph{homeomorphic}, maybe up to some finite cover, to the GH compactification, i.e., $\overline{\mathcal{M}_M}^{GH} \cong \overline{\mathcal{M}_M}^{ALG}$.

The only situation where this identification is known, is the case of Del Pezzo quartics, studied by T. Mabuchi and S. Mukai \cite{MM90}. In this Thesis we  discuss the GH compactification of KE Del Pezzo surfaces of lower degrees. Even if we do not provide a definitive answer, we show what the GH compactification should be in the case of cubics surfaces and degree two Del Pezzo surfaces. 
We should remark that the exact types of singularities appearing in GH degenerations of Del Pezzo surfaces are still not known: establishing a link between some notion of stability (maybe an ``ad hoc'' stability notion related to the geometric realization of a Del Pezzo surface of a given degree) and the GH compactification can be useful to \emph{classify} the singular KE Del Pezzo surfaces (and their singularities) appearing in the ``boundary'' of the GH compactified moduli space.

Another important theme of the Thesis, which is somehow connected with the previous discussion about the relations between metric and algebraic compactifications, is the study of the (metric) deformation theory of singular KE Del Pezzo surfaces (or more generally of singular Fano varieties). This investigation is essential in the study of the KE moduli space at its boundary points. In this direction we prove that generic deformations (partial-smoothings) of nodal Del Pezzo surfaces with discrete automorphism group always admit KE metrics (close to the original one). The proof of the Theorem is based on a gluing construction which ``reverses'' the GH degeneration picture. Similar gluing constructions have been considered by many authors, e.g., D. Joyce \cite{J96}, C. Arezzo and F. Pacard \cite{AP06}, S. Donaldson \cite{D10g} and O. Biquard and V. Minerbe \cite{BM11}.  We also describe the expected general picture in the case of deformations of KE Del Pezzo orbifolds with  more general singularities and with continuous family of automorphisms. These considerations should clarify what the GH compactification is in the case of Del Pezzo surfaces.  

Finally, we  discuss a higher dimensional example of the relations between GH degenerations and stability notions: the case of the intersections of two quadrics in $\p^5$.

In the next section we describe the principal topics discussed in this Thesis.

\section*{Brief outline of Chapters content}

In \emph{Chapter 1} we describe some basic properties of the moduli space of KE Fano structures. We begin by pointing out the precise set-theoretical relations between complex structures and KE metrics. For example, we observe that in general the metric structure of a KE Fano manifold determines the complex structure up to complex conjugation. Then we recall how to  compactify the space of KE Fano structures using the GH distance.
After stating a guideline conjecture on the relation between GH metric compactifications and algebraic compactifications, we will discuss how to study the local picture of the KE Fano (compactified) moduli problem.  We also recall some results concerning ``embedded'' degenerations of KE Fano manifolds, i.e., relations between GH limits and flat limits in a fixed projective space.

In \emph{Chapter 2} we focus our attention on the $2$-dimensional case, i.e., (singular) Del Pezzo surfaces, giving some examples of  KE metrics on singular Del Pezzo surfaces and showing some of their properties. Then we discuss what should happen around smoothable KE Del Pezzo surfaces. Finally we construct two (toric) examples of smoothable KE Del Pezzo surfaces in degree one and two which have strictly Kawamata log-terminal singularities (i.e., log-terminal but not canonical) and we discuss the classification of $\Q$-Gorenstein smoothable KE toric log Del Pezzo surfaces.

 \emph{Chapter 3} contains the main result of the Thesis. We prove that generic (partial)-smoothings of KE Del Pezzo surfaces with discrete automorphism group and only nodal singularities admit (orbifold) KE metrics which are close in the GH sense to the original singular metric. As a first application we construct new examples of singular KE Del Pezzo surfaces (e.g., on the deformations of the Cayley cubic $xyz+yzt+ztx+txy=0$).

In \emph{Chapter 4} we discuss compactified moduli spaces of Del Pezzo surfaces. We recall the Mabuchi-Mukai example of Del Pezzo quartics and we carefully study the geometry of their compactification. Then we discuss the case of cubic surfaces and we show how to use the deformation theory of singular Del Pezzo surfaces to understand the GH compactification. Finally, we briefly discuss the two remaining cases of Del Pezzo surfaces of degree one and two.

In the final \emph{Chapter 5} we consider possible higher dimensional generalizations of the results obtained in the previous Chapters. We briefly discuss existence (and obstruction to existence) of special classes of singular metrics on singular Fano varieties. Then we study the moduli space of intersections of two quadrics in $\p^5$ and its (possible) relation with GH degenerations.

\end{chapter}
\begin{chapter}{Moduli of K\"ahler-Einstein Fano structures}

In the first part of this Chapter we describe the theory of moduli spaces of KE Fano structures. After giving the main definitions, we study the set-theoretic relations between Fano complex structures and KE Riemannian metrics. Then we recall the notion of Gromov-Hausdorff (GH) distance and we explain its importance for metrically ``compactifying'' the moduli space of KE Fano structures. We also very briefly recall the main results concerning the ``structure'' of GH limits of KE Fano varieties.

Next we state a guideline Conjecture relating the GH compactified moduli space to the ``mythical'' algebraic compactification of the set of ``stable'' Fano varieties.

We end the Chapter by discussing the local structure of the KE (compactified) moduli space.

\section{Main definitions and first properties}

Let $M^{2n}$ be a smooth compact manifold of real dimension $2n$. We start by defining the set of \textit{Fano structures }on $M$, i.e.,
$$\mathcal{F}_M:=\left\{J \in End(TM)\, | \, J\, \mbox{ integrable complex structure,}\, K_{(M,J)}^{-1} \mbox{ is ample} \right\}.$$
We remark that, in reality, we are primarily interested in the connected components of the above space, that is in the subsets of complex structures that can be joined by smooth paths (w.r.t. the natural smooth topology induced on the tensor field after fixing any Riemannian metric on $M$).

Now we define the set of \textit{K\"ahler-Einstein (KE) Fano structures} as a subset of the above space:
$$\mathcal{F}_M^{KE}:=\left\{J \in \mathcal{F}_M \,|\, \exists\, \omega\in 2\pi c_1(M,J) \mbox{ s.t. } \mbox{Ric}(\omega)=\omega \right\}$$

For our purpose, it is essential to relate the set of KE Fano structures $\mathcal{F}_M^{KE}$ to the space of Riemannian metrics. Before establishing such connection, let us remark that the group of  diffeomorphisms $\mbox{Diff}(M)$ acts in the obvious way on the set $\mathcal{F}_M^{KE}$ as well as on the set $\mathcal{G}_M$ of Riemannian metric tensors on $M$. Then we can state the following:

\begin{prop}\label{D} There exists a well-defined map
$$D: \mathcal{F}^{KE}_M / \mbox{Diff}(M) \longrightarrow \mathcal{G}_M / \mbox{Diff(M)},$$
i.e., up to diffeomorphisms, to each $J\in\mathcal{F}^{KE}_M $ there exists a \textit{unique} K\"ahler-Einstein metric structure.
\end{prop}

\begin{dimo} Let $(M,J_1,\omega_1, g_1)$ and $(M,J_2,\omega_2, g_2)$ be two KE Fano structures such that $\phi^*{J_2}={J_1}$ for $\phi \in \mbox{Diff}(M)$. We claim that there exists a Riemannian isometry between the Riemannian metrics $g_2$ and $g_1$.

By the invariance of the Einstein equation under the action of the diffeomorphism group, $\phi^* \omega_2 \in c_1(M,\phi^* J_2)=c_1(M,J_1)$ is another K\"ahler (Einstein) form in $c_1(M,J_1)$. Then the uniqueness result for  KE metrics on Fano manifolds of S. Bando and T. Mabuchi \cite{BM87} implies that there is $\psi \in \mbox{Aut}_0(M,J_1)$ (the connected component of the biholomorphism group) such that $\psi^{*} (\phi^{*} \omega_2)=\omega_1$.

Thus $\phi \circ \psi \in \mbox{Diff}(M)$ is the desired  Riemannian isometry between $g_2$ and $g_1$.\qed \\

\end{dimo}

Before proceeding, some important remarks are needed. It is well known that the space of complex structures modulo diffeomorphisms is from the topological viewpoint very badly behaved. In particular, it often fails to be Hausdorff. On the other hand, it is a fact (compare \cite{B87}) that the space of Riemannian metrics modulo diffeomorphism satisfies the Hausdorff property. It is then interesting and natural to ask if the set of KE Fano structures modulo diffeomorphism $\mathcal{F}_M^{KE}/ \mbox{Diff}(M)$ is still not separated. We conjecture that  $\mathcal{F}_M^{KE}/ \mbox{Diff}(M)$ is actually Hausdorff: this is connected to the fact that the existence of a KE metric should be related to a stability condition of the complex structure, which, not too unrealistically, could be used to construct Hausdorff moduli spaces of complex varieties (maybe as a GIT quotient).

The second remark, which is somewhat connected to the above discussion, is that Proposition \ref{D} is not saying that the map $D$ is continuous. As we will explain later in this Chapter (compare Theorem \ref{D1}), it is relatively  easy to show continuity at point $[M,J]$ where $Aut_0(M,J)={1}$. However in the presence of  continuous families of automorphisms the situation is more subtle and linked again to stability considerations (Theorem \ref{D2}).

As we pointed out,  $D(\mathcal{F}_M^{KE}/ \mbox{Diff}(M))$ inherits a natural Hausdorff topology from the (quotient) space $\mathcal{G}_M / \mbox{Diff}(M)$. It follows by works of N. Koiso \cite{K78}, \cite{K83} on the local structure of the Moduli Space of Einstein metrics that  the space $D(\mathcal{F}_M^{KE}/ \mbox{Diff}(M))$ is locally given by the quotient of a finite dimensional real analytic variety by the action of a compact group. To be more precise, we have the following:

\begin{prop} [N. Koiso] \label{K} $D(\mathcal{F}_M^{KE}/ \mbox{Diff}(M))$ is naturally a Hausdorff topological space, in general non-compact, which locally at $D[(M,J,\omega,g)]$ is given by
$$\mathcal{Z} / \mbox{Iso}(M,g),$$
where $\mathcal{Z}$ is the zero locus of real analytic functions defined on a finite dimensional vector space.
\end{prop}

A detailed proof of the above statement can be found in A. Besse \cite{B87} (see Chapter $12$) or in the original papers of N. Koiso.  The main ingredient of the proof consists in showing that, after breaking the gauge symmetries (Ebin's slice theorem), the linearization of the (K\"ahler) Einstein equation is elliptic. Then the result follows by obstruction considerations and by the elimination of the residual symmetry given by the isometry group.

From the moduli point of view it is important to study the fibers of the map $D$ defined in Proposition \ref{D}, i.e., to understand when two different complex structures $I$ and $J$ are associated to the same Riemannian metric $g$.  For example the metric structure $g$ is fixed by the involution sending $J$ to $-J$. In the next Proposition we show  that this is essentially the only possibility for two KE Fano complex structures to be associated with the same metric structure.

\begin{thm}[KE Fano Splitting] \label{S}
Let $(M,g)$ be a KE Fano manifold with respect to two complex structures $I$, $J$ where $I \neq\pm J$. Then $M$ splits holomorphically and isometrically as
$$M=M_1\times \dots \times M_k$$
for $k \geq 2$, where $M_i$ are lower dimensional KE Fano manifolds.
\end{thm}

\begin{dimo}

Let $\omega_I:= g(I-,-)$ and $\omega_J:= g(J-,-)$ be the two real K\"ahler forms associated to $I$ and $J$. 

$\mathbf{Claim }$: $\omega_J$ is of type $(1,1)$ w.r.t. $I$, i.e., $\omega_J(I-,I-)= \omega_J(-,-)$ (and analogously for $\omega_I$).

Since $\omega_J$ is a K\"ahler form, $\Delta_g \omega_J=0$. Thus, decomposing $\omega_J$ as
$$\omega_J= \omega_J^{(2,0)_I}+ \omega_J^{(1,1)_I}+\omega_J^{(0,2)_I},$$
it follows by the K\"ahler identities that  the $\omega_J^{(2,0)_I}$ component is also harmonic (and similary for $\omega_J^{(0,2)_I}$). However we have that $\mathcal{H}^{(2,0)_I}(M)=H^0(M_I, \Omega_{M_I}^2)=0,$ by Kodaira's Vanishing. Hence  $\omega_J^{(2,0)_I}=\omega_J^{(0,2)_I}=0$.

Now we define the following ``splitting'' endomorphism $S \in End(TM)$:
$$S_p(v):= (i_v \omega_I)^{\sharp_{\omega_J}},$$
i.e., $\omega_J(S_p(v),w)=i_v\omega_I(w)=\omega_I(v,w)$ for all $v,w \in T_pM$.

$S$ satisfies the following properties:
\begin{itemize}
\item $S\neq \pm 1$, which is an immediate consequence of $\omega_I \neq \pm \omega_J$;
\item $\mbox{ker}\,S_p=\{0\}$, by the non-degeneracy of the K\"ahler forms;
\item $\nabla S=0$, being $\omega_i$ parallel; 
\item $[S,I]=0$. Using the fact that $\omega_J$ is of type $(1,1)_I$, we have that for all $v,w \in T_pM$:
$$\begin{array}{ccl}
\omega_J(S(Iv),w) &=&\omega_I(Iv,w) \\
&=&-\omega_I(v,Iw)\\
&=&-\omega_J(Sv,Iw)\\
&=&\omega_J(I(Sv),w).
\end{array}
$$
\item $g(S_pv,w)=g(v,S_pw)$ for all $v,w \in T_pM$, since
$$\begin{array}{ccl}
g(Sv,w) &=&\omega_J(Sv,Jw) \\
&=&\omega_I(v,Jw)\\
&=&-\omega_I(Jv,w) \\
&=&\omega_I(w,Jv)\\
&=&\omega_J(Sw,Jv)\\
&=&g(Sw,v)\\
&=& g(v,Sw).
\end{array}
$$

\end{itemize}

The above properties imply that $S_p$ is diagonalizable and that there must be \textit{at least two} distinct non-trivial orthogonal eigenspaces, which are also complex subspaces of $T_pM$. Moreover, since $S$ is parallel, the holonomy representation preserves these eigenspaces, that is $\mathcal{H}ol_p(g)$ is \textit{reducible} and
$$\mathcal{H}ol_p(g)\subseteq U(n_1)\times \dots \times U(n_k),$$
where $n_i$ denotes the complex dimension of the eigenspaces and $k\geq2$. 

Finally, recalling that every Fano manifold is simply connected (compare, for example, \cite{B06} Theorem 6.12), the De Rham  Theorem  (e.g., \cite{B06} Theorem 4.76) implies  that $M$ decomposes globally isometrically (and holomorphically) as 
$$M=M_1\times \dots \times M_k,$$
where, of course, each factor is a KE Fano manifold.
\qed \\
\end{dimo}

\begin{rmk}
We believe that the above Theorem should be known to experts in the field. However we have decided to give a complete proof, since we were unable to find such statement in the literature.  
\end{rmk}

 We define an \emph{irreducible Fano manifold} to be a Fano manifold that does not admit a holomorphic splitting into two lower dimensional Fano manifolds. Then,  as a corollary of the above Theorem, we immediately have
 
 \begin{cor} On the locus of irreducible KE Fano manifolds, the map $D$ defined in Proposition \ref{D} is ``generically" $2:1$.
 \end{cor}
 
  The map $D$ could have some ``ramification" locus, that is there are KE Fano complex structures $J$ for which there exists diffeomorphisms $\varphi \in \mbox{Diff}(M)$ satisfying $-J= \varphi^{*} J$ (i.e. $\varphi$ is anti-holomorphic). This is the reason why the map $D$ will be ``generically''  $2:1$.   For example, if $(M,J,g)$ admits an embedding in $\C\p^n$ such that its defining ideal admits as generators polynomials with real coefficients, then the restriction of complex conjugation of the homogeneous coordinates is an example of an anti-holomorphic diffeomorphism. 

In Chapter $4$ we will discuss carefully an example where the map $D$ can be completely described (Del Pezzo surfaces of degree $4$). In particular it can be shown that the involution sending $J$ to $-J$ is antiholomorphic with respect a natural complex structure defined on the moduli space of (smooth) degree $4$ Del Pezzo surfaces.

\section{Metric compactification}

As we have previously pointed out, the space $D(\mathcal{F}_M^{KE}/ \mbox{Diff}(M))$ is in general non-compact. Nevertheless, it admits in a certain sense a \textit{natural} compactification. First of all, we continuously embed the space of Riemannian metric structures $\mathcal{G}_M/ \mbox{Diff}(M)$ in the bigger space of isometry classes of compact metric spaces  $\texttt{Met}$, equipped with the so-called Gromov-Hausdorff (GH) distance. Secondly, we  prove that under this embedding the subspace  $D(\mathcal{F}_M^{KE}/ \mbox{Diff}(M))$ is pre-compact, i.e., its closure is a compact subspace of $\texttt{Met}$.

Let us start by recalling the definitions of Gromov-Hausdorff  distance and some basic properties (compare \cite{BBI01} for more information). Let $(X,d_X)$ and $(Y,d_Y)$ be two compact metric spaces. Then
$$d_{GH}(X,Y):= \inf\{ d_{H}^{Z}(X,Y) \, | \, X,Y\hookrightarrow Z \, \mbox{ isometrically} \},$$
where $d_{H}^Z$ denotes the Hausdorff distance between subset of the metric space $(Z,d)$.
It is a well-known fact that the above (pseudo)-distance defines a metric on the set of isometry classes of compact metric spaces  $\texttt{Met}$.

In general, it is  important to detect when two metric spaces are close to each other in the GH sense without explicitly compute the distance. In particular the following well-known criterion for GH closeness is very useful.

\begin{lem} \label{GHCC} Let $(X,d_X)$ and $(Y,d_Y)$ be two compact metric spaces.  If $d_{GH}(X,Y) \leq \epsilon$ then there exists a $3\epsilon$-quasi isometry $F:X\rightarrow Y$, i.e.,  a non necessarily continuous map $F:X \longrightarrow Y$  satisfying:
        \begin{itemize}
         \item $|d_X(p,q)-d_Y(F(p),F(q))| \leq 3\epsilon$ for all $p,q$ in $X$
         \item $F(X)$ is $3\epsilon$-dense in $Y$.
        \end{itemize}
Conversely, if there exists a $\epsilon$-quasi isometry $F:X\rightarrow Y$ then $d_{GH}(X,Y) \leq 3\epsilon$.
 
\end{lem}

An important application of the above Lemma is given by the following (well-known) observation.

\begin{lem}The natural map
$$i:\mathcal{G}_M/ \mbox{Diff}(M) \longrightarrow \texttt{Met},$$
given by associating to a Riemannian metric $g$ its distance function $d_g$, is a continuous embedding.
\end{lem}

\begin{dimo} Suppose that $[g_i]\rightarrow[g]$ in $\mathcal{G}_M/ \mbox{Diff}(M)$, i.e., there exists a sequence $(\phi_i) \subseteq \mbox{Diff}(M)$ such that $||\phi_i^* g_i-g||_{C^\infty}\rightarrow 0$ (computed by taking an arbitrary fixed background Riemannian metric). Then it follows that the identity map on $M$ is very closed to be an isometry  between $(M,d_{\phi_i^* g_i})$ and $(M, d_g)$.
Thus $i((M,g_i)) \rightarrow i((M,g))$ in the Gromov-Hausdorff sense.

The fact that the map $i$ is an embedding (i.e. is injective), is the classical Theorem of S. Myers, N. Stenrood and R. Palais (for a proof, compare \cite{KN96}, Chapter IV, Theorem 3.10) which says that metric isometries between smooth Riemannian metrics are Riemannian isometries.\qed \\

\end{dimo}

Observe that we are not claiming that $i$ is also open. However $i$ turns out to be open when restricted to special classes of metrics (i.e., metrics satisfying some equations). For example this is the case for positive Einstein metrics with injectivity radius bounded from below \cite{AC92}. 

The embedding of Riemannian manifolds in the larger space of metric spaces enables us to talk about the convergence of manifolds to singular spaces. In this direction,     a foundational result is the following Theorem due to M. Gromov (for a proof, see \cite{BBI01} Chapter $7$, Section $4$):

\begin{thm}[Gromov's precompactness] \label{GP} Let $\texttt{S}\subseteq \texttt{Met}$ be a uniformly totally bounded class of compact metric spaces, i.e.,
\begin{itemize}
 \item $\exists \,C >0$ such that $\mbox{Diam(X)}\leq C$ for all $X \in \texttt{S}$;
 \item for all $\epsilon>0$ there is $n=n(\epsilon)$ such that all $X\in \texttt{S}$ admit an  $\epsilon$-dense set with cardinality less then  $n(\epsilon)$.
\end{itemize}
 Then $\texttt{S}$ is GH-precompact. 
\end{thm}

We are now ready to recall the fundamental result concerning the pre-compactness of the set of Einstein metrics of positive scalar curvature. In our setting, this well-known result can be state as follow:  

\begin{prop} Using the previous notations:
$$i(D(\mathcal{F}_M^{KE}/ \mbox{Diff}(M)))$$
is Gromov-Hausdorff pre-compact.
\end{prop}

\begin{dimo} We need to show that $i(D(\mathcal{F}_M^{KE}))$ is uniformly totally bounded. Then the statement follows by Gromov's pre-compactness \ref{GP}.
By Myers' Theorem, the diameters are uniformly bounded. Next, let $S_{2\epsilon}$ be a finite $2\epsilon$-dense set in $(M,J,\omega)$. W.l.o.g. we may assume that $B(p,\frac{\epsilon}{2}) \cap B(q,\frac{\epsilon}{2}) =\emptyset$ for all $p,q \in S_{2\epsilon}$. All we have to show is a uniform bound on the cardinality of $S_{2\epsilon}$. 
Using Bishop-Gromov volume comparison, we find
$$|S_{2\epsilon}| \leq \frac{Vol(M)}{Vol(B(p, \frac{\epsilon}{2}))} \leq \frac{Vol(S^{2n}(C))}{Vol(B(\frac{\epsilon}{2}))},$$
where the last term in the above expression depends only on $\epsilon$ and on the diameter bound.

\qed\\
\end{dimo}
To fix some notations, let us call
\begin{itemize}
\item $\mathcal{M}_{M}:=i(D(\mathcal{F}_M^{KE}/ \mbox{Diff}(M)))$;
\item $\overline{\mathcal{M}_{M}}^{GH}:=\overline{i(D(\mathcal{F}_M^{KE}/ \mbox{Diff}(M)))}^{GH}$;
\item $\partial{\overline{\mathcal{M}_{M}}^{GH}}:=\overline{\mathcal{M}_{M}}^{GH} \setminus \mathcal{M}_{M}$.
\end{itemize}
On these spaces we consider the natural topology induced by the GH distance. By  definition, ${\overline{\mathcal{M}_{M}}^{GH}}$ is Hausdorff and compact. We remark that the ``boundary''  $\partial{\overline{\mathcal{M}_{M}}^{GH}}$ shouldn't be considered a real boundary of the compactified moduli space: for reasons which will be clear later, the boundary should have real codimension at least two in $\overline{\mathcal{M}_{M}}^{GH}$. Thus a better (and correct) way to think about the boundary, is given by considering it  a ``divisor'', in the sense of Algebraic geometry, in the compactified moduli space.

The two fundamental questions to ask on the compactified moduli space $\overline{\mathcal{M}_{M}}^{GH}$ are the following:
\begin{itemize}
 \item \emph{What are the objects parameterized by $\partial{\overline{\mathcal{M}_{M}}^{GH}}$?} A priori we only know that they are compact (length) spaces. However, since a sequence of $n$-dimensional  KE Fano manifolds is always non-collapsing (this follows  immediately by  Bishop-Gromov's volume comparison), it is not unrealistic to  believe (and in fact it is true \cite{DS12}) that the GH limits also carry the structure of an algebraic variety of dimension $n$ which is compatible with the metric structure (e.g., K\"ahler away from the singularities).
 \item \emph{Can we understand ${\overline{\mathcal{M}_{M}}^{GH}}$ ``globally''?} For example, it is interesting to explicitly identify this topological space and classify the varieties parametrized by it.
\end{itemize}

We will discuss more carefully the second question in the next section. For the moment, we restrict our attention to the first problem, i.e., to the study of the ``structure'' of GH limits. 

The answer to this question has been known for a long time in the case of complex dimension two, i.e., for KE Del Pezzo surfaces:

\begin{thm}[M. Anderson \cite{A89}, G. Tian \cite{T90} Theorem 5.2] \label{AT90} Let $(X_i,\omega_i)$ be a sequence of KE Del Pezzo surfaces. Assume that $$(X_i,\omega_i) \longrightarrow X_\infty,$$
in the GH topology. Then $(X_\infty, \omega_\infty)$ is an irreducible KE Del Pezzo orbifold with isolated singularities of the form $\C^2 / \Gamma$, where $\Gamma$ is a finite group of $U(2)$ acting freely on $\C^2$  minus the origin. 
Moreover these orbifold singularities must be of type $T$ (compare Chapter $2$ Theorem \ref{TS} for a definition).
\end{thm}

The proof of this Theorem is based on an $\epsilon$-regularity result which provides $L^\infty$ controls of the full curvature tensor on balls of radius $\frac{r}{2}$ in term of the $L^2$ norm of the curvature on balls of radius $r$, once this $L^2$ norm is small enough. The crucial fact is that  the total $L^2$ norm of the curvature is a topological constant for KE Del Pezzo surfaces. Then by a standard covering argument one can show that the singularities are isolated. Finally one shows, using a removable singularities Theorem (\cite{T90}, Section $4$), that the metric extends smoothly  (in the sense of orbifolds).

However little is known about the exact types of singularities of GH limits. It is known \cite{T90} that the possible singularities need to be of type $T$ (i.e., they must admit a $\Q$-Gorenstein smoothing). We will discuss  these singularities (and recall the definitions) in the next Chapter.

In higher dimension the situation is much more complicated. We point out that GH limits could have  singularities which are not of orbifold type: to be more precise, orbifold singularities should appear only in complex codimension two. Using the Cheeger-Colding-Tian theory of non-collapsing Riemannian manifolds with Ricci bounded from below \cite{CCT02}, S. Donaldson and S. Sun recently proved in \cite{DS12} that GH limits of KE Fano manifolds are singular Fano varieties ($\Q$-Fano, see the definition in the next section). Very little is known about the behavior of the metric of GH limits near the singularities (especially in complex codimension bigger then or equal to three).

\section{Main Conjecture}
As we have already pointed out several times, the existence of a KE metric on a Fano manifold  is conjectured to depend only on some algebraic geometric properties  (stability)  of the complex structure. With this expectation in mind, it is very  natural to  conjecture that the set of all KE Fano metrics on a fixed smooth manifold, together with its compactification made of GH degenerations, forms a space which admits an algebraic structure. Moreover it is not too unrealistic to  believe that this set can be ``explicitly computed'' in some situations (via some finite dimensional GIT construction).  

In order to make this Conjecture more precise, let us recall some definitions. 

\begin{defi} A reduced, irreducible, normal variety $X$ over $\C$ is called a $\Q$-\emph{Fano variety}  if the following properties hold:
\begin{itemize}
\item $-K_X$ is an ample $\Q$-Cartier divisor;
\item $X$ has at most (Kawamata) log-terminal singularities (compare Appendix B).
\end{itemize}
\end{defi}

\begin{defi} We say that a $\Q$-Fano variety $X$ has Hilbert polynomial equal to $h_X$ if
$$ \mbox{dim}\,H^0(X,-lK_X)= h_X(l),$$
for sufficiently big $l$.
\end{defi}

Recall that $ \mbox{dim}\,H^0(X,-lK_X)= h(l)=\mathcal{X}(-lK_X)$ for all $l$, by Kawamata-Viehweg vanishing Theorem (compare the first Chapter of \cite{IP99}). Moreover note that in dimension two and three the Hilbert polynomial of a (smooth) Fano depends only on the degree, $deg(X):=(-K_X)^n$:
\begin{itemize}
 \item $n=2$, $h_X(k)=\frac{k(k+1)}{2} deg(X) +1$;
 \item $n=3$, $h_X(k)=\left(\frac{k^3}{6}+\frac{k^2}{4}+\frac{k}{12}\right) deg (X) +2 k +1.$
\end{itemize}
 The above formulas follow immediately by standard Riemann-Roch computations.

We are now ready to state the guideline Conjecture (basically a reformulation of the YTD Conjecture \ref{YTD}): 

\begin{conj}\label{MC}
Let $\overline{\mathcal{M}(M^n,h)}^{GH}$ be the GH compactification of the space of isomorphism classes of KE metrics on Fano manifolds with Hilbert polynomial equal to $h$ and with $M^n$ fixed underlying differentiable manifold. 

Then there exists a natural continuous ramified covering map
$$D: \mathcal{X} \longrightarrow  \overline{\mathcal{M}(M^n,h)}^{GH},$$
where $\mathcal{X}$ is a \emph{compact projective algebraic variety}.

Moreover the following properties  hold:
\begin{enumerate}
 \item $\exists \, \mathcal{U} \subseteq \mathcal{X}$ Zariski open such that $D_{|_{\mathcal{U}} }$ is onto the set of smooth KE Fano structures on $M^n$;
  \item $\mathcal{X} \setminus \mathcal{U}$ is a subvariety and  $D_{|_{\mathcal{X} \setminus \mathcal{U}} }$ is onto the set of singular KE Fano varieties appearing as GH degenerations. More precisely, every GH degeneration admits the structure of a $\Q$-Fano variety;
\item $\mathcal{X}$ should be a coarse moduli space for ``\emph{stable}'' (e.g. K-polystable) $\Q$-Gorenstein smoothable $\Q$-Fano varieties with Hilbert polynomial equal to $h$, i.e., $$\mathcal{X}=: \overline{\mathcal{M}(M^n,h)}^{ALG}.$$
\end{enumerate}

\end{conj}

It is important to make some comments on this Conjecture and to explain a bit more precisely the objects involved in it.

\begin{itemize}
\item The map $D$ in the above Conjecture  should be given \emph{a-posteriori} by associating to each point $p \in \mathcal{X}$ the metric space $[X_p,d_g]$ induced by the unique (up to automorphism)  KE metric on the variety $X_p$.  This map should be continuous with  respect to the  analytic topology of the moduli space $\mathcal{X}$ and the GH topology of the target.

\item The reason why  we believe that in general the map $D$ is only a covering map (with finite fibers), and not an homeomorphism, basically has to do with the fact that the metric structure is insensitive to the operation of complex conjugation (i.e., sending the complex structure $J$ to $-J$); compare Proposition \ref{S}.   A way to repair this covering issue could be to consider ``marked'' KE Fano metric spaces (i.e., one should also remember the complex structure).

\item Point $1$ in the Conjecture is related to the fact that it is believed that the set of \emph{smooth}  KE structures is itself of algebraic nature (quasi-projective).

\item Point $2$ is related with the structure theory of GH  degenerations of Fano manifolds briefly discussed at the end of the previous section.

\item Finally, point $3$ is essentially the $YTD$ Conjecture \ref{YTD} relating the  existence of KE metrics to some algebraic notion of stability (e.g., K-polystability). It is natural to believe that the Conjecture extends to singular Fano varieties. One central problem is to detect which is the exact class of ``smoothable'' $\Q$-Fano varieties which should be added to the smooth ones in order to have a compact moduli space. We remark that the characterization of singular Fano varieties appearing in the boundary of the moduli space is unknown even in the case of Del Pezzo surfaces. The next Chapters will contain some results and conjectures about the characterizations of such degenerate Del Pezzo surfaces.
\end{itemize}

In the remaining part of this section, we describe a general heuristic recipe to relate a set of ``stable'' $\Q$-Fano varieties and GH compactifications. Here our arguments we will be of quite vague nature. A more formal discussion of these arguments will be given in Chapters $4$ and $5$, where we will apply the ideas outlined in this section to some specific examples. 

Although the YTD Conjecture relates the existence of a KE metric with the ``infinite dimensional'' notion of K-polystability,   it may happen that in many concrete situations this abstract notion of stability reduces to a  classical GIT stability notion (e.g., Chow stability for a fixed projective embedding). To fix notation and conventions, let us recall some basic definitions (\cite{M77} for more information).

Let $G$ be a complex reductive group  acting \emph{linearly} on a projective variety $X \subseteq \mathbb{P}^n$, that is through a representation of $G$ in $SL(n+1,\C)$. 
\begin{defi} Let $x\in X$ and $\hat x \in \C^{n+1}$ be a point lying over it. Then 
\begin{itemize}
\item If $0 \notin \overline{G.\hat x}$, then $x$ is called \emph{semistable} (and we write $X^{ss}$ for the $G$-invariant set of semistable points);
\item If $G.\hat x$ is closed in $\C^{n+1}$ and the stabilizer of $\hat x$ is finite, then $x$ is called \emph{stable} ($X^s$ set of stable points).
\end{itemize}
 $x$ is called \textit{unstable} if $x$ is not semistable. $x$ is called polystable if the orbits of $\hat x$ is closed (and no assumptions on the stabilizer). 
\end{defi}
The crucial result in GIT is the following \cite{M77}:
\begin{thm}(D. Mumford) Let $G \curvearrowright  X \subseteq \mathbb{P}^n$ and let $Y$ be the projective variety defined by $X//G:= Proj \bigoplus_k H^{0}(X, \mathcal{O}_X(k))^G$ (which is well-defined since the algebra of invariants is finitely generated for a reductive $G$). Then we have the following:
$$ \exists \varphi: X^{ss}\subseteq X \longrightarrow Y $$ surjective, such that 
\begin{itemize}
\item
$x,y \in X^{ss}, \, \varphi(x)=\varphi(y)\, \Leftrightarrow\, \overline{G.x} \cap \overline{G.y} \cap X^{ss} \neq \emptyset$
\item $x,y \in X^s \subseteq X^{ss}, \, \varphi(x)=\varphi(y)\, \Leftrightarrow\,    y=g.x, \, \, g\in G$
\end{itemize}

\end{thm}

This Theorem tells us that, to obtain a ``good" geometric quotient (namely the projective variety $X//G$), we have to remove the unstable points and identify some semistable orbits. Moreover, the resulting projective variety can be thought  as a natural compactification of the set of stable orbits. By the Kempf-Ness Theorem it is known that every semistable orbit has an \emph{unique} polystable representative. Hence, 
 set theoretically, the GIT quotient is in natural bijection with the polystable orbits. 

A canonical way to form moduli spaces using classical GIT is via Chow or Hilbert stabilities. We refer the reader to \cite{M77} for the precise definitions. Here we just recall that the main point consists in associating to a variety $X\subseteq \p^n$ inside a given family (which may consists of varieties of the same degree or with equal Hilbert polynomial),  a \emph{point} $[X]\in \mathcal{X}\subseteq \p^N$ (the Chow point, or m-th Hilbert point depending on the stability notion considered). Then the GIT quotient of the Chow (Hilbert) variety $\mathcal{X}$ by the induced representation of $SL(n+1)$ can be considered to be the moduli space of (polystable) varieties $X\subseteq \p^n$. A typical example is given by hypersurfaces in projective space. Here to each hypersurface we  associated a point (the Chow point) in a projective space (the Chow variety) simply by considering the point defined by the coefficients of the equation of the hypersurfaces.

Suppose now that we have a moduli of Fano varieties that is ``naturally'' associated to some classical ``ad hoc''  GIT stability notion. For example cubic surfaces in $\p^3$, intersection of two  quadrics in $\p^n$ or other situations (we will encounter several  examples in the next sections and Chapters). Thus we can apply GIT machinery to construct a compact projective moduli space of  GIT ``polystable'' Fano varieties ${\overline{\mathcal{M}}}^{ALG}$. On the other hand, we may consider the problem of existence of KE metrics on such stable Fano varieties and then its associated GH compactification ${\overline{\mathcal{M}}}^{GH}$.

The natural question is the following: when is, up to some finite quotient, ${\overline{\mathcal{M}}}^{ALG} \cong {\overline{\mathcal{M}}}^{GH}$?

Suppose that we know that a subset of ``polystable'' Fano manifolds admit KE metrics.  A possible strategy towards this identification of moduli spaces can be the following:

\begin{enumerate}
 \item Prove that GH limits of ``KE points'' in $\mathcal{M}$ are (poly)stable varieties which appear in the  algebraic compactification considered (recall that, in contrast with the GH compactification, the algebraic compactification is non-canonical).
 \item Study the local structure of the moduli space around possibly singular KE $\Q$-Fano varieties.
  \item Argue, using the previous two steps, that the identification  between   ${\overline{\mathcal{M}}}^{ALG}$  and ${\overline{\mathcal{M}}}^{GIT}$ holds.
\end{enumerate}
  
We will discuss  the second point in the next section. The third point is related to some topological argument that one can apply in order to obtain the desired identification. Essentially it is based on exploiting the Hausdorff and compactness properties of the two moduli spaces together with a kind of ``continuity method'' where the ``openness'' part is related to the study of the local structure of the KE moduli space. A formal discussion is given when we study the GH compactification of cubic surfaces in Chapter $4$.  

In the remaining part of the section, we focus our attention on the first point. First of all note that an algebraic compactified moduli space of Fano varieties  ${\overline{\mathcal{M}}}^{ALG}$ need to parameterize irreducible Fano varieties with at most log-terminal singularities in order to be a good candidate for being equal to the GH compactification. This follows by the theory of GH degenerations of KE Fano manifolds \cite{DS12}. Suppose that we find an algebraic compactification $\overline{\mathcal{M}}^{ALG}$ which parametrizes also non-normal Fano varieties. Then it may be possible to perform some birational modification of ${\overline{\mathcal{M}}}^{ALG}$ (compare the discussion of Del Pezzo of degree $2$ in Chapter $4$) to find a new compactification which parametrizes $\Q$-Fano varieties (which we can call geometric compactification). However we should stress that, in general, there can be more than one geometric compactification.

Suppose for simplicity that ${\overline{\mathcal{M}}}^{ALG}$ is realized as GIT quotient of a (component) of  the Chow variety in a \emph{fixed} projective space $\p^n$ and assume that we have a non empty set of KE polystable smooth Fano manifolds $(X_i,\omega_i)$ with ${\mathcal{O}(1)}_{X_i}=-K_{X}^k$. We may assume that $(X_i,\omega_i) \rightarrow G$ in the GH topology. Then the crucial point consists in showing that $G$ is isomorphic to a polystable  variety $X_\infty \subseteq \p^n$ (w.r.t. the same stability condition considered for constructing ${\overline{\mathcal{M}}}^{ALG}$).

In order to relate $G$ with subvarieties of $\p^n$ it is useful to consider the so called \emph{Tian's embeddings}.  We  (re)-embed the variety $X_i$ by using sections which are $L^2$ orthonormal w.r.t. the KE metric, i.e., $T(X_i)\subseteq \p^n$ (we will also refer to these embeddings as ``KE Fano manifolds in Tian's gauge''). Since these bases of sections are defined only up to the \emph{compact} group $U(n+1)$  we may take a limit inside the (compact) Chow variety, i.e., $T(X_{i_j}) \rightarrow W$ in the analytic topology of the Chow variety. We call $W$ the \emph{flat limit}. Moreover, it follows by $L^\infty$ estimates on KE orthonormal sections \cite{T90} that we have a ``map'':
$$T_{\infty}:G\setminus U_\epsilon \dashrightarrow W \subseteq \p^n,$$ 
where $U_\epsilon$ is some small neighborhood of $Sing(G)$. In the case of Del Pezzo surfaces it is known \cite{T90} that $T_{\infty}$ is actually a \emph{rational map} between the GH limit $G$ and the flat limit $W$, and it is given by $L^2$ orthonormal sections (with respect to the orbifold KE metric of the GH limit) of the sheaf $i_{*} \left(-K_{G\setminus Sing(G)}^k  \right)$ (which may be identified with the sheaf of orbifold sections by Hartog's Theorem).  

Moreover we should recall that $G$ must have the same degree as $X_i$ and the same Hilbert polynomial (since antiplurigenera are preserved under GH degenerations \cite{T90}). Similar results hold in higher dimensions (compare the very recent paper \cite{DS12}). 

The crucial point consists in showing that the above map $T_\infty$ is an isomorphism and $W$ is actually a (polystable) variety. In general the map $T_\infty$ is only rational with ``image'' that does not need to be of the same dimension as $G$. Moreover $W$ can be a priori quite nasty, e.g., not irreducible. It has been shown by G. Tian in \cite{T90} and by S. Donaldson and S. Sun in \cite{DS12} that $T_\infty$ becomes an isomorphism if we consider (re)-embeddings of $X_i$ given by \emph{sufficiently high powers} of the anticanonical bundles. These (re)-embeddings in higher dimensional projective spaces create of course severe  problems for the relations between GH degeneration and the stability notion considered. 

If  ${\overline{\mathcal{M}}}^{ALG}$ already parameterizes ``nice'' varieties (i.e., it is a geometric compactification), we may expect that (re)-embeddings in higher dimensional projective spaces are not needed. This can be verified in some situations, using the particular geometry of the Fano manifolds considered. We will discuss examples of such situations in Chapter $4$.

The above picture fits particularly well with the notion of Chow stability (in a given embedding). Recall that it is known that all \emph{smooth} Fano manifolds have uniform embeddings in fixed projective spaces \cite{KMM92}. It would be interesting to know what happens to Chow stability in these ``minimal'' embeddings (in particular which singular varieties are added to compactify the locus of smooth Fano varieties). We should also recall that Chow stability is also naturally associated to the KE problem: by the work of S. Donaldson \cite{D01}, it is known that if a KE $X$ has discrete automorphisms then it must be \emph{asymptotically} Chow stable. However we should remark that there are several examples of KE varieties (with automorphisms or with singularities) which are not asymptotically Chow stable (e.g. \cite{O10}).
Thus \emph{in general} one should be careful in believing  that Chow polystability captures the existence of a KE metric (and the GH compactification).

 It may be that at the end the abstract notion of $K$-polystability is really the correct one for the existence of KE metrics. Then one expects  the algebraic compactification ${\overline{\mathcal{M}}}^{ALG}:={\overline{\mathcal{M}}}^{K}_{h}$ to be a \emph{coarse moduli space parametrizing  $\Q$-Gorenstein smoothable K-polystable $\Q$-Fano varieties of fixed Hilbert polynomial $h$}. However it is not known that the \emph{set} ${\overline{\mathcal{M}}}^{K}_{h}$ actually carries the structure of some type of algebraic space. In particular it is not known if it is ``Hausdorff'' and ``compact'' (properties which hold for the metric compactification).

\section{Local Moduli problem}

The central question discussed in this section is the following: let $(X_0,\omega_0)$ be a (smooth) KE Fano variety and let $X_t$ be a small deformations of $X_0$ (i.e. a variation of the complex structure). Under which conditions does $X_t$ carry a KE metric? If this metric exists, is it GH close to the metric on $X_0$?

We can decompose the problem into four cases:

\begin{enumerate}
 \item $(X_0,\omega_0)$ is smooth and $Aut(X_0)$ is discrete;
 \item $(X_0,\omega_0)$ is smooth and $\mbox{dim}_{\C}Aut(X_0)>0$;
 \item $(X_0,\omega_0)$ is ``singular'' and $Aut(X_0)$ is discrete;
 \item $(X_0,\omega_0)$ is ``singular'' and $\mbox{dim}_{\C}Aut(X_0)>0$. 
\end{enumerate}

Some remarks are needed. In the next paragraphs we will discuss the known results in the case of smooth Fano manifolds (cases $1$ and $2$ above). Then we will focus our attention on the singular cases. However we will  describe only the general picture:  more precise statements are postponed to the next Chapters (e.g., $2$, $3$ and $5$).

In order to make  the  metric deformation theory of KE Fano manifolds precise,  let us recall very briefly the deformation theory of smooth Fano varieties. By the classical theory of Kodaira-Spencer-Kuranishi (compare \cite{MK06}), given a complex manifold $X$ one has that its versal deformation space is given by the analytic germ $Def(X):= ob_X^{-1}(0)$, where $ob_X$ denotes the obstruction holomorphic map:
$$ob_X: \T^1_X:=H^1(X,\Theta_{X}) \longrightarrow \T^2_X:=H^2(X,\Theta_{X}).$$
Here $\Theta_{X}$ denotes the sheaf of holomorphic vector fields on $X$.
If $X$ is a Fano manifold, we have the following.

\begin{lem} Let $X^n$ be a n-dimensional Fano manifold. Then
 $$\T^2_X=0.$$
In particular its versal deformation space $Def(X)$ can be identified with an open neighborhood of the origin in $\T^1_X$. 

Moreover every $X_t$ with $t\in Def(X)$ sufficiently small is itself an $n$-dimensional Fano manifold. 
\end{lem} \label{DF}
\begin{dimo}
 By Serre duality $\T^2_X=H^{2}(\Theta_X)=H^{n-2}(\Theta_{X}^\checkmark \otimes K_X)^\checkmark$. Moreover by Kodaira-Nakano vanishing (e.g.,\cite{MK06} Chapter $7$)
$H^{n-2}(\Theta_{X}^\checkmark \otimes K_X)=H^{n-2}(\Omega^{1}(K_X))=0$ since $K_X$ is negative.

In order to prove that small deformations of $X$ are Fano, let us recall that we may assume that $X$ is isomorphic to a fiber $X_0$ of a fibration $$\pi: \mathcal{X} \rightarrow Def(X),$$
i.e. $X \cong \pi^{-1} (0)=:X_0$, $\mathcal{X}$ is a smooth complex manifold and $\pi$ an holomorphic submersion.  Let $-K_{\mathcal{X}/Def(X)}$ be the relative anticanonical bundle. By hypothesis ${-K_{\mathcal{X}/Def(X)}}_{|{X_0}}\cong -K_{X_0} $ is ample. By  \cite{KM98}, Proposition 1.41, it follows that ampleness is an open condition in a family, i.e., for $t$ small enough $X_t$ has ample anticanonical bundle.
\qed
\end{dimo}

Recall that if $Aut(X)$ is discrete then $Def(X)$ is  actually a universal deformation. Hence $Def(X)$ may be identified, up to some finite symmetries, with a neighborhood of the moduli space.

With the above observation and Lemma \ref{DF} in mind, it makes sense to define
$$\mbox{Virdim}\,\mathcal{M}_{[X]}:= h^{1}(\Theta_X)-\mbox{dim}_{\C} \, Aut(X),$$
to be the the \emph{virtual dimension} of the moduli space around $X$.

However recall that the action of the non-compact group $Aut(X)$ on the versal  deformation space $Def(X)$ can be quite nasty. In particular the space $Def(X)/Aut(X)$ equipped with the quotient topology won't  in general be Hausdorff.

Since smooth Fano threefolds are classified (by V. Iskovskikh, S. Mukai,   etc\dots compare \cite{IP99}), it is interesting to compute the virtual dimension of their moduli spaces as a  function of (some of) the principal invariants of the classification (see the Appendix of \cite{IP99} for a complete list of  Fano 3-folds in terms of their principal invariants).

\begin{prop}\label{VDM} Let $X$ be a smooth $n$-dimensional Fano manifold. Then
\begin{itemize}
 \item for $n=2$, $\mbox{Virdim}\,\mathcal{M}_{[X]}=10-2\,deg(X)$;
 \item for $n=3$, $\mbox{Virdim}\,\mathcal{M}_{[X]}= 18-\rho(X)+\frac{1}{2}\left(b_3(X)-deg(X)\right)$;
\item in general, $\mbox{Virdim}\,\mathcal{M}_{[X]}=-td(X).ch(X)$.
\end{itemize}
Here $deg(X)=(-K_X)^n$,  $\rho(X)=b_2(X)$ (the Picard rank), $td(X)$  the Todd class and $ch(X)$ the Chern character.
\end{prop}

\begin{dimo}
 By Serre duality and the  Kodaira-Nakano vanishing Theorem $h^{i}(\Theta_X)=h^{n-i}(\Theta_{X}^\checkmark \otimes K_X)=h^{n-i}(\Omega^{1}(K_X))=0$ if $i\geq2$ since $-K_X$ is ample.
Thus it follows by Hirzebruch-Riemann-Roch that $-\mbox{Virdim}\,\mathcal{M}_{[X]}=\chi(\Theta_X)=td(X).ch(X).$
Recall that $td(X)=1+\frac{c_1}{2}+\frac{c_1^2+c_2}{12}+\frac{c_1c_2}{24}+\dots$ and $ch(X)=n+c_1+\frac{c_1^2}{2}-c_2+\frac{c_1^3}{6}-\frac{c_1 c_2}{2}+\frac{c_3}{2}+\dots,$ where $c_i$ denotes the ith-Chern class of the tangent bundle.   A computation shows that

\begin{itemize}
 \item $n=2$: $\mbox{Virdim}\,\mathcal{M}_{[X]}=\frac{5c_2-7c_1^2}{6}$. Since $-K_X$ is ample, it follows by Kodaira vanishing that $1=\chi(\mathcal{O}_X)=td(X)$, that is $c_2=12-c_1^2$. Thus $\mbox{Virdim}\,\mathcal{M}_{[X]}=10-2c_1^2$ as desired.
 \item $n=3$: $\mbox{Virdim}\,\mathcal{M}_{[X]}=-\frac{c_3}{2}-\frac{c_1^3}{2}+\frac{19}{24}c_1c_2$. By Kodaira vanishing, $c_1c_2=24$.  Using $c_3=e(X)=2+2\rho-b_3$, we find $\mbox{Virdim}\,\mathcal{M}_{[X]}=18-\rho(X) +\frac{1}{2}(b_3(X)-c_1^2)$.
\end{itemize}

\qed
\end{dimo}

Of course, in the case of Del Pezzo surfaces the above formula for the dimension of the moduli space agrees with the count given by looking at the ways one can realize a Del Pezzo surface as a blow up of points in $\p^2$. Browsing the list of Fano $3$-folds, it emerges that the virtual dimension is positive in many cases. For example in the case of Mukai-Umemura 3-folds we observe that the virtual dimension of their moduli space is $6$ (later in this section we will discuss a possible description for a  $6$ dimensional compactified moduli space of ``stable'' (KE?)  Mukai-Umemura varieties).
We may also point out that if the degree of a Fano 3-fold is sufficiently big (i.e., bigger than $38$) the virtual dimension is negative (i.e., there are holomorphic vector fields). Around $1/3$ of the deformation classes of Fano $3$-folds are in this range.

\subsection*{Case $1$}

For smooth Fano manifolds with no continuous families of automorphisms, the study of deformations of KE metrics is well known:

\begin{thm} \label{D1} Let $(X_0,\omega_0)$ be a smooth KE Fano manifold with  $Aut(X_0)$  discrete. Then for all $t\in Def(X)$ sufficiently small, $X_t$ admits a (unique) KE metric $\omega_t$ which is GH close to $\omega_0$. 
\end{thm}

\begin{dimo} We want to solve the equation $$E(t,\varphi):=(\beta_t + \D_t \varphi)^n- e^{f_t-\varphi} \beta_t^n=0,$$
for $t$ small. Here $\beta_t$ denotes a background K\"ahler form in $c_1(X_t)$ and $f_t$ its Ricci potential, i.e., $Ric(\beta_t)=\beta_t+ \Dt  f_t$. 

W.l.o.g. we may assume that  we have fixed the background smooth structure, say $M$, and $X_t$ is a K\"ahler manifold for a smooth family of  integrable complex structures $J_t$ where $t$ is now a real parameter and $\beta_0=\omega_0$ (and $f_0=0$).

Then $$\partial_\varphi E_{(0,0)}=\Delta_{\omega_0}+1,$$
with kernel isomorphic to the Lie algebra of holomorphic vector fields on $X_0$, via $\phi \rightarrow (\overline{\partial} \phi)^{\sharp_{\omega_0}}$ (compare \cite{T97} or the proof of Proposition \ref{ker}). By our hypothesis we know that this linear operator has no kernel. Note also that it is formally self-adjoint.

Thus realizing the operator $E(t,\varphi)$ in suitable Banach spaces (e.g., $C^{2,\gamma}_{\beta_t}$), it follows by the implicit function Theorem that $E(t,\varphi(t))=0$  for some smooth $\varphi_t$ and small $t$. That is $\omega_t:= \beta_t + \Dt \varphi(t)$ is the desired KE form.

In order to see that $\omega_t$ is GH close to $\omega_0$, it is sufficient to observe that, on the level of smooth  metric tensors on $M$, $||g_0-g_t||_{g_0} \rightarrow 0 $ as $t \rightarrow 0$ and apply Lemma \ref{GHCC}.

\qed 
\end{dimo}

For example the above Theorem can be used to prove that small deformations of KE Del Pezzo surfaces of degree $d\leq 4$ admit KE metrics. In particular, up to a finite group action, we may identify $Def(X)$ with a neighborhood of $[X_0,\omega_0]$  in the space of KE metrics.

\subsection*{Case $2$}
In the  presence of a continuous family of automorphisms the situation is more subtle. In fact, as we can see by looking into the proof of  Theorem \ref{D1}, non-trivial holomorphic vector fields give obstructions  to the invertibility of the linearized equation. However, thanks to the work of G. Sz\'ekelyhidi \cite{S10} (and the PhD Thesis of T. Br\"onnle  \cite{TB11}), also in this case the behavior of the KE metric under small deformations is completely understood.   Sz\'ekelyhidi's Theorem holds more generally for K-polystable deformations of cscK metric. In the Fano framework, his Theorem assumes the following form:

\begin{thm}[G. Sz\'ekelyhidi]\label{D2} Let $(X_0,\omega_0)$ be a smooth KE Fano manifold with  $dim_{\C} Aut(X_0)>0$. Then $X_t$ for $t\in Def(X_0)$ sufficiently small admits a (unique up to automorphisms) KE metric $\omega_t$ if and only if $t$ is polystable (i.e., the orbit is closed) with respect to the linear action  of $Aut(X_0)$  on the Kuranishi space $\T^1_{X_0}$. 

Moreover the KE metrics on the polystable deformations are close to the metric $\omega_0$ in the GH sense.
\end{thm}

\begin{dimo} The proof of Sz\'ekelyhidi uses the description of the cscK equation as a moment map for the action of (exact) symplectomorphisms on the infinite dimensional space of complex structures $\mathcal{J}_\beta$ compatible with a fixed symplectic form $\beta$. We outline very briefly the main steps in the proof and we explain why in the Fano case the Theorem has the form we stated.

 The crucial point in his argument is to show that there is a finite dimensional reduction of the problem, i.e., an equivariant map for the action of complexification of  the group of Hamiltonian isometries $K^{\C}$ of $(X_0,\omega_0=\beta)$
$$ \Phi: B \subseteq \T^1_X \longrightarrow  \mathcal{J}_\beta.$$
 Then one shows that if $v\in B$,  sufficiently small, is polystable for the linear action of $K^{\C} $,  then one can deform $v$ in its $K^{\C}$ orbit to $\tilde{v}$ with $\mu(\tilde{v})=0$, where $\mu$ is the cscK equation (moment map). 

If $(X_0,\omega_0)$ is a smooth KE Fano manifold, then we may fix $\omega_0$ as background fixed symplectic form and look at the complex structures compatible with it. This is not a restriction, since every small Fano deformation of $X_0$ has $c_1(X_t)=c_1(X_0) \in H^2(M,\Z)$ (where $M$ is the underlying smooth manifold). Then by  the general theory of deformations previously explained, we may assume that $B$ is contained in the space of versal deformations $Def(X_0)$. 

Since $X_0$, being Fano,  is simply-connected (compare, for example, \cite{B06} Theorem 6.12), every isometry (in the connected component of the identity) is Hamiltonian. Moreover it follows by Matsushima's theorem that   the complexification of the isometry group is the group of holomorphic automorphisms. Finally, it is sufficient to recall that the cscK equation reduces to the Einstein equation when we are considering (anticanonical) polarizations.

In order to see that the condition of being polystable is not only sufficient for the existence of an Einstein metric but also necessary, observe that if $X_t$ is not polystable we can find a $1$-ps (parameter subgroup) degeneration to a polystable (hence KE) Fano $X_s$ for $s\in Def(X_0)$. In particular it follows by a result of G. Tian on analytic stability of smooth Fano manifolds  \cite{T97} that $X_s$ cannot carry a KE metric. For another point of view, observe that the Donaldson-Futaki invariant of the considered $1$-ps  is $0$, but the induced test configuration is not trivial, i.e., $X_s$ is not $K$-polystable.

The fact that polystable $(X_t,\omega_t)$ and $(X_0,\omega_0)$ are close in the GH sense follows by the way  the  KE metric $\omega_t$ is obtained from an implicit function Theorem argument.

\qed
\end{dimo}

As a consequence we may identify, up to a finite quotient, a neighborhood of $[X_0,\omega_0]$  in the space of KE metrics with a neighborhood in the GIT quotient $Def(X_0)// Aut_0(X_0)$.

A well-known example where the above Theorem \ref{D2} applies is the case of the Mukai-Umemura Fano threefold $U_{22}$ (compare \cite{T97} and \cite{D08}). $U_{22}$ admits a KE metric and has  $PSL(2,\C)$ as automorphism group. Then $Def(U_{22})$ can be identified as an open neighborhood of the origin in the $9$-dimensional space $\T^1_{U_{22}}$. Some of the orbits of the action of $PSL(2,\C)$ on $\T^1_{U_{22}}$ are not closed (hence the corresponding ``Iskovskikh threefolds'', i.e., smooth Fano threefolds of degree $22$ and Picard rank one, do not admit  KE metrics). On the other hand, KE metrics exist  on all small deformations corresponding to the  closed orbits.   

As a small digression, we can mention that the way Iskovskikh threefolds are constructed, is naturally associated with a classical GIT framework.  Recall that smooth Iskovskikh threefolds are given by the zero set of three generic sections $s_i$ of $\Lambda^2 U^\checkmark \rightarrow Gr_3(V)$, where $U$ is the tautological vector bundle over the Grassmannian. Then the classical six dimensional GIT quotient
$$\overline{\mathcal{M}_{I}}^{ALG}:= Gr_3(\Lambda^2 V^\checkmark) // SL(V),$$
where $V$ is a $7$-dim vector space and the linearization of $SL(V)$  is induced by Pl\"ucker embedding, can be a good candidate for being the GH compactified set of Iskovskikh threefolds. Thus
\begin{quest}
 Does $\overline{\mathcal{M}_{I}}^{ALG}$ parameterizes only $\Q$-Fano varieties (``degenerate Iskovskikh threefolds'')? If it is so, is $\overline{\mathcal{M}_{I}}^{ALG}$ homeomorphic to the GH compactification of the set of smooth KE Iskovskikh manifolds (which is not empty, since it contains  Mukai-Umemura threefold $U_{22}$ and some of its deformations)?
\end{quest}

\subsection*{Cases $3$ and $4$}

The above two cases $1$ and $2$  explain completely what happens  around smooth KE Fano manifolds from the metric view point. However, to study the GH compactification,  it is also important to look at the deformation theory of singular KE spaces. Before starting to talk about deformation theory, it is important to say something on what  a ``singular'' KE Fano manifold is.

The underlying complex variety must be a $\Q$-Fano variety $X$ (i.e., irreducible, normal, with log-terminal singularities and ample $\Q$-Cartier anticanonical divisor). A KE metric $\omega$ on $X$ must be, first of all, a smooth (incomplete) KE metric on the smooth part $X \setminus Sing(X)$. The crucial problem consists in describing what happens at the singularities. A first basic property which is desired,  is that the metric completion of $(X \setminus Sing(X), \omega)$ is topologically equal to the variety $X$ itself.

In the case of quotient singularities, the natural class of metrics to consider is of course the class of orbifold metrics. However recall that, by virtue of the results recalled in the previous section and Schlessinger's rigidity Theorem on isolated quotient singularities of dimension bigger then or equal to three \cite{S71}, the KE Fano orbifolds which are expected to appear in GH compactifications should have singularities locally of the form $$\C^{n-2}\times \C^2/ \Gamma,$$ with $\Gamma \subseteq U(2)$ finite acting freely away from the origin. For example, in the study of KE compactified moduli spaces of Fano 3-folds, we can avoid the study of rigid, in the sense of smoothability, KE Fano varieties with isolated singularities of the form $\C^3 /  \Gamma$, with $\Gamma \in U(3)$ acting freely away from the origin.
However we can not avoid  Fano manifolds which are singular along complex curves. In the case the singular locus is a \emph{smooth} complex curve, then the KE metric to consider is of orbifold type (we will show an example of this kind of singular KE Fano variety at the end of Chapter $5$).
  
On the other hand it is also unavoidable to consider singularities which are not of orbifold type. The simplest singularity of this kind is the node $\{z_1^2+\dots +z_n^2=0\} \subseteq \C^n$. In general the existence theory of KE metrics on these singular spaces is quite complicated. In \cite{BBEGZ11} the authors consider ``weak'' solutions to the KE equation and give sufficient condition for the existence of such singular KE metrics.  However nothing can be said at present about the ``asymptotic'' behavior of the metric near the singularities. As we will explain later in Chapter $5$, we expect that the KE metric has some asymptotic rate of decay to some Calabi-Yau (CY) cone metric at the singularities. To the author knowledge, no examples where  decays  can be computed are known. This lack of information makes the study of the KE local moduli space at this singular space hard to tackle. 

As we will carefully explain in Chapters $2,3$ and $5$, a way to study the deformation theory of KE metrics at singular KE Fano varieties is via ``gluing techniques''. The rough idea (at least for isolated singularities) is to reverse the GH degeneration picture by ``gluing'' an asymptotically conical (AC)  CY metric on a small deformation of the singularity and ``deforming'' the approximate solution to a genuine KE metric. Thus a crucial input in the construction is given by CY metrics on smoothings of the singularities which satisfy some curvature decay at infinity. In the case of two dimensional orbifold singularities one can consider, for example,   Kronheimer's ALE metrics \cite{K89}. In higher dimension some examples of these AC CY metrics are known (e.g., Stenzel's metric on the smoothing of a node \cite{S93} and the examples in the recent paper \cite{CH12}).

In chapter $2$ we will describe the expected picture in the case of Del Pezzo orbifolds. In Chapter $3$ we will study the deformations of nodal Del Pezzo surfaces, verifying a case of the Conjecture made in Chapter $2$. In the last Chapter we will briefly discuss deformations of nodal Fano varieties.

\end{chapter}
\begin{chapter}{KE metrics on (singular) Del Pezzo surfaces and deformations}

We begin this Chapter by briefly discussing the problem of existence of KE (orbifold) metrics on singular Del Pezzo surfaces. In particular, we describe a relation between the degree of a singular KE Del Pezzo surface and the (maximal) order of the orbifold groups. Then we show an easy example of a countable family of topologically distinct singular  KE Del Pezzo surfaces which we prove to be not asymptotically Chow semistable. These observations, which we believe should be well-known to the experts,  illuminate the difference between the case of smooth and singular Del Pezzo surfaces. In particular, it is clear that without establishing a precise link with some stability (e.g., K-polystability), the problem of classifying singular KE Del Pezzo surfaces seems to be hopeless.   

The subsequent two sections are devoted to the study of some special classes of singular KE Del Pezzo surfaces, i.e., Del Pezzo surfaces which admit smoothings. This investigation is of course crucial in understanding how the KE compactified moduli space looks at its boundary points. After reviewing the deformation theory of complex spaces, with emphasis on the case of log-terminal Del Pezzo surfaces, we make a Conjecture on the expected behavior of the KE metrics under $\Q$-Gorenstein smoothings. Finally we show an example of two KE Del Pezzo surfaces,  of degree one and two, which admit $\Q$-Gorenstein smoothings and have genuine  log-terminal singularities (i.e., non-canonical). In view of a Conjecture of Tian (\cite{T06}, pag. 387), these examples look quite unexpected.  As a byproduct of the construction of these examples, we get a self-contained easy proof of the classification of $\Q$-Gorenstein smoothable toric KE Del Pezzo surfaces of Picard rank one.

\begin{rmk} In this Chapter we use several times some results from Toric Geometry. In Appendix A we have collected the principal ``Toric Propositions'' we used.    
\end{rmk}
\section{The landscape of KE log-terminal Del Pezzo surfaces: an overview}

It is classically very well-known that there are only $10$ topological types of \textit{smooth} Del Pezzo surfaces, namely $\p^2$, $\p^1 \times \p^1$, and the blowup of $\p^2$ in strictly less then $9$ points (in very general position). However, it is important from the viewpoint of both algebraic and differential geometry to enlarge the category of smooth Del Pezzo surfaces by including some kind of singular varieties. A natural class to consider is given by the collection of the so called \textit{log-terminal Del Pezzo surfaces}.
\begin{defi} A normal complex two dimensional variety $X$ is called a log-terminal Del Pezzo surface if $-K_X$ is an ample $\Q$-Cartier divisor and the singularities of $X$ are at most log-terminal, i.e., $\mbox{Discrep}(X) > -1$.
\end{defi}

Remarkably, in complex dimension two the concept of log-terminal singularity coincide with the concept of quotient singularity (orbifold), i.e., singularity of the form $\C^2/ \Gamma$ with $\Gamma$ finite subgroup of $U(2)$ acting freely on $\C^2 \setminus \{0\}$. The important class of canonical singularities ($\mbox{Discrep}(p)\geq 0$), also known as ADE singularities,  is given by quotients of $\C^2$ by a finite subgroup of $SU(2)$.

In contrast to the case of smooth Del Pezzo surfaces, the number of topological types of log-terminal Del Pezzo surfaces is unbounded, i.e., there are infinitely many non homeomorphic singular Del Pezzo surfaces. An easy example of this phenomenon is given by considering the log-terminal Del Pezzo $X_n$ defined by blowing down the $(-n)$-curve in the Hirzebruch surface $\mathbb{F}_n:=\p\left(\mathcal{O}\oplus\mathcal{O}(-n) \right)$, i.e., $X_n$ is equal to the weighted projective plane  $\p(1,1,n)$. Observe that $\mbox{Discrep}(X_n)=-1+\frac{2}{n}>-1$ and $\mbox{deg}(X_n)=K_{X_n}.K_{X_n}= \frac{(n+2)^2}{n} \rightarrow +\infty$ (recall that for smooth Del Pezzo surfaces the anticanonical degree is always an integer between $1$ and $9$). 

From a differential geometric perspective it is natural to ask which log-terminal Del Pezzo surfaces admit orbifold K\"ahler-Einstein metrics. The answer is completely known for smooth Del Pezzo \cite{T90} surfaces, but it remains rather mysterious for the singular ones: it is known that some of them  admit KE metrics (by $\alpha$-invariant computation, compare for example \cite{C08}, \cite{C09}) and that some do not (by Futaki invariant computation \cite{DT92}). However an algebraic characterization is still missing. Of course, it is quite natural to believe that there should exist a notion of stability (e.g., K-polystability) capturing the existence of KE orbifold metrics.

The technique used by Tian to study the existence of KE metrics in the smooth case is based on a continuity method through  deformation of the complex structure on a topologically fixed smooth Del Pezzo. This technique has two  disadvantages. Firstly, it requires a starting point, that is a KE metric for some special complex structure. Secondly,  a case-by-case analysis using numerical properties of the anticanonical divisors  is needed in order to prove closedness. These observations, together with the previous discussion on the topological type of log-terminal Del Pezzo surfaces,  suggest the difficulty to extend Tian's approach to orbifold Del Pezzo surfaces.

 Finally it is worth noting that a new proof of the existence of KE metric on smooth Del Pezzo has been found using the K\"ahler-Ricci Flow \cite{W10}. This method has the advantage of not requiring a starting KE metric and of generalizing partially in a natural way to the orbifold setting. However a case-by-case analysis (again based on refined $\alpha$-invariant computations) is still needed to prove the convergence of the flow.   

In the next paragraphs, we describe some properties of singular KE Del Pezzo surfaces.

We begin by showing that there are restrictions on the anticanonical degree of KE Del Pezzo surfaces. Similar results were  also obtained by Tian for GH degenerations of smooth KE Del Pezzo surfaces (compare Theorem $5.2$ in \cite{T90}). His estimate has the same shape as our estimate given in the below Proposition. However, his estimate contains a $48$ instead of our $12$. We believe that $12$ is the correct number (see the proof below). Having this new control on the degree will be important in the discussion of Chapter $4$ concerning GH degenerations of Del Pezzo surfaces. 
\begin{thm} \label{NDP}Let $X$ be a KE log-terminal Del Pezzo surface. Define $|\Gamma_{\max}|:=\max_{p \in Sing(X)} |\Gamma_p|$, then 
$$|\Gamma_{\max}|\mbox{deg}(X) < 12.$$
\end{thm}

\begin{dimo} W.l.o.g. we can assume that the metric satisfies $Ric (g)= 3 g$, so that, by the orbifold version of Myers' Theorem \cite{B93},  $\mbox{diam}(X) \leq \pi$.  Using the Bishop-Gromov volume comparison for orbifolds \cite{B93}, for all $p \in X$ the function
$$r \mapsto \frac{Vol(B(p,r))}{Vol(\overline{B}(r))},$$
is non-increasing with limit as $r \rightarrow 0$ equal to $\frac{1}{|\Gamma_p|}$, where $\overline{B}(r)$ denotes the ball of radius r in the standard sphere of radius one $S^4(1)$. 

Combining the diameter estimate with the volume comparison, we find that
$$Vol(X) |\Gamma_p| \leq Vol (S^4(1)).$$
 The condition $Ric (g)= 3 g$ implies that the associated K\"ahler form $\omega$ satisfies $[w] = \frac{2 \pi}{3} c_1^{orb}(X)$.
 Then
 $$\begin{array}{ccl}
Vol(X) &=&\int_X \frac{\omega^2}{2}\\
&=&\frac{2 \pi^2}{9} \int_X c_1^2\\
&=& \frac{2 \pi^2}{9} \mbox{deg}(X).
\end{array}
$$
Recalling that $Vol(S^4(1))=\frac{8}{3} \pi^2$ and taking $p$ equal to the point where the order of the orbifold group is maximal, we find the desired estimate with ``$\leq$'' instead of ``$<$''. To remove the equality case, it is sufficient to note that if the equality is achieved, then $X$ is a finite quotient of the sphere (compare \cite{B93}). But this is impossible, since $X$ is K\"ahler.  \qed \\

\end{dimo}
Since for singular log-terminal Del Pezzo surfaces $|\Gamma_{max}| \geq 2$, we immediately find that the degree must be bounded from above. More precisely we have:

\begin{cor} Let $X$ be a \emph{singular} KE log-terminal Del Pezzo surface. Then
$$\mbox{deg}(X) < 6.$$
\end{cor}

As an easy application of the above corollary we observe that the surfaces $\p(1,1,n)$ do not admit KE orbifold metrics, since $\mbox{deg}(X_n)\geq 8$ for $n\geq 2$ (compare also \cite{GMSY07} and \cite{RT11}).

Another consequence is the following restriction on the type of singularities of  KE Del Pezzo surfaces with canonical singularities (see Appendix B):

\begin{cor}\label{KECDP}
 Let $(X,\omega)$ be a KE Del Pezzo surface with canonical singularities. Then, denoting with $d:=(-K_X)^2$ the degree,
\begin{itemize}
 \item $d=4$: Sing($X$) consists (at most) of points of type $A_1$;
 \item $d=3$: Sing($X$) consists (at most) of point of type $A_1$ and $A_2$;
 \item $d=2$: Sing($X$) consists (at most) of point of type $A_1$, $A_2$, $A_3$ and $A_4$;
 \item $d=1$: Sing($X$) consists (at most) of point of type $A_1$, $A_2$, $A_3$, $A_4$, $A_5$, $A_6$, $A_7$, $A_8$  and $D_4$.
\end{itemize}
\end{cor}
Note that in degree one case singularities of type $A_9$, $A_{10}$ and $A_{11}$ cannot exist by obvious Picard rank considerations. The above Corollary should be compared with Conjecture 1.19 in \cite{CK10}. However, we will show later that the above restriction on the singularity type is only a necessary condition. In fact, there exist Del Pezzo surfaces of degree four with only one nodal singularity which do not admit (orbifold) KE metric (see Corollary \ref{NKE}). The natural Conjecture about existence of KE metrics on Del Pezzo with canonical singularities is the following: a Del Pezzo with canonical singularities admits a KE if and only if is K-polystable.

Next, it is natural to ask if the set of KE Del Pezzo orbifolds consists of finitely many distinct topological types. However it is easy to construct a sequence of distinct KE Del Pezzo orbifolds with $deg(Y_n)\rightarrow 0$ (which is possible since the anticanonical divisor is only $\Q$-Cartier).

\begin{prop} \label{Yn} Let $(Y_n)_{n\in \N}$ be the following family of toric surfaces defined by the complete fan whose minimal generators of the rays are
$$v_1:=(-n,1), v_2=(n,1), v_3=(n,-1),v_4=(-n,-1)$$ 
Then $Y_n$ are  \emph{K\"ahler-Einstein} log-terminal Del Pezzo surfaces satisfying
\begin{itemize}
\item $\mbox{Sing}(Y_n)=2\C/\Z_{2n}+2A_{2n-1}$, $\mbox{deg}(Y_n)=\frac{4}{n}$ and $\mbox{Discrep}(Y_n)=-1+\frac{1}{n}$; 
\end{itemize}
Here with $\C/\Z_{2n}$ we denote the singularity given by the cone over the  Veronese embedding of $\p^1$ of degree $2n$.

In particular, the number of topological types of KE log-terminal Del Pezzo surfaces is unbounded.
\end{prop}

\begin{dimo}
 The proof follows by elementary computations in toric geometry together with the fact that for toric varieties the existence of a KE metric is equivalent to the vanishing of the Futaki invariant \cite{SZ11}, i.e., $\mbox{Bar}(P^\checkmark)=0$ where $P^\checkmark$ denotes the moment polytope (see Appendix A).
\qed \\
\end{dimo}

The final observation about singular KE Del Pezzo surfaces is related to ``stability''. It follows by  Mumford's criterion for instability of singular varieties \cite{M77}  that the above KE varieties $Y_n$ are not asymptotic Chow semistable. Similar examples, in the CY and negative first Chern class case, have been considered in \cite{O10}.

\begin{cor}Let $Y_n$ be the singular KE Del Pezzo surface defined in Proposition \ref{Yn}. Then $Y_n$ is not asymptotically Chow semistable for $n\geq 4$.
\end{cor}

\begin{dimo}

 Mumford's instability criterion states that if a variety $X$ is asymptotic Chow semistable then $$\mbox{Mult(p)} \leq (\mbox{dim}_{\C}(X)+1)!,$$
for all $p \in X$, where $\mbox{Mult}(p)$ (the multiplicity) denotes the coefficient of maximal degree of the Hilbert-Samuel polynomial, i.e.,
$$\mbox{dim}_{\C}\left(\frac{\mathcal{O}_{X,p}} {m_p^k}\right)= \frac{\mbox{Mult(p)}}{n!} k^n+\mathcal{O}(k^{n-1}),$$
where $n=\mbox{dim}(X)$ and $m_p$ is the maximal ideal at $p$.

 The value of the  multiplicities of the singularities we are considering is well-known \cite{A66}. For sake of completeness, we provide an elementary proof. 
The singularities of $Y_n$ are given by two canonical singularities of type $A_k$ and two cones over Veronese  embeddings. Since canonical singularities in dimension two are defined by a single polynomial in $\C^3$ of multiplicity two at the origin, $\mbox{Mult}(p)=2$ for $p$  an $A_k$ point.
To compute the multiplicity of the cone $V_d$ over the degree d Veronese embedding of $\p^1$, it is sufficient to note that $$\frac{m_p^j}{m_p^{j+1}} = H^0\left(\p^1, \mathcal{O}(jd) \right) \;\mbox{ and} \; \frac{\mathcal{O}_{V_d,p}} {{m_p}^{j+1}}=\frac{\mathcal{O}_{V_d,p}} {m_p^j} \oplus \frac{m_p^j}{{m_p}^{j+1}}.$$ Then a simple computation shows that
$$\mbox{dim}_{\C}\left(\frac{\mathcal{O}_{V_d,p}} {m_p^k}\right)= \frac{d}{2} k^2+\mathcal{O}(k),$$
that is $\mbox{Mult(p)}=d$. 

Applying the above computations to our surfaces $Y_n$, we find that $Y_n$ has singularities of multiplicity equal to $2n$. Hence  Mumford's instability criterion is satisfied as long as $n\geq4$. 

\qed \\

\end{dimo}

In the next sections of the Chapter we are going to study the  class of singular KE Del Pezzo surfaces which admit a smoothing. This class is of course particularly important for the study of the metric compactification of the moduli of smooth KE Del Pezzo surfaces.

\section{Deformations of Del Pezzo orbifolds}

Let $(X_0,\omega_0) $ be a KE  Del Pezzo surface with log-terminal singularities. It is then natural to ask when (partial) smoothings $X_t$ of $X_0$ carry  KE (orbifold) metrics $\omega_t$ which are close in the GH sense to the metric $\omega_0$.
As we have already pointed out in the first Chapter, the investigation of the metric behavior under small deformations (smoothings) of the complex structure is crucial for understanding  the local structure of the GH compactification at the ``boundary''.  

Before stating a conjectural answer to the above question, we need to collect  some facts about deformations of Del Pezzo orbifolds. The material should be well-known, but is sometimes quite difficult to find in the literature. Our main references are \cite{G74}, \cite{P76}, \cite{M95}, \cite{HP10}, \cite{C11}.

By the work of Grauert on the general theory of deformations of complex spaces \cite{G74}, it is known that for every $X$ there exists a \emph{versal} (or semiuniversal) deformation, i.e., a deformation (flat family) $\mathcal{U} \rightarrow Def(X)$ which has the property that every other deformation  $\mathcal{X} \rightarrow S$ of $X$ is obtained from it by pulling back by a (not necessary unique) map $S \rightarrow Def(X)$.  
Recall that if $X$ is a smooth variety, $Def(X) \subseteq \T^1_X \cong H^1(X,\Theta_X)$ is a germ of an analytic space with equations given by an obstruction map $ob: \T^1_X \rightarrow \T^2_X= H^2(X,\Theta_X)$. In particular if $\T^2_X=0$ (as in the case of Fano manifolds) then $Def(X)$ is an open set of $\T^1_X$. 

To study the deformation theory of a singular $X$, it is important to consider the following exact sequence  (the $5$-term exact sequence associated to the Grothendieck local-to-global spectral sequence),
$$0 \rightarrow H^1(X,T^0_X) \rightarrow \T^1_X \rightarrow H^0(X,T^1_X)\rightarrow H^2(X,T^0_X) \rightarrow \T^2_X
$$
where, for notational simplicity, we have defined $\T^i_X := Ext^{i}\left(\Omega_X^1, \mathcal{O}_X\right)$ and
$T^i_X:= \mathcal{E}xt^{i}(\Omega_X^1, \mathcal{O}_X).$  Here $\Omega^1_X$ denotes the sheaf of K\"ahler differentials.

The following result is well-known \cite{HP10}. For reader convenience we sketch a simple proof.

\begin{lem} Let $X$ be a Del Pezzo surface with log-terminal singularities. Then
$$H^{2}\left(X,T^0_X\right)=0.$$
 \qed
\end{lem}

\begin{dimo} By hypothesis $X$ admits at most isolated quotient singularities. By Hartog's Theorem, the sheaf  $i_{*} \left(\mathcal{H}om(\Omega_{X}^1,\mathcal{O}_{X_0})\right)_{|X \setminus \mbox{Sing}(X)}$, where $i: X\setminus \mbox{Sing}(X) \hookrightarrow X$ denotes the inclusion,  is equal to the sheaf $\Theta_{X}^{orb}$ of orbifolds holomorphic vector fields.
 Since $\mbox{depth}( \mathcal{H}om(\Omega_{X}^1,\mathcal{O}_{X}))=2$,  we have $\mathcal{H}om(\Omega_{X}^1,\mathcal{O}_{X})=i_{*} \left(\mathcal{H}om(\Omega_{X}^1,\mathcal{O}_{X_0})\right)_{|X \setminus \mbox{Sing}(X)}$. Thus $h^i\left( \mathcal{H}om(\Omega_{X}^1,\mathcal{O}_{X})\right)=h^{i}(\Theta_X^{orb})$ if $i\geq 0$.
 
Now the proof in the orbifold case is identical to the prove in the smooth case. By (the orbifold versions of) Serre duality and Kodaira-Nakano vanishing Theorem, $h^2\left(\mathcal{H}om(\Omega_{X}^1,\mathcal{O}_{X})\right)=h^{2}(\Theta_X^{orb})=h^{0}(\Omega^{1}_{orb}(K_{X}^{orb}))=0$ since $K_{X}^{orb}$ is negative.
\qed
\end{dimo}

By the above Lemma it follows that the space of infinitesimal deformations $\T^1_X$ of a Del Pezzo orbifold splits (non-canonically) as 
$$\T^1_X\cong H^{1}(T^0_X)\oplus H^{0}(T^1_X)\cong H^{1}(\Theta_X^{orb})\oplus \bigoplus_{p\in Sing(X)}\C^{n(p)},$$
where $n(p)$ denotes the dimension as a $\C$-vector space of the stalk ${T^{1}_{X}}_{p}$ (note that on the smooth locus $U=X \setminus \mbox{Sing}(X)$,  $(T^1_X)_{| U}=0$ since $\Omega^1_X$ is locally free).

One could think that, in the case of Del Pezzo orbifolds, $Def(X)$ is an open subset of $\T^1_X$, since $H^{2}\left(X,T^0_X\right)=0$. However the situation is a bit more subtle (obstructions are in $\T^2_X$). 

Let $V_p$ be the germ of a singularity of $X$ and let $Def(V_p)$ be  its versal deformations space. Then one can consider the natural map 
$$\phi:Def(X) \rightarrow \times_{p\in \mbox{Sing}(X)} Def(V_p).$$
The condition on the vanishing of $H^2(X,T_X^0)$ implies  that there are no local-to-global obstructions to deform $X$. More precisely:
\begin{prop}\label{Manetti}  Let $X$ be a surface with log-terminal singularities. If $H^2(X,T_X^0)=0$, then
$$\phi:Def(X) \rightarrow \times_{p\in \mbox{Sing}(X)} Def(V_p),$$
is surjective.
\end{prop}  
A proof can be found in Manetti's PhD Thesis \cite{M95} Chapter $3$. 

In order to understand the deformation theory of orbifolds, it is then essential to study the (local) deformation theory of quotient singularities. Let $V_p$ be an isolated quotient singularity. Then it is known (compare \cite{C11}) that $Def(V_p)$ can be identified with an analytic subset of the finite dimensional $\C$-vector space $\mathcal{E}xt^{1}(\Omega^1_{V_p}, \mathcal{O}_{V_p})_{p}$ with equations given by an obstruction map
$$ob:\mathcal{E}xt^{1}(\Omega^1_{V_p}, \mathcal{O}_{V_p})_{p} \rightarrow \mathcal{E}xt^{2}(\Omega^1_{V_p}, \mathcal{O}_{V_p})_{p}.$$
Before proceeding, let us see two examples of two singularities which are important for us. Let $V_p=\{f=0\} \subseteq \C^3$ be the germ of a canonical singularity. An important Theorem of G. Tjurina \cite{T69} and, independently, A. Kas and M. Schlessinger \cite{KS72} says that for two dimensional canonical singularities there are no obstructions to the deformations and $Def(V_p)$ can be identified with an open subset of the stalk of $\mathcal{E}xt^1$, which in this case is simply given by  
$$Def(V_p) \subseteq \mathcal{E}xt^{1}(\Omega^1_{V_p}, \mathcal{O}_{V_p})_{p} \cong \mathcal{O}_{\C^3,0}/ \left(f,\frac{\partial f}{\partial x_1},\frac{\partial f}{\partial x_2},\frac{\partial f}{\partial x_3}\right).$$
For example if $f(x_1,x_2,x_3)=x_1^2+x_2^2+x_3^{k}$, i.e., an $A_{k-1}$ singularity, then we can take $\{x_{3}^{k-2},x_{3}^{k-3},\dots, x_3, 1 \}$ as a $\C$-basis of $\mathcal{E}xt^{1}(\Omega^1_{V_p}, \mathcal{O}_{V_p})_{p}$.
In the above case the total family of the deformation space is given by $\{F(x,t)=0\} \subseteq \C^3_x\times \C^{k}_t$ where
$$F(x,t)= x_1^2+x_2^2+x_3^{k}-t_{1}x_{3}^{k-2}-t_2x_{3}^{k-3}-\dots -t_k.$$ 

The other interesting example is given by considering the cone over the degree $4$ Veronese embedding of $\p^1$. This is a quotient singularity given by the action of $\Z_4$ on $\C^2$ generated by $iId$. As first observed by H. Pinkham, this singularity admits a reducible (hence singular) deformation space $Def(\C^2/ \Z_4)$. $Def(\C^2/ \Z_4) \subseteq \mathcal{E}xt^{1}(\Omega^1_{V_p}, \mathcal{O}_{V_p})_{p} $ consists of two (smooth) irreducible components, of dimensions $1$ and $3$.
The one dimensional component can be described as follows. Consider the ordinary node $x_1^2+x_2^2+x_3^2=0$ and its smoothing $x_1^2+x_2^2+x_3^2=t$. These spaces (and the total family of the deformations) admit a holomorphic $\Z_2$ action given by $x_i\rightarrow -x_i$. It is easy to check that $\C^2/\Z_4 \cong (\C^2/\Z_2)/\Z_2$. Thus the $1$ dimensional component parameterizes  the smoothing of $\C^2 /\Z_4$ given by the $\Z_2$-quotient of the smoothing of the node.

The three dimensional component admits the following description (compare \cite{S03}, first Chapter). We can embed  $\C^2/\Z_4$ in the affine space, using invariant functions, as a cone over the degree four Veronese embedding of $\p^1$. Its equations are given by
$$\mbox{rk} \left( \begin{array}{cccc} z_0 & z_1 & z_2 & z_3 \\ z_1 & z_2 & z_3 & z_4 \end{array} \right) \leq 1.$$
The three dimensional component parametrizes deformations of the form
$$\mbox{rk} \left( \begin{array}{cccc} z_0 & z_1+s_1 & z_2+s_2 & z_3+s_3 \\ z_1 & z_2 & z_3 & z_4 \end{array} \right) \leq 1.$$
 $\mathcal{E}xt^{1}(\Omega^1_{V_p}, \mathcal{O}_{V_p})_{p}$ is a four dimensional vector space. We can take as coordinates $t$ and $s_i$ to be the deformation parameters defined above. Then it is known that the deformation space is singular (in particular it is reducible) and given by equations:
$$Def(\C^2/\Z_4)=\{ts_1=ts_2=ts_3=0\}\subseteq \C^4 \cong \mathcal{E}xt^{1}(\Omega^1_{V_p}, \mathcal{O}_{V_p})_{p}.$$

The above example is a special case of an important class of surface singularities, i.e., singularities which admit a $\Q$-\emph{Gorenstein smoothing}.  By definition a two dimensional quotient singularity admits a $\Q$-Gorenstein smoothing if  it admits a deformation over the disc $\mathcal{X} \rightarrow \Delta \subseteq \C$ with smooth generic fiber and with $\Q$-Gorenstein total space $\mathcal{X}$. Such singularities have been classified. We survey the results in the next Theorem.

\begin{thm}[J. Koll\'ar, N. Shepherd-Barron \cite{KSB88} or M. Manetti \cite{M90}]\label{TS} Let $V_p$ be a germ of a two dimensional quotient singularity. Then $V_p$ admits a $\Q$-Gorenstein smoothing if and only if either $V_p$ is given by the quotient of $\C^2$ by the action of $\Z_{dn^2}$ with weights $(1,dna-1)$ with $n,a$ coprime or $V_p$ is canonical. These singularities are known as $T$-singularities.
 
Moreover every $T$-singularity which is not canonical, is obtained by the quotient of a $A_{dn-1}$ singularity, i.e., $xy=z^{dn}$, by the  $\Z_n$ action: $(x,y,z)\rightarrow (\mu x, \mu^{-1} y, \mu^{a} z)$, with $\mu^n=1$.

Finally the (versal) family of  $\Q$-Gorenstein deformations is given by the $\Z_n$ quotient of the hypersurface
$$\{xy=z^{dn}+\sum_{i=1}^{d }t_i z^{(d-i)n}\} \subseteq \C^3\times \C^d,$$ 
where the action is generated by $(x,y,z,t_1,\dots,t_d)\rightarrow (\mu x, \mu^{-1} y, \mu^a z, t_1, \dots, t_d)$. In particular the vector space $\C^d_t$ can be identified as a vector subspace of $\mathcal{E}xt^{1}(\Omega^1_{V_p}, \mathcal{O}_{V_p})_{p}$.
We denote this subspace by $\G ({T^1_X}_p)$.
\end{thm}

\begin{rmk}
We say that a T-singularity is of type $(d,n,a)$ if it is given by the quotient of $\C^2$ by the action of $\Z_{dn^2}$ with weights $(1,dna-1)$ with $n,a$ coprime. 
\end{rmk}

We can  state the final result on the  deformation theory of (smoothable) Del Pezzo orbifolds:

\begin{thm}\label{DDP} Let $X$ be a $\Q$-Gorenstein smoothable Del Pezzo orbifold, i.e.,  $X$ has at most $T$-singularities. Then there exists an open subset $\G Def(X)$  of a vector subspace of $\T^1_X$ which can be identified with the versal deformation space parameterizing $\Q$-Gorenstein deformations of $X$. More precisely 
 $$\G Def(X) \subseteq \G(\T^1_X) \subseteq \T^1_X,$$
where $\G(\T^1_X):= H^{1}(\Theta_X^{orb})\oplus \bigoplus_{p\in \mbox{Sing}(X)} \G ({T^1_X}_p) $.
\end{thm}
Note that if the singularities of $X$ are canonical, then $\G(\T^1_X)=\T^1_X$. Moreover, every $\Q$-Gorenstein deformation of $X$ is a Del Pezzo surface of the same degree (compare Proposition 1.41 in \cite{KM98}).
Finally we remark on one way to think about the above Theorem \ref{DDP}. The part $H^{1}(\Theta_X^{orb})$ parametrizes the deformations of the orbifold which fix the singularity types, while $\bigoplus_{p\in \mbox{Sing}(X)}( \G {T^1_X}_p) $ describes (local) deformations of the singularities (i.e., smoothings or jumps of the singularity type).

\section{Deformations of KE Del Pezzo orbifolds}

In this section we state a Conjecture about the behavior of the KE metric under small deformations of smoothable KE Del Pezzo surfaces. Our Conjecture is motivated by two results. The first one is  Sz\'ekelyhidi's Theorem on the deformation of smooth cscK manifold (\cite{S10}).  The second one is the existence theory of ``ALE'' CY metrics on $\Q$-Gorenstein deformations of T-singularities due to P. Kronheimer for the canonical (hyperk\"ahler) case \cite{K89} and  to Y. Suvaina for the remaining situations \cite{S12}:
\begin{thm}[P. Kronheimer and Y. Suvaina] Let $\mathcal{X} \rightarrow \Delta$ be  a $\Q$-Gorenstein smoothing of a T-singularity. Then $X_t$ admits an ALE CY metric (for each K\"ahler class).
 \end{thm}
 Note that the metric on the non-canonical $T$-singularities is simply given by taking the quotient by the $\Z_n$-\emph{isometric} action on the Kronheimer ALE metric on the canonical cover (see \ref{TS}).

Thus our  expectation is that ``stable'' deformations of a $\Q$-Gorenstein smoothable Del Pezzo surface carry KE (orbifold) metrics which, after rescaling, are modeled on the ALE CY metrics locally near the deformation of the singularities. 
More precisely, we formulate the following natural Conjecture:
 
\begin{conj}[KE Deformations Conjecture]  \label{KEDC} Let $(X_0,\omega_0)$ be a $\Q$-Gorenstein smoothable KE orbifold Del Pezzo surface and let 
$$\begin{array}{ccc}
   X_0 & \hookrightarrow & \mathcal{X} \\
       &                 & \downarrow \\
       &                 & \G Def(X_0) \subseteq \G(\T^1_{X_0}) \subseteq \T^1_{X_0}
  \end{array}
$$
 be the versal $\Q$-Gorenstein deformation family, where $\G(\T^1_{X_0})$ denotes the subspace of $\T^1_{X_0}:=Ext^1(\Omega_{X_0}^1,\mathcal{O}_{X_0})$ parametrizing $\Q$-Gorenstein deformations $X_t \subseteq \mathcal{X}$.

Consider the natural action of the reductive group $Aut(X_0)$ on $\G(\T^1_{X_0})$. Then there exists an open set $B \subseteq \G Def(X) $ containing the origin such that for all $t \in \G Def(X)$ we have:
\begin{itemize}
 \item $t$ is polystable $\Rightarrow$ $X_t$ admits a (unique) KE orbifold metric $\omega_t$.
  \item $t$ is not-polystable $\Rightarrow$ $X_t$ cannot admit a KE orbifold metric.
\end{itemize}

Moreover there exists a continuous function $f:B\rightarrow \R^{\geq 0}$ with $f(0)=0$ such that
$$d_{GH}\left( (X_t,\omega_t), (X_0,\omega_0) \right) \leq f(t),$$
for all polystable $t$.
\end{conj}

If the above Conjecture is true, we will have a complete understanding of the local structure of the compactified moduli of KE Del Pezzo surfaces at the boundary points. In the next Chapter we are going to show that a special case of the Conjecture actually holds. We will show that sufficiently generic partial smoothings of nodal Del Pezzo surfaces (i.e. with $A_1$ singularities) with discrete automorphism group admit KE metrics. 
In Chapter $4$ we will describe some possible applications of the Conjecture to the study of the KE compactification problem. For example, we will discuss carefully the case of Del Pezzo surfaces of degree $3$ (cubic surfaces).

At this point it is interesting to recall a Conjecture of Tian:

\begin{conj}[G. Tian](\cite{T06}, pag. 387)
 GH limits of smooth KE Del Pezzo surfaces admit at worst canonical singularities.
\end{conj}

As we will recall in the next Chapter, this Conjecture is known to be true for the case of degree $4$ Del Pezzo surfaces. As we will explain later, we believe that also the space of KE Del Pezzo cubics is compactified by adding singular Del Pezzo surfaces with only canonical singularities.
However, we think that the situation starts to be different in the case of Del Pezzos of degrees $2$ and $1$. The main reason for this belief  is the following observation:    

\begin{prop}\label{12}
There exist $\Q$-Gorenstein smoothable KE Del Pezzo orbifolds of degree $1$ and degree $2$ with log-terminal (but not canonical) singularities.
\end{prop}

\begin{dimo}
Consider the KE toric Del Pezzo surface $Y_2$ given in Proposition \ref{Yn}. The surface $Y_2$ has degree $2$ and four singularities. Two of them are canonical (i.e., of type $A_3$) and two are $T$-singularities of type $(1,2,1)$ (i.e., the cone over the Veronese embedding of $\p^1$ of degree $4$ previously described). By Theorem \ref{DDP}, $Y_2$ admits $\Q$-Gorenstein smoothings.
It is easy to see that the versal space of $\Q$-Gorenstein deformations is given by
$$Def(Y_2)\subseteq 2\C^{3} \oplus 2 \C,$$
where $\C^3$ parametrize deformations of the $A_3$ singularities and $\C$ is the $1$-dimensional $\Q$-Gorenstein deformation component of the smoothings of the $(1,2,1)$ singularities. 

Observe that the virtual dimension of the moduli space at this point is $6$ (you need to subtract from the dimension of $Def(Y_2)$, which is $8$, the dimension of the automorphism group, which is $2$). This is of course equal to the dimension of the moduli space of (smooth) degree $2$ Del Pezzo surfaces.

The degree one example $Y_1$ comes from the classification of $\Q$-Gorenstein smoothable toric KE Del Pezzo orbifolds of Picard rank one (which will we prove in the next paragraphs). Let $Y_1:=X_P$ as defined in Theorem \ref{TP1}. Since all the singularities are of type $T$, $Y_1$ admits $\Q$-Gorenstein smoothings. In this case
$$Def(Y_1) \subseteq \C^8 \oplus \C \oplus \C,$$
where $\C^8$ is the space of versal deformations of the $A_8$ singularity and the two $\C$ are the space of $\Q$-Gorenstein deformations of the other $T$-singularities.
The virtual dimension of the moduli space is $8$, as is the dimension of the moduli of smooth degree $1$ Del Pezzo surfaces.

\qed \\
\end{dimo}

Finally, we give the classification of $\Q$-Gorenstein smoothable toric KE Del Pezzo surfaces with Picard rank $1$. The degree $1$ example in the previous Proposition emerges from this classification. A similar result could be obtained (but we didn't carry out the computations) by browsing the list of smoothable Del Pezzo surfaces of Picard rank $1$ in P. Hacking and Y. Prokhorov \cite{HP10} and computing the Futaki invariant in the toric case. However we decided to present an easy, self-contained proof which makes use of the KE condition from the very beginning.

First of all, we need a combinatorial Lemma:

\begin{lem} \label{TD} Let $P$ be a Del Pezzo polytope whose vertices are given by the primitive integer vectors $\{v_0,..., v_{n-1}\}$, in a clockwise order. Then the degree of the associated toric variety $X_P$ is given by 
 $$\mbox{deg}(X_{P})=2 \sum_{i=1}^{n}\frac{1}{|det(v_{i-1},v_i)|} + \sum_{i=1}^{n}\frac{|det(v_{i-1},v_{i+1})|}{|det(v_{i-1},v_i)||det(v_{i},v_{i+1})|}$$
where $v_0=v_n$ and $(v_{i-1},v_i)$ denotes the $2\times 2$ matrix with $v_i$ as column vectors.

\end{lem}

\begin{dimo}
Since $\mbox{deg}(X_{P})=2 \mbox{Vol}(P^{\checkmark})$, we need to express the volume of the dual polytope in term of the primitive generators $v_i$ of the rays of the fan. 

It follows by definition that the vertices of the dual polytope are given by $w_i=\left(\frac{v_{i-1}^y-v_{i}^y}{ |det(v_{i-1},v_i)|}, \frac{v_{i}^x-v_{i-1}^x}{ |det(v_{i-1},v_i)|} \right)$. Since
$2\mbox {Vol}(P^{\checkmark})= \sum_{i=0}^{n} |det(w_{i-1},w_i)|,$ 
the Lemma follows by a trivial computation.
\qed
\end{dimo}

With this Lemma in mind, it is easy to obtained the desired classification (and the degree one example we need).
\begin{thm}\label{TP1}
A $\Q$-Gorenstein smoothable KE toric Del Pezzo orbifold of Picard rank one is one of the following:
  \begin{itemize}
  \item $\p^2$ (degree $9$ case);
  \item $\{xyz=t^3\} \subseteq \p^3$ with $Sing(X)=3A_2$ (degree $3$ case);
  \item a variety $X_P \subseteq \p^6$ with $Sing(X_P)= A_8+2(1,3,1)$ (degree $1$ case).
  \end{itemize} 
\end{thm}

\begin{dimo}
 The Picard rank condition $\rho(X_{P})=1$ implies that $P$ must be a triangle, $P=\mbox{Convex}\{v_0,v_1,v_2\}$. By \cite{SZ11} we know that $X_{P}$ admits a KE (orbifold) metric iff the barycenter $Bar(P^{\checkmark})=0$, i.e., if
$w_0+w_1+w_2=0$, where $w_i$ are the vertices of the dual polytope $P^{\checkmark}$. It is easy to see that the balancing condition is equivalent to $v_0+v_1+v_2=0$ for the vertices of the Del Pezzo polytope $P$.

Therefore we may assume w.l.o.g. that $$v_0=(0,1) \; v_1=(-k,-l) \; v_2=(k, l-1),$$
for positive integers $k,l$. Applying Lemma \ref{TD} to our polytope, we find that
$$deg(X_p)=2\left( \frac{1}{k}+\frac{1}{k}+\frac{1}{k}\right)+\frac{3}{k}=\frac{9}{k}.$$
Since the degree is preserved under $\Q$-Gorenstein smoothings and the degree of a smooth Del Pezzo surface must be a positive integer less than or equal to $9$,  only the three possibilities $k=1,3,9$ survive.

\begin{itemize}
 \item Degree $9$ case ($k=1$). The unimodular matrix
 $$  
 \left(
  \begin{array}{cc}
    1 & 0 \\
    1-l & 1 \\
  \end{array}
\right)
 $$ 
 transforms the standard polytope of $\p^2$, given by $$ \mbox{Convex}\{(0,1),(1,0),(-1,-1)\},$$ into our polytope  $P= \mbox{Convex}\{(0,1),(-1,-l),(1, l-1)\}$. Thus the only (smoothable) toric degree $9$ KE Del Pezzo  is  $\p^2$.

\item Degree $3$ case ($k=3$). Since the vertices need to be primitive integer vectors,  $l\equiv 2$ mod $3$. Then the unimodular matrix
 $$  
 \left(
  \begin{array}{cc}
    1 & 0 \\
    \frac{2-l}{3} & 1 \\
  \end{array}
\right)
 $$ 
 transforms  every polytope  $P= \mbox{Convex}\{(0,1),(-3,-l),(3, l-1)\}$ into the polytope given by taking $l=2$. Using the unimodular transformation $  
 \left(
  \begin{array}{cc}
    0 & -1 \\
    1 & -1 \\
  \end{array}
\right),
 $ we see that the polytope for $l=2$ is sent to a polytope whose dual is given by $P^{\checkmark}=\mbox{Convex}\{(0,1),(1,0),(-1,-1)\}$. Thus taking Proj of the algebra associated to the cone over $P^{\checkmark}$, we find that
$X_P= \mbox{Proj} \left(\C[x,y,z,t] / <xyz=t^3>\right).$ Notice that $Sing(X_P)=3A_2$.

\item Degree $1$ case ($k=9$). The vertices of the polytopes are primitive only if:
\begin{enumerate}
 \item $l\equiv 2$ mod $9$;
\item  $l\equiv 5$ mod $9$;
\item  $l\equiv 8$ mod $9$.
\end{enumerate}
As before, we can see that up unimodular transformations only three cases remains: $l=2,5,8$. However it is easy to see that these three cases are also equivalent.  The polytope for $l=5$ (resp. $l=8$) is transformed into the polytope for $l=2$ by the matrix $\left(
  \begin{array}{cc}
    4 & -9 \\
    1 & -2\\
  \end{array}
\right),$ (resp.$\left( \begin{array}{cc}
    -8 & 9 \\
    -1 & 1\\
  \end{array}
\right).$)

It is easy to verify that the singularities are given by one $A_8$ singularities and of type T with weights  $(1,3,1)$ (in particular are all $\Z_9$ quotients). Since these singularities are locally smoothable and since there are no local-to-global obstructions, we can apply Theorem \ref{DDP}  and conclude that the above toric variety is smoothable.

Finally, a simple computation shows that the dual polytope of this toric variety is  $P^{\checkmark}=\mbox{Convex}\{(0,-1),(-\frac{1}{3},2),(\frac{1}{3},-1)\}$. By standard toric geometry it follows that $-K_{X_P}$ is only $\Q$-Cartier and that $-3K_{X_P}$ is a very ample Cartier divisor whose linear system has dimension $\mbox{dim}\,H^{0}(X_P,-3K_{X_P})=7$.

\end{itemize}

\qed
\end{dimo}

It is natural to ask what happens if we remove the restriction on the Picard rank in the classification of $\Q$-Gorenstein smoothable toric KE log Del Pezzo surfaces. We asked A.  Kasprzyk to look inside his (singular) toric Del Pezzo database \cite{GRDb}.  We found that the examples of Del Pezzo surfaces of degree $1$ and $2$ previously considered, together with the ``obvious'' Gorenstein ones, are indeed all the possible $\Q$-Gorenstein smoothable KE toric log del Pezzo. More precisely: 
\begin{thm}
 Let $(X,\omega)$ be a $\Q$-Gorenstein smoothable toric KE log Del Pezzo surface. Then $X$ is isomorphic to one of the following:
\begin{itemize}
 \item $deg(X)=9$: $X \cong \p^2$, (with Picard rk $\rho(X)=1$);
 \item $deg(X)=8$: $X \cong \p^1 \times \p^1$, ($\rho(X)=2$);
 \item $deg(X)=6$: $X \cong Bl_{[1:0:0],[0:1:0],[0:0:1]}\p^2$, ($\rho(X)=4$);
 \item $deg(X)=4$: $\{xy=z^2=tw\} \subseteq \p^4$, ($\rho(X)=2$);
 \item $deg(X)=3$: $X \cong \{xyz=t^3\} \subseteq \p^3$, ($\rho(X)=1$);
 \item $deg(X)=2$: $ X \cong Y_2$ defined as in Proposition \ref{12}, ($\rho(X)=2$);
 \item $deg(X)=1$: $X \cong X_P$ defined as in  Theorem \ref{TP1}, ($\rho(X)=1$).
\end{itemize}
\end{thm}

\begin{dimo}
It follows by Theorem \ref{NDP} and by the classification of $\Q$-Gorenstein smoothable singularities Theorem \ref{TS}, that the index, i.e. the smallest positive integer such that $-lK_X$ is Cartier, is bounded (actually $\leq 6$). Then it is known that the set of toric log Del Pezzo surfaces (with T-singularities) of index at most $l$ is finite (see \cite{KKN10}). The statement follows by computing the barycenter (i.e., the Futaki invariant) of all the finitely many toric log Del Pezzo surfaces of index at most six, collected in the database \cite{GRDb}.
\qed
\end{dimo}

\begin{rmk}
 The proof of the above Theorem has the disadvantage of requiring a case-by-case check of the Futaki invariant. It may be interesting to find a ``conceptual'' proof of the above result (e.g., a proof similar to the one we gave in the Picard rank one case).
\end{rmk}

\end{chapter}
\begin{chapter}{Deformations of nodal KE Del Pezzo surfaces with discrete automorphism groups}
 
The Theorem proved in this Chapter can be considered a special case of  Conjecture \ref{KEDC}. We show that generic (partial)-smoothings of a nodal KE Del Pezzo surface $X_0$ (i.e., with singularities  at most of type $A_1$) with discrete automorphism group admit orbifold KE metrics close in the GH sense to the original singular KE metric on $X_0$.
As an application, we show that KE orbifold metrics exist on some cubics with two or three $A_1$ singularities.
 
\section{Statement of the Theorem, remarks and first applications}
The main result of this Chapter is the following:

\begin{thm}\label{MT} Let  $(X_0,\omega_0)$  be a KE Del Pezzo orbifold. Suppose that
\begin{itemize}
\item All the singularities are nodes (i.e., locally of the form $z_1^2+z_2^2+z_3^2=0)$;
\item $\sharp (Aut(X_0)) < \infty$.
\end{itemize}
Then, if $\pi:\mathcal{X} \rightarrow \Delta \subseteq \C_t$ is a generic (partial) smoothing of $\pi^{-1}(0)=X_0$, $X_t$ admits an (orbifold) KE metric $\omega_t$ for $|t|<<1$. Moreover $(X_t,\omega_t) \rightarrow (X_0,\omega_0)$ in the Gromov-Hausdorff (GH) sense.
\end{thm} 

By the word ``generic'' we mean that if $t$ is the parameter of the base of the smoothing family $\mathcal{X}$ and if $z_1^2+z_2^2+z_3^2=s \in \C^{3}\times \C_s$ is the total family of the versal deformation of the node \cite{KS72}, then $s$ and $t$ are related by $$s=s(t)=Ct+\mathcal{O}(t^2)$$ with $C\neq 0$.

The reason why we need the above assumption is essentially technical and it has mainly to do with our choice of function spaces in which to perform the gluing construction. Removing the above genericity assumption gives rise to a concentration of the error term on the variety $X_t$ away from the singularities (more precisely the Ricci potential of the pre-glued metric becomes too big to be controlled). Nevertheless, we believe that the Theorem should be true even without the genericity assumption.  

Before giving the proof of the main Theorem, it is interesting to make some remarks.
First of all, in our considerations on moduli spaces we are principally interested in Del Pezzo orbifolds which appear in GH degenerations of smooth KE Del Pezzo surfaces, that is in smoothable Del Pezzos. In the next Proposition we see that such smoothable Del Pezzo surfaces with discrete automorphism group cannot have many nodal singularities. 
\begin{prop}Let $X_0$ be a smoothable log Del Pezzo surface with discrete automorphism group and only nodal singularities. Then
$$\sharp\{\mbox{nodes}\}\leq 10-2\, deg(X_0).$$
\end{prop}

\begin{dimo}
Let $\mathcal{X} \rightarrow \Delta_t \subseteq\C$ be a full smoothing of $X_0$ and define 
$$\varphi(t):= \mbox{dim}\,Ext^{0}\left(\Omega_{X_t}^1,\mathcal{O}_{X_t}\right)-\mbox{dim}\,Ext^{1}\left(\Omega_{X_t}^1,\mathcal{O}_{X_t}\right)+\mbox{dim}\,Ext^{2}\left(\Omega_{X_t}^1,\mathcal{O}_{X_t}\right).$$ 
Since $\varphi$ is upper semi-continuous and $X_t$ is a smooth Del Pezzo of degree $deg(X_0)$ for all $t\neq 0$, we find 
 \begin{itemize}
     \item  $\varphi(0) \geq \varphi(t)=\chi (\Theta_{X_t})=2\, deg(X_0)-10$ by Proposition \ref{VDM};
      \item $Ext^{0}\left(\Omega_{X_0}^1,\mathcal{O}_{X_0}\right)=0$ since $Aut(X_0)$ is finite and the singularities are normal. 
    \end{itemize}
Thus $$ \mbox{dim}\,Ext^{1}\left(\Omega_{X_0}^1,\mathcal{O}_{X_0}\right) \leq 10-2\, deg(X_0).$$
Finally we can estimate the dimension of the above $Ext^1$ from below by the number of nodes. Using the Grothendieck's local-to-global Ext spectral sequence, the vanishing of $H^{2}(\mathcal{H}om(X,\Omega_{X_0},\mathcal{O}_{X_0}))$, and recalling that $\mathcal{E}xt^{1}\left(\Omega_{X_0}^1,\mathcal{O}_{X_0}\right)$ is a skyscraper sheaf supported on the singular locus with stalk isomorphic to the space of versal deformation of the singularities, we have
$$\begin{array}{rl} \mbox{dim}\,Ext^{1}\left(\Omega_{X_0}^1,\mathcal{O}_{X_0}\right) & \geq  \mbox{dim} \, H^0 \left( \mathcal{E}xt^{1}\left(\Omega_{X_0}^1,\mathcal{O}_{X_0}\right)\right) \\ & =\mbox{dim}\, \left( \oplus_{p \in Sing(X_0)} \C_p  \right) \\ &= \sharp\{\mbox{nodes}\}.\end{array}$$  
\qed
\end{dimo}

\begin{rmk}
 The above Proposition follows also by the \emph{classification} of Del Pezzo surfaces with canonical singularities \cite{BW79}, \cite{D80}, \cite{HW81}, \cite{U83}, \cite{FU86}, \cite{Z88}, \cite{MZ88} and \cite{MZ93}.
\end{rmk}

What can we say about smoothable KE Del Pezzo orbifolds (with the above properties) for each interesting degree $d=(-K_{X_0})^2$?  It is known that there are smoothable del Pezzo quartics, i.e., degree equal to $4$, with no holomorphic vector fields and with only one or two nodal singularities. On the other hand,  it follows by the Mabuchi-Mukai result \cite{MM90} (compare the next Chapter) that all KE Del Pezzo surfaces appearing in the GH compactification must have holomorphic vector fields. Even if in this case we cannot use our Theorem to study the local behavior of the moduli space at the boundary points, we can combine it with the Mabuchi-Mukai result to prove the following observation:

\begin{cor} \label{NKE} Let $X$ be a Del Pezzo quartic with only  a nodal singularity and with discrete automorphism group. Then $X$ does not admit a KE metric.
\end{cor} 
   
\begin{dimo} 
 By the classification of intersections of two quadrics (compare for example \cite{AL00}), we can assume that $X_0$ is of the form
\begin{displaymath}
\left\{ \begin{array}{ll}
2 x_0 x_1+x_2^2 +x_3^2+x_4^2=0; & \\
 2 x_0 x_1 +x_0^2 +  x_2^2+a x_3^2+b x_4^2=0, & \end{array} \right.
\end{displaymath}
for generic $a,b\in \C$. Observe that $Aut(X_0)$ is finite. In order to see this, note that the above $X_0$ can be degenerated, by the one-parameter subgroup $T_s$, given by $x_0 \rightarrow  s x_0$, $x_1 \rightarrow s^{-1} x_1$,  to a KE Del Pezzo quartic $X_c$ with two nodes, where the two defining equations are both diagonalized (after an obvious change of coordinates). The automorphism group of $X_c$ can be easily computed by looking at the dimension of the Lie group of matrix in $SL(5,\C)$ which fix the defining equations (note that this group is the full automorphism group since $-K_{X_C}$ is very ample). By a simple computation, it follows that the Lie algebra is one dimensional. Hence $\mbox{dim}_{\C}Aut(X_0)< \mbox{dim}_{\C}Aut(X_c)=1$.   

A generic deformation (smoothing) $X_t$ is given by
\begin{displaymath}
\left\{ \begin{array}{ll}
2x_0 x_1+x_2^2 +x_3^2+x_4^2=0; & \\
2x_0 x_1 +x_0^2 +t x_1^2+  x_2^2+a x_3^2+b x_4^2=0. & \end{array} \right.
\end{displaymath}
If $X_0$  has a KE metric, then by Theorem \ref{MT}  $(X_t,\omega_t)$  degenerates in the GH sense to $X_0$. In particular $X_0$ must appear in the GH compactification. However this is forbidden by the Mabuchi-Mukai result (see next Chapter).\qed
\end{dimo}

As we pointed out in the proof of the previous Corollary, an $X_0$ satisfying the hypothesis of the Corollary can be degenerated by a one parameter subgroup $T_s$ to a KE Del Pezzo quartic $X_c$ with two nodes (and with $\C^{*} \subseteq Aut(X_c)$). In particular the Donaldson-Futaki invariant $DF(T_s.X_0)=0$ (since the central fiber admits a KE metric) and $T_s.X_0$ is  a non-trivial test configuration.
We should recall that Stoppa's argument for proving that the existence of a KE (csck) metric implies K-polystability \cite{S09} doesn't extend easily in the case of singular spaces. However, in \cite{B12} the author showed that ``KE implies K-polystability'' is indeed true for $\Q$-Fano varieties. Thus Corollary \ref{NKE} follows also by this more general result.
 
It is more interesting to see what happens for degree three Del Pezzos, i.e., cubic surfaces. As we will recall in the next Chapter, it follows by the well-known classical GIT pictures for cubic surfaces that cubics with nodal singularities have no holomorphic vector fields. The (unique) cubic with the maximal number of nodal singularities is the so-called Cayley's cubic given by the equation $$C_0:=\{xyz+yzw+zwx+wxy=0\} \subseteq \p^3.$$
Moreover it known that  that $C_0$ admits an orbifold KE metric (e.g., \cite{C09}, or simply by noting that it is a quotient of the KE Del Pezzo surface of degree $6$). Thus we can use our Theorem to prove the following:
\begin{cor}\label{ex}Some cubic surfaces with two or three nodal singularities admit KE orbifold metrics. 
\end{cor}
\begin{dimo}
 It is sufficient to consider deformations of $C_0$ of the form:
$$C_t:=\{xyz+yzw+zwx+wxy+tx^3+(ty^3)=0\} \subseteq \p^3.$$\qed
\end{dimo}

We also recall that there are several examples of Del Pezzo surfaces of degree $1$ and $2$ to which our Theorem applies.

The proof of the Theorem is based on a gluing construction. In the next section, we explain how to construct an approximate solution by gluing an Eguchi-Hanson space. Then we will apply some analytical argument in order to deform the approximate solution to a genuine KE metric. Finally, we show that the deformation process is continuous in the GH topology.

\section{Construction of an approximate solution}

Let $(X_0, \omega_0)$ be a nodal KE Del Pezzo orbifold and let $p$ be a singular point. Take $(\zeta_1,\zeta_2)$ coordinates on $\C^2/\Z_2$ around $p$ so that locally $$\omega_0=\Do (\varphi^1_0)=\Do (|\zeta|^2+\mathcal{O}(|\zeta|^4)).$$ We say that a function $f(\zeta,\bar{\zeta})$ is   $\mathcal{O}(|\zeta|^k)$ if
$ ||\nabla^j f || \leq C |\zeta|^{k-j}$ for $j\geq 0$. 

Identifying as usual  $\C^2/\Z_2= \{z_1^2+z_2^2+z_3^2=0\}$, where $z_i=z_i(\zeta_1,\zeta_2)$ are quadratic expression in $\zeta_i$, it is easy to see that the metric $\omega_0$  can be written as  $$\omega_0=\D (|z|+\mathcal{O}(|z|^2))_{|V_0}$$ in a neighborhood $V_0 \subseteq X_0$ of the origin.

 By  the theory of versal deformations of hypersurface singularities \cite{KS72}, we have that the general deformation of the node is given by  $ z_1^2+z_2^2+z_3^2=t$, with $t \in \C$. Observe that we may assume $t$ to be real (by ``rotating'' the coordinates). 

Thus, using our genericity assumption, we  \emph{identify} a portion of $z_1^2+z_2^2+z_3^2=t \in \R$ (say $|z|\leq C$ for some positive constant C) with a subset $V_t$ of $X_t$ for some real $t$.

Moreover, using the relative anticanonical sections,  we can assume that $\mathcal{X}$ embeds in $\p^n\times \Delta$ and the map $\pi$ is given by composing the anticanonical embedding with the projection onto the base. Since  $\mathcal{O}_{\p^n}(1)_{|X_t}=K_{X_t}^{-k}$ (for notational simplicity we assume $k=1$, since the power $k$ does not enter in our argument),  we have constructed a background K\"ahler metric  ${\beta_t}_{|X_t}:={FS_t}_{|X_t}$ in $c_1(X_t)$.

The goal of this section is to construct a  pre-glued metric $\omega_t$, i.e., a metric on $X_t$  GH close to the original metric $\omega_0$, of the form $$\omega_t=\beta_t+\Dt \phi_t \in c_1(X_t).$$
Of course the argument above can be applied to all singular points $p \in Sing(X_0)$.

\subsection*{A Map between the smoothing and the central fiber}

Assume for simplicity that $Sing(X_0)=\{p\}$.
Let  $V_0=\{z_1^2+z_2^2+z_3^2=0\}$ and $ V_{t}=\{w_1^2+w_2^2+w_3^2=t\}\subseteq X_t$ and consider the diffeomorphism 
$$\begin{array}{cccc}
F_t:& V_0\setminus \{|z|^2 \leq \frac{t}{2} \}  & \longrightarrow &V_{t} \setminus \{|w|^2= t\} \\
                                      &  z_i &  \longmapsto & w_i=z_i+\frac{t}{2|z|^2}\bar{z}_i 
\end{array},$$
between  the singular space and its smoothing (away from the singularities and defined on a small region of $X_0$). Observe that topologically, $$  V_0\setminus \{|z|^2 \leq \tfrac{t}{2} \}  \cong V_{t} \setminus \{|w|^2= t\}  \cong \R\times \R \p^3,$$
and that $L_t:= V_{t} \cap \{|w|^2= t\}\cong S^2$. Moreover note that $|z|\mapsto |w|= \sqrt{|z|^2+\frac{t^2}{4|z|^2}}$.

Let us point out that we can assume that $F_t$ extends to a diffeomorphism  $\psi_t:X_0 \setminus \{|z|^2 \leq \frac{t}{2} \}   \longrightarrow  X_t \setminus \{|w|^2= t\}$.

\begin{lem} There exists a  diffeomorphism $$\psi_t:X_0 \setminus \{|z|^2 \leq \tfrac{t}{2} \}   \longrightarrow  X_t \setminus \{|w|^2= t\},$$
which coincides with the above $F_t$ for $|z|^2 \leq 4$. 

Moreover the complex structures $J_t$ of $X_t$ and $J_0$ of $X_0$  satisfy on  $X_0 \setminus \{ |z|^2 \leq 1 \}$, $||\nabla^k_{\omega_0}( \psi^{*}_t J_t - J_0)||_{\omega_0}=\mathcal{O}(t)$ for all $k$ (or with respect to any other metric equivalent to $\omega_0$). 

\end{lem}

\begin{dimo} As we have previously remarked, we may assume that $\mathcal{X} \hookrightarrow \p^n\times \Delta$. Then $\mathcal{X}$ inherits a K\"ahler metric from the metric $\omega_{FS}+\D |t|^2$ on $\p^n \times \Delta$. We define a connection on $\mathcal{X}\setminus Sing(X_0)$ simply by taking the normal directions to the fiber of the map $pr_2:\mathcal{X} \subseteq \p^n\times \Delta \rightarrow  \Delta$. 

Considering the flow given by the lifting of the radial vector field $X:= -\frac{\nabla |t|^2}{|\nabla |t||^2}$, we find a smooth diffeomorphism  from $X_{|t|} \setminus K_{|t|}$ to $X_0 \setminus \{|z|^2 \leq 4\}$, for some small compact subsets of $K_{|t|}$ of $X_{|t|}$. Let $G_t:X_0 \setminus \{|z|^2 \leq 4\} \rightarrow X_t$ be the inverse of the previous map (assume w.l.o.g that $t$ is a point on the real path $[0,1]\subseteq \Delta$).

It follows that on the strip $S:=\{4\leq |z|^2 \leq 5\} \subseteq X_0$ we have defined two diffeomorphisms $F_t$ and $G_t$ onto some region of $X_t$. Let $S^{'} \subseteq S$ be some smaller strip and use the map $F_t$ to identify  $S$ with its image inside $X_t$. Under the above identification, $G_t$ can be seen as a family of diffeomorphisms onto their image
$$\tilde{G}_t: S^{'} \rightarrow S \cong [0,1] \times \R\p^3,$$
which are by construction smoothly isotopic to the identity on $S^{'}$, i.e. $$\tilde{G}_t \rightarrow Id$$ as $t \rightarrow 0$. In particular $\tilde{G}_t$ is given by the flow of a time dependent vector field $v_t$ defined on $S^{'}$, i.e.,
$$\frac{d}{dt}(\tilde{G}_t)= v_t(\tilde{G}_t).$$
Let $\tau: S \rightarrow [0,1]$ be a smooth cut-off function given by 
$$\tau(z)=\chi(|z|^2),$$
where $\chi$ is a smooth increasing function equal to zero for $|z|^2 \leq 4$ and equal to one for $|z|^2 \geq 5$. Let $w_t$ be the vector field defined by $w_t:=\tau v_t$ and $H_t$ its associated diffeomorphism (for small $t$).
By construction $H_t$ is equal to the identity for $|z|^2 \leq 4$ and equal to $\tilde{G_t}$ for $|z|^2 \geq 5$. Thus the diffeomorphism
$$
\psi_t:=  \left\{ \begin{array}{rl}
 F_t\circ H_t & \mbox{ if $|z|^2 \leq 5$}; \\
  G_t &\mbox{ if $|z|^2 \geq 5$},
       \end{array} \right. 
$$
satisfies the desired property. 

The estimate on the complex structures follows by our construction of the diffeomorphism $\psi_t$ which is given by the flow of a smooth vector field. More precisely, we can cover the total family $\mathcal{X}$ ``away'' from a neighborhood of the singularity with a finite number (by compactness) of charts $\mathcal{U}_j$ of the form $(\zeta_1^j,\zeta_2^j,t)$. Since  by construction $\phi_t$ is isotopic to the identity, in the chart $\mathcal{U}_j$ $\psi_t$ is given by $\zeta_i \mapsto \zeta_i +f^j(t,\zeta,\overline{\zeta})$ for a \emph{smooth} $f^j$ satisfying $f^j(0,\zeta,\overline{\zeta})=0$. Then the claim follows simply by considering the Taylor's expansion of $f^j$ in the variable $t$.

 \qed
\end{dimo}

\subsection*{Pre-glued metric close to the singularities}

On $V_1:=\{z_1^2+z_2^2+z_3^2=1\}\cong TS^2$ we have the CY conical Eguchi Hanson metric: $$\eta_1:= \D(\sqrt{|z|^2+1})_{|V_1}.$$ Note that $(\sqrt{|z|^2+1} -|z|) \leq C \frac{1}{|z|}$ for $|z|>>1$. Pulling-back by the holomorphic map $w_i=\sqrt{t} z_i$ and scaling the metric, we have a (scaling) family of  RF metrics on $V_t\subseteq  \{w_1^2+w_2^2+w_3^2=t\} \subseteq X_t$,
$$\eta_{\delta,t}=  \frac{\delta^2}{\sqrt{t}} \,\D \left(\sqrt{|w|^2+t}\right)_{|V_t},$$
satisfying $\mbox{Diam}_{\eta_{\delta,t}}(L_t)=\delta$.

Pulling-back at the level of the potential the KE metric $$\omega_0=\D (|z|+\mathcal{O}(|z|^2))_{|V_0}= \D (\varphi^1_0)_{|V_0}$$ using the (inverse) of $\psi_t$, we find a metric $\omega^1_t$ for small $t$, degenerate in a small neighborhood of $L_t$, given by 

$$\omega^1_t:= \D \left( {\psi^{-1}_t}^{*} (\varphi^1_0) \right)_{|V_t}.$$
%$$\omega^1_t:= \Dt\left(\sqrt{\frac{|w|^2+\sqrt{|w|^4-t^2}}{2}}+ \mathcal{O}\left(\frac{|w|^2+\sqrt{|w|^4-t^2}}{2}\right)\right).$$\\
We show how to glue $\eta_{\delta,t}$ to $\omega^1_t$. Take $\delta^2=\sqrt{t}$ and consider the gluing region $ \delta^\alpha\leq |w|_{|V_t} \leq 2 \delta^\alpha$ (topologically $[0,1] \times \R\p^3)$ for $\alpha \in[0,2)$ (observe that $\{|w|=\delta^2\}\cong L_t$). Define
\begin{itemize}
\item $\varphi_\delta^1(w):=  {\psi^{-1}_t}^{*} (\varphi^1_0)$;
\item $\varphi_\delta^2(w):= \sqrt{|w|^2+\delta^4}.$
\end{itemize}

\begin{lem} On the annulus $ \delta^\alpha\leq |w|_{|V_t} \leq 2 \delta^\alpha$, where $\alpha \in[0,2)$, we have:
 $$|\nabla_{\eta_\delta}^{k}(\varphi^1_\delta-\varphi^2_\delta)|_{\eta_\delta} =\mathcal{O}(\delta^{8-3\alpha-\frac{k\alpha}{2}})+ \mathcal{O}(\delta^{2\alpha-\frac{k\alpha}{2}})+ \mathcal{O}(\delta^{4-\alpha-\frac{k\alpha}{2}}).$$
 The error is minimized for $\alpha=\frac{4}{3}$: $|\nabla_{\eta_\delta}^{k}(\varphi^1_\delta-\varphi^2_\delta)|_{\eta_\delta} =\mathcal{O}(\delta^{\frac{2}{3}(4-k)})$.
\end{lem}

\begin{dimo}
First of all  note that on the annulus region $\eta_\delta \approx \Dt (|w|)_{|V_{\delta^4}}$, as the  metric (i.e., the scaled Eguchi-Hanson metric) is equivalent to the restriction of the (singular at the origin) metric $\D  (|w|)$). In particular this implies that  $|\nabla_{\eta_\delta}^k(|w|^j)|_{{\eta_\delta}} \leq C |w|^{j-\frac{k}{2}}.$

By the definition of $\varphi^1_\delta$ and $\varphi^2_\delta$ and since $|z|\mapsto |w|= \sqrt{|z|^2+\frac{t^2}{4|z|^2}}$ under the diffeomorphism $\psi_t$, it follows that

\begin{itemize}
 \item $|\nabla_{\eta_\delta}^{k}(\varphi^2_\delta- |w|)|_{\eta_\delta} \leq C \delta^4 |w|^{-1-\frac{k}{2}}$;
 \item $|\nabla_{\eta_\delta}^{k}(\varphi^1_\delta- \sqrt{\frac{|w|^2+\sqrt{|w|^4-\delta^8}}{2}} )|_{\eta_\delta} \leq C |w|^{2-{\frac{k}{2}}}$.
\end{itemize}
Thus
$$|\nabla_{\eta_\delta}^{k}(\varphi^1_\delta-\varphi^2_\delta)|_{\eta_\delta}\leq   C \left(|w|^{2-{\frac{k}{2}}} +  \delta^4 |w|^{-1-\frac{k}{2}}\right) + $$ $$+|\nabla_{\eta_\delta}^{k}\left(|w|- \sqrt{\frac{|w|^2+\sqrt{|w|^4-\delta^8}}{2}}\right)|_{\eta_\delta},$$
i.e.,
$$|\nabla_{\eta_\delta}^{k}(\varphi^1_\delta-\varphi^2_\delta)|_{\eta_\delta}\leq   C \left(|w|^{2-{\frac{k}{2}}} +  \delta^4 |w|^{-1-\frac{k}{2}} + \delta^8 |w|^{-3-{\frac{k}{2}}}\right).$$

\qed
\end{dimo}

\begin{lem}Define for $\sqrt{t}=\delta^2$ and $|w|_{|X_t} \leq 2$ $$\tilde{\omega}_{t,\delta}^1:=\Dt \left( \chi_\delta \varphi_\delta^1+(1-\chi_\delta)\varphi_\delta^2 \right),$$
where $\chi_{\delta}:=\chi\left(\delta^{-\frac{4}{3}}|w|\right)$ is a smooth increasing cut-off function supported in $|w| \geq \delta^{\frac{4}{3}}$, identically equal to one for $|w| \geq 2\delta^{\frac{4}{3}}$. Then for $\delta$ (hence $t$) sufficiently small
\begin{itemize}
\item $||\nabla_{\delta}^k (\tilde{\omega}_{t,\delta}^1-\eta_{\delta,t})||_{\eta_{\delta,t}} = \mathcal{O}(\delta^{\frac{4-2k}{3}})$;
\item $\tilde{\omega}_{t,\delta}^1>0$.
\end{itemize}

\end{lem}

\begin{dimo}
It follows immediately from the previous Lemma observing that $||\nabla^k_\delta \chi_\delta||_{\eta_{\delta,t}}=\mathcal{O}(\delta^{-\frac{2k}{3}})$ on the strip $|w| \in [\delta^{\frac{4}{3}},2 \delta^{\frac{4}{3}}]$ and zero otherwise.
\qed\\
\end{dimo}

\subsection*{The pre-glued metric away from the singularities and matching}

First of all we take the following  global K\"ahler potential for the singular metric $\omega_0$. Let $s$ be any  non-vanishing local section of $K_{X_0}^{-1}$. Then
$$\varphi_0:=\log \frac{|s|^2_{\beta_0}}{|s|^2_{\omega_0}},$$
is a well defined smooth function on $X_0 \setminus Sing(X_0)$. Observe that $\varphi_0$ is just the (log of) the ratio of the two volume forms.
Later we will need to fix a precise section of $K_{X_0}^{-1}$ near the singularities. The correct choice is to take  the section $\hat \Omega_0$ where $\hat\Omega_0$ is the section which pulls-back using the orbifold chart to $\partial_{\zeta_1} \wedge \partial_{\zeta_2}$ with $(\zeta_1,\zeta_2)$ the local charts where  $\omega_0= \delta_{ij}+\mathcal{O}(|\zeta|^2)$.

Since $\omega_0$ is K\"ahler-Einstein, it is  easy to see that
$$\omega_0=\beta_0+\Do \varphi_0.$$

Now we are ready to define the approximate KE metric away from the singularities (i.e., away from $|w|_{X_t} \leq 1$ on $X_t$):
$$\tilde{\omega}_{t,\delta}^2:=\beta_t +\Dt {\psi_t^{-1}}^{*} \varphi_0>0$$
for $t$ sufficiently small and where $\psi_t$ is the map between the smooth and singular fiber constructed before.

The next goal is to match the metrics $\tilde{\omega}_{t,\delta}^1$ and $\tilde{\omega}_{t,\delta}^2$ and construct a metric
$$\tilde{\omega}_{t,\delta} =\beta_t +\Dt \phi_t \in c_1(X_t).$$ 

Since $\eta_{\delta,t}$ is a CY metric, we know that there exists a (explicit) nowhere vanishing holomorphic $(2,0)$ form $\Omega_t$ such that for  $\sqrt{t}=\delta^2$
$$\eta_{\delta,t}^2= C \Omega_t\wedge\overline{\Omega}_t.$$
Let $\hat{\Omega}_t$ be its dual. Then $\hat{\Omega}_t$ is a trivialization of $K_{X_t}^{-1}$ on $V_t$ which satisfies $|\hat{\Omega}_t|_{\eta_{\delta,t}}=1$ (after normalization).

Define $b_t \in C^{\infty}(V_t)$ to be the function given by 
$$b_t:= |\hat{\Omega}_t|_{\beta_t}^2.$$
Then $\beta_t=-\Dt \log b_t$ on $V_t$. Thus $\tilde{\omega}_{t,\delta}^2=\Dt \left(- \log b_t + {\psi_t^{-1}}^{*} \varphi_0\right)$ on $V_t \cap \{1 \leq |w| \leq 2\}$. Since we have that $\tilde{\omega}_{t,\delta}^1=\Dt \varphi^1_\delta$ on the same region, it is natural to match the two metrics at the level of the potential using a cut-off function $\tau_t:=\tau(|w|_{V_t})$, a smooth decreasing cut-off function supported in $|w| \leq 2 $  identically equal to one for $|w| \leq 1$ (notice that here the strip where the cut-off function is non constant is of ``fixed shape'').

\begin{prop}\label{PG} There exists a pluriharmonic function $p_t$ on $V_t$ such that for $\sqrt{t}=\delta^2$
$$\tilde{\omega}_{t,\delta}:= \left\{ \begin{array}{ll}
\tilde{\omega}_{t,\delta}^1  & |w|_{|V_t} \leq 1 \\
\Dt \left( \tau_t\left(\varphi^1_\delta-p_t\right)+(1-\tau_t) \left(- \log b_t + {\psi_t^{-1}}^{*} \varphi_0\right) \right) &  1\leq |w|_{|V_t} \leq 2\\
\tilde{\omega}_{t,\delta}^2 &  \mbox{otherwise}\\
\end{array}  \right.$$
is a K\"ahler metric in $c_1(X_t)$. More precisely $\tilde{\omega}_{t,\delta}-\beta_t=\Dt \phi_t$ where
$$\phi_t=\tau_t \left( \left( \chi_\delta \varphi_\delta^1+(1-\chi_\delta)\varphi_\delta^2 \right) -p_t +\log b_t \right)+(1-\tau_t) \left( {\psi_t^{-1}}^{*} \varphi_0\right) \in C^\infty(X_t,\R).$$
Moreover we can choose $p_t$ to satisfy $|\nabla_{\eta_\delta}^{k} p_t|_{\eta_{\delta,t}}\leq C |w|^{1-\frac{k}{2}}$. 
\end{prop}

\begin{dimo}

Let $q$ be the singular point of $X_0$. By our hypothesis on the metric $\omega_0$ we know that, in the orbifold chart centered in $q$, $\omega_0$ can be expressed for $|\zeta| \leq \sqrt{2}$ as
$$\omega_0=\Do \varphi^1_0=-\Do \log (|\hat{\Omega}_0|^2_{\omega_0})=Ric (\omega_0),$$
where $\varphi^1_0=|\zeta|^2+\mathcal{O}(|\zeta|^4)$ is the local K\"ahler potential. Thus
$$p_0:= \varphi^1_0+ \log (|\hat{\Omega}_{0}|^2_{\omega_0}),$$
is a $\Z_2$-invariant pluriharmonic real function on  $|\zeta| \leq \sqrt{2}$  vanishing at the origin.
We claim that $p_0= \mathfrak{Re} (h_0)$ where $h_0$ is a $\Z_2$-invariant holomorphic function. Since $d \partial_0 p_0= \bar{\partial}_0 \partial_0 p_0=0$, by the Poincar\'e Lemma there exists a $\Z_2$-invariant holomorphic function $h_0$ vanishing at the origin satisfying
$$d\frac{h_0}{2} =\partial_0 \frac{h_0}{2}=\partial_0 p_0.$$
Then, being $p_0$ real,
$$d  (\mathfrak{Re} (h_0)-p_0)= \partial_0 \frac{h_0}{2} +\bar{\partial}_0 \frac{\bar{h}_0}{2} -\partial_0 p_0 - \bar{\partial}_0 p_0=0.$$
Hence $p_0= \mathfrak{Re} (h_0)$ .

Since $h_0$ is a holomorphic function on $\C^2 / \Z_2$ vanishing at the origin, identifying $\C^2 /\Z_2$ with $w_0^2+w_1^2+w_2^2=0 \in \C^3$, we may assume that $h_0=H(w_1,w_2,w_3)_{|X_0}$ where $H$ is an holomorphic function on $\C^3$ vanishing at the origin (since the node is a normal singularity). We define for sufficiently small $t$:  
$$p_t:= \mathfrak{Re} H_{|V_t},$$
which is a pluriharmonic function on the local smoothing for $|w|_{|V_t} \leq 2$, satisfying the desired estimates.

Now define $a_t:= \varphi^1_\delta-p_t$ and $c_t:=- \log b_t + {\psi_t^{-1}}^{*} \varphi_0$. It is evident that the closed $(1,1)$-form $\tilde{\omega}_{t,\delta}$ is well-defined and that it  would be positive definite as long as 
$$||a_t-c_t||+||d(a_t-c_t)||_{\tilde{\omega}_{t,\delta}^2} \rightarrow 0$$ as $t \rightarrow 0$ on   $1\leq |w|_{|V_t} \leq 2$.

By the definition and the estimates of the diffeomorphism $\psi_t$ onto the (non-collapsing) regions $1\leq |w|_{|V_t} \leq 2$:
$$|a_t-c_t| =   |(\log b_t - {\psi_t^{-1}}^{*}\log b_0)+ ({\psi_t^{-1}}^{*}p_0 -p_t)|=\mathcal{O}(t),$$
(and similarly for the derivatives). Hence $\tilde{\omega}_{t,\delta}>0$.

Finally it follows by a simple computation that the global K\"ahler potential of  $\tilde{\omega}_{t,\delta}$ with respect the background metric $\beta_t$ is given by $\phi_t$ defined in the statement of the Proposition.
\qed \\
\end{dimo}

We should observe that Proposition \ref{PG} has an obvious generalization when we consider smoothings, or partial smoothings, with more than one singularity. In fact, it is sufficient to repeat our argument around all singularities which are smoothed out and note that around points which remain singular we can simply glue the orbifold metric $\omega_0$ to  $\tilde{\omega}_{t,\delta}^2 $ on $1\leq |w|_{V_t} \leq 2$ with a cut-off function $\tau_t$ as we did in the proof of the above Proposition. Note that in this last case a pluriharmonic correction is not needed, since $\varphi_0^1$ can be taken  equal to $-\log |\hat{\Omega}|^2_0$.

Then, for sake of notational simplicity, we continue to argue assuming that $Sing(X_0)=\{p\}$.

Now we define the following function (the Ricci potential) for $\sqrt{t}=\delta^2$:
$$f_{\delta,t}:=\log \left(\frac{|s_t|^2_{\beta_t}}{|s_t|^2_{\tilde{\omega}_{t,\delta}}}\right)-\phi_t,$$
 where $\phi_t$ is the global potential of the metric defined in \ref{PG} and $s_t$ a non-vanishing local section of $K_{X_t}^{-1}$. On $V_t$ we take $s_t$ to be equal to $\hat{\Omega}_t$ previously considered.
Then we have the following Proposition:

\begin{prop}\label{E} The function $f_{\delta,t}$ defined above satisfies:
\begin{itemize}
\item $Ric \, \tilde{\omega}_{\delta,t}=\tilde{\omega}_{t,\delta}+\Dt f_{\delta,t}$;
\item Estimates:

 \begin{itemize}\item $||\nabla^{k}_{\tilde{\omega}_{t,\delta}} f_{\delta,t}||_{\tilde{\omega}_{t,\delta}}= \mathcal{O}(\delta^4)$, on $X_t \setminus \{|w| \leq 1 \cap V_t\}$ (recall $\sqrt{t}=\delta^2)$;
\item   $||\nabla^{k}_{\tilde{\omega}_{t,\delta}} f_{\delta,t}||_{\tilde{\omega}_{t,\delta}}= \mathcal{O}(\delta^{4-2\alpha-\frac{k\alpha}{2}})$, on $\{\delta^{\alpha} \leq |w| \leq 2 \delta^{\alpha}\} \cap V_t$ for $\alpha \in [0,\frac{4}{3}]$;
\item  $||\nabla^{k}_{\tilde{\omega}_{t,\delta}} f_{\delta,t}||_{\tilde{\omega}_{t,\delta}}= \mathcal{O}(\delta^{\alpha-\frac{k\alpha}{2}})$, on $\{\delta^{\alpha}\leq  |w| \leq 2 \delta^{\alpha}\} \cap V_t$ for $\alpha \in [\frac{4}{3},2].$
\end{itemize}
\end{itemize}
\end{prop}

\begin{dimo}
It is immediate to check that $f_{\delta,t}$ defined as above is a Ricci potential, i.e., it satisfies the equation $Ric \, \tilde{\omega}_{t,\delta}=\tilde{\omega}_{t,\delta}+\Dt f_{\delta,t}$.

Now we compute how $f_{\delta,t}$ behaves as we approach the singular fiber:

\begin{itemize}
\item Region $\delta^2 \leq |w|_{|V_t} < \delta^{\frac{4}{3}}$. It follows by definition that
$$f_{\delta,t}=\log b_t -\log|\hat{\Omega}_t|^2_{\eta_{\delta,t}}-\varphi^2_\delta+p_t-\log b_t.$$
Since $\eta_{\delta,t}$ is Ricci Flat, $\log|\hat{\Omega}_t|^2_{\eta_{\delta,t}}=0$. Recalling that $|\nabla_{\eta_\delta}^{k}(\varphi^2_\delta- |w|)|_{\eta_\delta} \leq C \delta^4 |w|^{-1-\frac{k}{2}}$ and $|\nabla_{\eta_\delta}^{k} p_t|_{\eta_{\delta,t}}\leq C |w|^{1-\frac{k}{2}}$, we find that
$$ ||\nabla^{k}_{\tilde{\omega}_{t,\delta}} f_{\delta,t}||_{\tilde{\omega}_{t,\delta}}= \mathcal{O}(\delta^{\alpha-\frac{k\alpha}{2}}),$$
on $\{\delta^{\alpha}\leq  |w| \leq 2 \delta^{\alpha}\} \cap V_t$ for $\alpha \in [\frac{4}{3},2]$.

\item Region $\delta^{\frac{4}{3}} \leq |w|_{|V_t} \leq 2 \delta^{\frac{4}{3}}$ (the glueing region). As before we have that
$$f_{\delta,t}= -\log|\hat{\Omega}_t|^2_{\tilde{\omega}_{\delta,t}}-\chi_\delta \varphi_\delta^1-(1-\chi_\delta)\varphi_\delta^2 +p_t.$$

It follows by  $||\nabla_{\delta}^k (\tilde{\omega}_{t,\delta}^1-\eta_{\delta,t})||_{\eta_{\delta,t}} = \mathcal{O}(\delta^{\frac{4-2k}{3}})$ that $$ \log|\hat{\Omega}_t|^2_{\tilde{\omega}_{\delta,t}}= \mathcal{O}(\delta^{\frac{4}{3}}),$$
Hence $$ ||\nabla^{k}_{\tilde{\omega}_{t,\delta}} f_{\delta,t}||_{\tilde{\omega}_{t,\delta}}= \mathcal{O}(\delta^{\frac{4-2k}{3}}).$$

\item  Region $ \delta^{\frac{4}{3}} < |w|_{|V_t} \leq 1$. We have that
$$f_{\delta,t}= -\log|\hat{\Omega}_t|^2_{\tilde{\omega}_{\delta,t}}-\varphi_\delta^1 +p_t.$$
Recalling that the diffeomorphism between the smoothings and the singular fiber is $\psi_t: w_i  \longmapsto w_i+\frac{t}{2|w|^2} \bar{w}_i$ and the definition of $p_t$, it follows that
$$||\nabla^k_\delta ( -\log|\hat{\Omega}_t|^2_{\tilde{\omega}_{\delta,t}}-\varphi_\delta^1 +p_t)||_{\eta_{\delta,t}}\leq Ct|w|^{-2-\frac{k}{2}}.$$
Thus, since $\sqrt{t}=\delta^2$,  for $|w|=\delta^{\alpha}$ with $0 \leq \alpha < \frac{4}{3}$,
$$||\nabla^{k}_{\tilde{\omega}_{t,\delta}} f_{\delta,t}||_{\tilde{\omega}_{t,\delta}}= \mathcal{O}(\delta^{4-2\alpha-\frac{k\alpha}{2}}).$$

\item Region $X_t \setminus \{ |w|_{|V_t} \leq 1\}$. By definition of the diffeomorphism $\psi_t$, we can assume we have complex structures $J_t$ on $X_t$ and $J_0$ on $X_0$  satisfying on  $X_t \setminus \{ |w|_{|V_t} \leq 1\}$, $||\nabla^k_\delta( J_t -{\psi_t^{-1}}^{*}J_0)||_{\tilde{\omega}_{t,\delta}}=\mathcal{O}(t)$ for all $k$. Then,  since $\sqrt{t}=\delta^2$,
$$||\nabla^{k}_{\tilde{\omega}_{t,\delta}} f_{\delta,t}||_{\tilde{\omega}_{t,\delta}}= \mathcal{O}(\delta^4),$$
for all $k$.

\end{itemize}
\qed \\
\end{dimo}

\section{Deformation to a genuine solution}

Let $\tilde{\omega}_{t,\delta}$ be the approximate KE metric on $X_t$ constructed before, where $t$ and $\delta$ satisfy the relation $\sqrt{t}=\delta^2$. Let us recall the basic properties of this metric needed in this section:
\begin{itemize}
\item $\tilde{\omega}_{t,\delta}=\beta_t+\Dt \phi_t \in c_1(X_t)$ (and $\tilde{\omega}_{t,\delta}$ is of orbifold type if some singularities are kept under partial-smoothing) ;
\item ${\tilde{\omega}}_{t,\delta}$ restricted to $ \{ |w| \leq \delta^{\frac{4}{3}} \} \cap V_t$ is equal to  $\eta_{t,\delta}$ (Eguchi Hanson metric) with vanishing cycle  satisfying  $\mbox{Diam}_{\tilde{\omega}_{t,\delta}}(L_t) \approx \delta$.
\item The Ricci potential $f_\delta$ satisfies:
\begin{itemize}
\item $||\nabla^{k}_{\tilde{\omega}_{t,\delta}} f_\delta||_{\tilde{\omega}_{t,\delta}}= \mathcal{O}(\delta^4)$, on $X_t \setminus \{|w| \leq 1 \cap V_t\}$ (recall $\sqrt{t}=\delta^2)$;
\item   $||\nabla^{k}_{\tilde{\omega}_{t,\delta}} f_\delta||_{\tilde{\omega}_{t,\delta}}= \mathcal{O}(\delta^{4-2\alpha-\frac{k\alpha}{2}})$, on $\{\delta^{\alpha} \leq |w| \leq 2 \delta^{\alpha}\} \cap V_t$ for $\alpha \in [0,\frac{4}{3}]$;
\item  $||\nabla^{k}_{\tilde{\omega}_{t,\delta}} f_\delta||_{\tilde{\omega}_{t,\delta}}= \mathcal{O}(\delta^{\alpha-\frac{k\alpha}{2}})$, on $\{\delta^{\alpha} \leq |w| \leq 2 \delta^{\alpha}\} \cap V_t$ for $\alpha \in [\frac{4}{3},2].$
\end{itemize}
Note that the error is bigger exactly at the gluing region $|w| \approx \delta^{\frac{4}{3}}$, where it behaves like  $\delta^{\frac{4-2k}{3}}$.
\end{itemize}

In order to find a KE metric, we need to solve the following equation:
$$E_{t,\delta}[\varphi]:= \frac{\left(\tilde{\omega}_{t,\delta}+\Dt \varphi\right)^2}{\tilde{\omega}_{t,\delta}^2}-e^{f_\delta -\varphi}=0,$$
where  $\sqrt{t}=\delta^2$ and $\varphi$ is a real valued smooth function on $X_t$. Then $\tilde{\omega}_{t,\delta}+\Dt \varphi>0$ would be the desired KE metric. We claim that it is possible to find a solution, provided that $\delta$  (hence the complex structure parameter $t$) is sufficiently small. 

In order to solve the problem we realize the operator $E_{t,\delta}$ in some weighted-H\"older spaces:

$$E_{t,\delta}: \mathcal{U} \subseteq \mathcal{C}^{2,\gamma}_{\tilde{\omega}_{t,\delta}, \beta}(X_t,\R) \longrightarrow  \mathcal{C}^{0,\gamma}_{\tilde{\omega}_{t,\delta}, \beta-2}(X_t,\R),$$
where $\beta$ is some real number in $(-2,0)$ and $\gamma\in(0,1)$ the H\"older exponent, and $\mathcal{U}$ a suitable neighborhood of the origin.

The Banach space $\mathcal{C}^{k,\gamma}_{\tilde{\omega}_{t,\delta}, \beta}(X_t,\R)$ is for $\beta <0$ defined as follows. As a vector space it is simply $\mathcal{C}^{k,\gamma}(X_t,\R)$. However, instead of using the usual H\"older norm,  we use a weighted norm. First of all we  define a weight function $\rho_t:X_t \rightarrow [\delta,1]$:

\begin{itemize}
\item $\rho_t(|w|_{|X_t})= \delta$, if $|w|_{|X_t}\leq 2 \delta^2$;
\item $\rho_t(|w|_{|X_t})= |w|^{\frac{1}{2}}_{|X_t}$, if $3 \delta^2 \leq|w|_{|X_t}\leq \frac{1}{2}$;
\item $\rho_t(|w|_{|X_t})=1$, if $|w|_{|X_t}\geq 1$;
\item $\rho_t$ is a smoothly increasing interpolation between the above values in the remaining regions. 
\end{itemize}

Then the weighted norm is:

$$||\varphi||_ {\mathcal{C}^{k,\gamma}_{\tilde{\omega}_{t,\delta}, \beta}}:= \sum_{j \leq k} || \rho_{t}^{-(\beta-j)} \nabla_{\tilde{\omega}_{t,\delta}}^j \varphi ||_{L^{\infty}_{\tilde{\omega}_{t,\delta}, \beta}} +  [\varphi]_ {\mathcal{C}^{k,\gamma}_{\tilde{\omega}_{t,\delta}, \beta}},$$
where
$$[\varphi]_ {\mathcal{C}^{k,\gamma}_{\tilde{\omega}_{t,\delta}, \beta}}:= \sup_{p,q | \; d_{t} (p,q) \leq inj_{t} \; p \neq q}\left( \min \{ \rho_t^{-(\beta-j-\gamma)}(p), \rho_t^{-(\beta-j-\gamma)}(q)\}\frac{ ||\nabla_{t}^k \varphi(p) - \nabla_{t}^k \varphi(q)||_{t}}{d_{t}^\gamma(p,q)}\right)$$ 
 (We compute the difference of the two derivatives by  parallel transport along the unique minimal geodesic).

Rewrite $E_{t,\delta}[\phi]$ as:
$$E_{t,\delta}[\varphi]=(1-e^{f_\delta})+\mathcal{D}_{t,\delta}[\varphi]+\mathcal{R}[\varphi].$$
Here $\mathcal{D}_{t,\delta}[\varphi]= \Delta_{\tilde{\omega}_{t,\delta}}\varphi+e^{f_{\delta}}\varphi$ and $\mathcal{R}[\phi]$ contains the non-linearities.

The first input we need is the following (scaled)-Schauder estimate.

\begin{lem} If $\delta$ (hence $t$) is sufficiently small  and $\beta \in (-2,0)$, then  
 $$ ||\varphi||_{\mathcal{C}^{2,\gamma}_{\tilde{\omega}_{t,\delta},\beta}} \leq C \left( ||\varphi||_{{L}^{\infty}_{\tilde{\omega}_{t,\delta}, \beta }}+ ||\mathcal{D}_{t,\delta}[\varphi]||_{\mathcal{C}^{0,\gamma}_{\tilde{\omega}_{t,\delta}, \beta-2}} \right)$$
for all $\varphi \in \mathcal{C}^{2,\gamma}_{\tilde{\omega}_{t,\delta},\beta}(X_t)$ with  a positive constant $C$ independent of $\delta$.
\end{lem}

\begin{dimo}
Take a small $c>0$  such that $B_{\omega_0}(p,c) \subseteq X_0$, with $p\in Sing(X_0)$ is metrically equivalent to the ball of radius $1$ in the flat $\C^2/\Z_2$.

Then, by standard Schauder estimates, we know that the desired estimate holds in the region $\psi_t(X_0 \setminus B_{\omega_0}(p,c)) \subseteq X_t$ for a constant $C=C(c)>0$ fixed. More precisely, we can cover the region with (a finite number) of domains  $D_t^i$ where the geometry is close to the euclidean for all $0\leq t \leq c$. By definition, the weighted H\"older norms are equivalent to the usual H\"older norms. Then the claim follows by the standard Schauder estimates since $\mathcal{D}_{t,\delta}$ is elliptic (with smoothly varying bounded coefficients in $t$). 

It remains to show what happens in the ``collapsing'' region $$X_t \setminus \psi_t(X_0 \setminus B_{\omega_0}(p,c)) \subseteq \{w_0^2+w_1^2+w_2^2=t\}.$$

First of all we pull-back the functions $\varphi_i(w)$ from $w_1^2+w_2^2 + w_3^2=t=\delta^4$ to $z_1^2+z_2^2 +z_3^2=1$ using the map $w_i=z_i \delta^2$. Then define the function
$$\tilde{\varphi}(z):=\delta^{-\beta} \varphi(w(z)),$$
and scale the metric $\tilde{\omega}_{t,\delta}$ so that the diameter of the cycle $L_1$ is equal to a constant  (i.e. consider the metric  $g_{t,\delta}:= \frac{1}{\delta^2} \tilde{\omega}_{t,\delta}$). 

Then our scaled Schauder estimates follow once we prove that
$$ ||\tilde{\varphi}||_{\mathcal{C}^{2,\gamma}_{g_{t,\delta},\beta}(\rho_{g_{t,\delta}} \leq \frac{c}{\delta})} \leq C \left( ||\tilde{\varphi}||_{{L}^{\infty}_{g_{t,\delta}, \beta }(\rho_{g_{t,\delta}} \leq \frac{2c}{\delta})}+ ||\Delta_{g_{t,\delta}}\tilde{\varphi}+\tilde{c}\delta^2\tilde{\varphi} ||_{\mathcal{C}^{0,\gamma}_{g_{t,\delta}, \beta-2} (\rho_{g_{t,\delta}} \leq \frac{2c}{\delta})} \right),$$
 for some constant $C$ independent of $\delta$. Here $c,\tilde{c}$ denote  positive constants and the weighted norm is essentially  given by the sum of seminorms $$[\tilde{\varphi}]_{j,\beta}= \sup || \tilde{\rho}_{g_{t,\delta}}^{-(\beta-j)} \nabla_{g} \tilde{\varphi}||_{g_{t,\delta}},$$
(and similarly for the H\"older seminorm), where $\tilde{\rho}_{g_{t,\delta}}$ is a weight function equal to $1$ on $\{x \in X_t \, | \,\rho_{g_{t,\delta}}(x)\leq 1\}$ ($\rho_{g_{t,\delta}}$ denotes the distance from $L_1$ w.r.t. $g_{t,\delta}$), exactly equal to the distance $\rho_{g_{t,\delta}}$ for  $\{x \in X_t \, |\, \rho_{g_{t,\delta}}(x) \geq 2\}$ and  a smooth increasing interpolating function on the remaining annulus.

On the compact piece  $\{x \in X_t \, | \,\rho_{g_{t,\delta}}(x)\leq R\}$ (for a fixed $R>0$),  the estimate holds by standard Schauder estimates (the Riemannian manifold we are considering is just a compact subset containing $L_1$ of the Eguchi Hanson space. Then we can argue as we did in the region away from the singularities). It remains to show what happens in the regions $R\leq \rho_{g_{t,\delta}}(x) \leq \frac{c}{\delta}$. Of course this region can be covered by balls $B_{g_{\delta},t}(p_i, \frac{r_i}{8})$, where $  R \leq r_i \leq \frac{c}{\delta}$. By scaling the usual Schauder estimates on the euclidean balls of radius $1$ and noting that, since $r_i \leq \frac{c}{\delta}$, the balls $B_{g_{\delta},t}(p_i, \frac{r_i}{8})$ are metrically equivalent to standard balls of radius $\frac{r_i}{8}$ in the euclidean space, we have an estimate of the form
$$|\tilde{\varphi}|_{L^\infty(B_{g_{t,\delta}}(p_i,\frac{r_i}{8}))}+ r_i^{\gamma}[\nabla \tilde{\varphi}]_{C^{0,\gamma}(B_{g_{t,\delta}}(p_i,\frac{r_i}{8}))}+\dots+ r_i^{2+\gamma} [\nabla^2 \tilde{\varphi}]_{C^{0,\gamma}(B_{g_{t,\delta}}(p_i,\frac{r_i}{8}))}$$ 
$$\leq C \left( |\tilde{\varphi}|_{L^\infty(B_{g_{t,\delta}}(p_i,\frac{r_i}{4}))}+ r_i^2 |\Delta_{g_{t,\delta}}\tilde{\varphi}+\tilde{c}\delta^2\tilde{\varphi} |_{L^{\infty}(B_{g_{t,\delta}}(p_i,\frac{r_i}{4}))}+ r_i^{2+\gamma} |\Delta_{g_{t,\delta}}\tilde{\varphi}+\tilde{c}\delta^2\tilde{\varphi} |_{C^{0,\gamma}}  \right)  $$  
for a constant $C$ independent of $t$ (and $\delta$). Then the desired estimate follows by multiplying the above inequality by $r_i^{-\beta}$ and observing that on $B_{g_{t,\delta}}(p_i,\frac{r_i}{4}))$ the weight function is bounded by $\frac{r_i}{2} \leq  \tilde{\rho}_{g_{t,\delta}} \leq 2 r_i$.

\qed \\
\end{dimo}

The following estimate is fundamental:
\begin{prop}\label{I} If $\delta$ (hence $t$) is sufficiently small  and $\beta \in (-2,0)$, then 
$$ ||\varphi||_{\mathcal{C}^{2,\gamma}_{\tilde{\omega}_{t,\delta},\beta}} \leq C \, ||\mathcal{D}_{t,\delta}[\varphi]||_{\mathcal{C}^{0,\gamma}_{\tilde{\omega}_{t,\delta}, \beta-2}} $$
for all $\varphi \in \mathcal{C}^{2,\gamma}_{\tilde{\omega}_{t,\delta},\beta}(X_t)$ with  a positive constant $C$ independent of $\delta$.
\end{prop}

\begin{dimo}
The proof is by contradiction. Assume that the above estimate does not hold. Then there exists a sequence $(\delta_i)$ going to zero and smooth functions $\varphi_i$ on $X_{t_i}$ satisfying $ \infty > ||\varphi_i||_{\mathcal{C}^{2,\gamma}_{i,\beta}}\geq C >0$ and $||\mathcal{D}_{i}[\varphi]||_{\mathcal{C}^{0,\gamma}_{i, \beta-2}}\rightarrow 0$.

Then by the above scaled Schauder estimates we must have that $||\varphi_i||_{L^{\infty}_{\beta}}$ is bounded above from zero. W.l.o.g. we can assume that there is a sequence of points $p_i \in X_{t_i}$ where
$$\rho^{-\beta}_i(p_i)|\varphi_i(p_i)|=1.$$
Now we have three cases depending on how (a subsequence of) $p_i$ converges.

$\mathbf{Case \, 1}$: Assume that $\rho_i(p_i)\geq C> 0$. Then $p_i \rightarrow p_0 \in X_0\setminus \mbox{Sing}(X_0)$. By the upper bound on $||\varphi_i||_{\mathcal{C}^{2,\gamma}_{i,\beta}}$ we can assume that on all compact subsets of $ X_0\setminus \mbox{Sing}(X_0)$ (a subsequence of) $\varphi_i \rightarrow \varphi_\infty$ in $\mathcal{C}^{2,\gamma-\epsilon}$ (Ascoli-Arzela) with $\varphi_\infty(p_0)=c>0$ and $|\rho_0^{-\beta}\varphi_\infty| \leq C$,  and similarly for the derivatives, near the singularity (where $\rho$ denotes the distance from the singularity w.r.t. the KE metric $\omega_0$ on $X_0$).

Moreover we have that by the $C^2$ convergence $$\Delta_0\varphi_\infty+\varphi_\infty=0,$$
on all compact  subsets of $ X_0\setminus \mbox{Sing}(X_0)$  and where $\Delta_0$ denotes the Laplacian w.r.t. the KE metric $\omega_0$. Pulling back locally $\varphi_\infty$ to $\C^2$ using the orbifold covering map and recalling the behavior of $\varphi_\infty$ at the origin, it follows that $\varphi_\infty$ is a weak-solution of $\Delta_0\varphi_\infty+\varphi_\infty=0$ as long as $\beta>-2$. In fact, for all $u \in C^{\infty}_{0}(B(0,R))$
$$\int_{B(0,R) \setminus B(0,\rho)} \varphi_\infty \,(\Delta_0 u + u) dV_{\omega_0}= \int_{B(0,R) \setminus B(0,\rho)} (\Delta_0 \varphi_\infty \, + \varphi_{\infty}) u dV_{\omega_0} + $$  $$+\int_{\partial B(0,\rho)} \varphi_{\infty}\frac{\partial u}{\partial \nu}-u\frac{\partial \varphi_{\infty}}{\partial \nu} d\Sigma, $$ 
where $\frac{\partial}{\partial \nu}$ denotes the normal derivative (at $\partial B(0,\rho)$ w.r.t. the metric $\omega_0$). Then
$$\int_{B(0,R) \setminus B(0,\rho)} \varphi_\infty \,(\Delta_0 u + u) dV_{\omega_0} \leq $$  $$ \leq C_1 \rho^3 \sup_{B(0,\rho)}|\nabla_0 u|\sup_{\partial B(0,\rho)} |\varphi_\infty|+ C_2 \rho^3 \sup_{B(0,\rho)}| u|\sup_{\partial B(0,\rho)} |\nabla_0 \varphi_\infty|\leq $$
$$\leq C\rho^{3}(\rho^{\beta}+\rho^{\beta-1}) \rightarrow 0,$$
as $\rho \rightarrow 0$ if $\beta >-2$. By standard regularity theory of elliptic operators, it follows that the weak solution $\varphi_{\infty}$ is actually (orbifold) smooth.

The above implies by a well-known Bochner identity for KE metric on Fano orbifolds (\cite{T00}) that $(\bar{\partial}\varphi_\infty)^{\sharp}$ must be a holomorphic vector field on $X_0$. Using now the discrete automorphism  hypothesis, we have that $ (\bar{\partial}\varphi_\infty)^{\sharp}=0$. Since $\varphi_{\infty}$ is real valued, it follows that $\varphi_{\infty}$ must be constant, hence identically zero by the equation. However this is in contradiction with  $|\varphi_\infty(p_0)|=c>0$.

Now we investigate the cases when $\rho_i(p_i)\rightarrow 0$.

$\mathbf{Case \, 2}: \frac{\delta_i}{\rho_i(p_i)} \rightarrow C > 0$.

First of all we pull-back the functions $\varphi_i(w)$ from $w_1^2+w_2^2+ w_3^2=t=\delta^4$ to $z_1^2+z_2^2 +z_3^2=1$ using the map $w_i=z_i \delta^2$. Then we define the functions
$$\tilde{\varphi_i}(z):=\delta_i^{-\beta} \varphi_i(w(z)).$$
Moreover, scaling the metric $\tilde{\omega}_{t,\delta}$ so that the diameter of $L_1$ is equal to a constant for all $i$ (i.e. blowing up the metric by $\frac{1}{\delta_i^2}$), we have that:
\begin{itemize}
\item $||\tilde{\varphi_i}||_{C^{2,\gamma}(K)} \leq C$, with respect to the norm induced by the blow up metric, for all compact subsets of $K$ of   $z_1^2+z_2^2 +z_3^2=1$.
\item  $||\tilde{\varphi_i}||_{L^\infty}\leq \frac{C}{1+|z|^{-\frac{\beta}{2}}}$;
\item $|\tilde{\varphi_i}(p_i)|= c>0$ for $p_i$ contained in a compact subset of  $z_1^2+z_2^2 +z_3^2=1$.
\end{itemize}

It  follows by Ascoli-Arzel\'a that (a subsequence of) $\tilde{\varphi_i}\rightarrow \tilde{\varphi}_\infty$ in ${C^{2,\gamma-\epsilon}(K)}$ for all compact subsets, $\lim_{|z|\rightarrow +\infty}   \tilde{\varphi}_\infty(z) = 0$ and $ |\tilde{\varphi}_\infty (p_0)|=c>0$ for some point in $z_1^2+z_2^2 +z_3^2=1$ at finite distance from $L_1$.

Moreover the hypothesis  $||\mathcal{D}_{i}[\varphi_i]||_{\mathcal{C}^{0,\gamma}_{i, \beta-2}}\rightarrow 0$ implies that
$$||\Delta_{\frac{1}{\delta_i^2}\tilde{\omega}_{t,\delta}} \tilde{\varphi_i}+\delta_i^2 e^{\tilde{f}_{\delta_i}}\tilde{\varphi_i}||_{L^\infty(K)}\rightarrow 0$$
for all compact $K$ of $z_1^2+z_2^2 +z_3^2=1$. Hence
$$||\Delta_{\frac{1}{\delta_i^2}\tilde{\omega}_{t,\delta}} \tilde{\varphi_i}||_{L^\infty(K)} \leq ||\Delta_{\frac{1}{\delta_i^2}\tilde{\omega}_{t,\delta}} \tilde{\varphi_i}+\delta_i^2 e^{\tilde{f}_{\delta_i}}\tilde{\varphi_i}||_{L^\infty(K)} +\delta_i^2 e^{\tilde{f}_{\delta_i}}||\tilde{\varphi_i}||_{L^\infty(K)}\rightarrow 0.$$
Since blowing-up $\tilde{\omega}_{t,\delta}$ gives in the limit the Eguchi-Hanson Ricci flat metric $\eta_1$, we obtain
 $$\Delta_{\eta_1}  \tilde{\varphi}_\infty=0,$$
that is $ \tilde{\varphi}_\infty$ is an harmonic function on $z_1^2+z_2^2 +z_3^2=1$ equipped with the Eguchi-Hanson metric. If $\beta<0$, $\lim_{|z| \rightarrow +\infty} \tilde{\varphi}_\infty=0$. Then, by the maximum principle, $\tilde{\varphi}_\infty$ must be identically zero, which is again in contradiction with 
$ |\tilde{\varphi}_\infty (p_0)|=c>0$.

$\mathbf{Case \, 3}: \frac{\delta_i}{\rho_i(p_i)} \rightarrow 0$ (and  $\rho_i(p_i)\rightarrow 0$).

We consider the function $\tilde{\varphi_i}$ defined as  we did in the previous case. By hypothesis it follows that $|\tilde{\varphi_i}|(z(p_i))=c>0$ for a sequence $|z(p_i)|=:R_i^2 \rightarrow +\infty$, where $R_i\delta_i \rightarrow 0$. Blowing-down the metrics by a factor $R_i^{-2}$, the new metrics $\frac{1}{\delta_i R_i} \tilde{\omega}_{t,\delta}$ converge to the flat metric $\eta_0$ on $\C^2 /\Z_2$. Arguing as in the previous section, the functions $\chi_i:= R_i^{-\beta} \tilde{\varphi}_i$ converge (up to a subsequence) to a smooth function $\chi_{\infty}$ on $\C^2 / \Z_2 \setminus \{0\}$ which is harmonic w.r.t. the flat metric. Moreover,
\begin{itemize}
 \item $|\chi_{\infty}(p_0)|=C > 0$, where $p_0$ is a point at distance $1$ from the origin;
  \item $|\chi_{\infty}(p)| \leq C d_{\eta_0}(0,p)^{\beta}$.
\end{itemize}
 The pull-back to $\C^2$ of $\chi_\infty$ is an harmonic function which goes to zero at infinity and which is less singular than the Green's function if $\beta \in(-2,0)$.
Thus it extends smoothly to all $\C^2$. Finally, the maximum principle implies that it must be identically zero. However this is in contradiction with $|\chi_{\infty}(p_0)|=C > 0$.

\qed \\
\end{dimo}

\begin{cor} \label{I2} If $\delta$ (hence $t$) is sufficiently small and $\beta \in (-2,0)$, $$\mathcal{D}_{t,\delta}: \mathcal{C}^{2,\gamma}_{\tilde{\omega}_{t,\delta}, \beta}(X_t,\R) \longrightarrow  \mathcal{C}^{0,\gamma}_{\tilde{\omega}_{t,\delta}, \beta-2}(X_t,\R)$$ is invertible with norm of the inverse independent of $\delta$.
\end{cor}
 
\begin{dimo}Observe that in the non-weighted norm the linearized operator is Fredholm of index zero (it is just Laplacian plus 1). Now for a fixed value of $\delta$ (and $\beta\leq0$), the spaces $ \mathcal{C}^{k,\gamma}_{\tilde{\omega}_{t,\delta}, \beta}(X_t,\R)$ and  $\mathcal{C}^{k,\gamma}_{\tilde{\omega}_{t,\delta}}(X_t,\R)$ are equivalent, i.e. they are the same vector space with equivalent norms. This implies that  $\mathcal{D}_{t,\delta}: \mathcal{C}^{2,\gamma}_{\tilde{\omega}_{t,\delta}, \beta}(X_t,\R) \longrightarrow  \mathcal{C}^{0,\gamma}_{\tilde{\omega}_{t,\delta}, \beta-2}(X_t,\R)$ is also Fredholm of index zero.

By the previous estimate we know that $\mathcal{D}_{t,\delta}$  has no kernel for sufficiently small $\delta$. Thus is also surjective (being of index zero). The fact that the norm of the inverse is uniform is again a consequence of  Proposition \ref{I}. 
\qed\\
\end{dimo}

 To solve the equation $E_{\delta,t}(\varphi)=0$ we use the following version of the Implicit Function Theorem (compare Lemma 1.3 in \cite{BM11}):

\begin{lem}\label{In} Let $E:X \rightarrow Y$ be a differentiable  map between Banach spaces and let $R(x):=E(x)-E(0)-D_0E(x)$ be the non linearities. Assume there exists some positive  constants $L$, $r_0$ and $C$ such that:
\begin{itemize}
\item $|R(x)-R(y)|_Y  \leq L(|x-y|_X)(|x|_X+|y|_X),$ for all $x,y$ in $B_X(0,r_0)$;
\item $D_0E$ is invertible with norm of the inverse bounded by $C$.
\end{itemize}
If  for an $r< \mbox{min}\{r_0,\frac{1}{2LC}\}$ the initial error $|E(0)|_Y \leq \frac{r}{2C}$, then there exists a unique solution of the equation $E(x)=0$ in $B_X(0,r)$.

\end{lem} 

\begin{dimo}

Fix   $r<\mbox{min}(r_0,\frac{1}{2LC})$ so that  $|E(0)|_Y \leq \frac{r}{2C}$. Finding a solution of the equation $E(x)=0$ is equivalent to find a fixed point of the map:
$$M(x):=D_0E^{-1}(-E(0)-R(x)).$$
 $M$ sends $B_X(0,r)$ into itself: since $|x|_X \leq r<\mbox{min}(r_0,\frac{1}{2LC})$,
$$|M(x)|_X \leq C\left(|E(0)|_Y+|R(x)|_Y\right) \leq C\left(\frac{r}{2C}+L|x|_X^2\right)\leq r.$$
Moreover,  $|M(x)-M(y)|_{X} < |x-y|_X$, for all $x,y$ in $B_X(0,r)$:
$$|M(x)-M(y)|_X\leq C|R(x)-R(y)| \leq CL(|x|+|y|)|x-y|\leq 2CLr |x-y| < |x-y|.$$
Then the result follows immediately by the Banach-Caccioppoli contraction Theorem.
\qed \\
\end{dimo}

In order to apply the above Lemma to our equation, we need to take a closer look at the initial error and at the non linearities.

\begin{lem}
The initial error $E_{\delta,t}(0)=1-e^{f_d}$  is estimated in the weighted norms as 
$$||E_{\delta,t}(0)||_{\mathcal{C}^{0,\gamma}_{\tilde{\omega}_{t,\delta}, \beta-2}(X_t)}=\mathcal{O}(\delta^{\frac{8-2\beta}{3}}),$$
for $\beta \in(-2,0)$.
\end{lem}
\begin{dimo}
We show how to estimate the norm $||E_{\delta,t}(0)||_{L^{\infty}_{\tilde{\omega}_{t,\delta}, \beta-2}(\delta^{\frac{4}{3}}\leq |w|_{|V_t}\leq 1)} $. The other estimates are similar.

 First of all note that at first order the error term is simply given by  the Ricci potential $f_\delta$ (which we know to be small in the point-wise norm by Proposition \ref{E}). Then according to the proof of Proposition \ref{E} and he definition of the weighted norm, we have
$$||E_{\delta,t}(0)||_{L^{\infty}_{\tilde{\omega}_{t,\delta}, \beta-2}(\delta^{\frac{4}{3}}\leq |w|_{|V_t}\leq 1)} \leq C \sup \{\rho_t^{-(\beta-2)}|f_\delta|\} \leq C\delta^4 \sup\{|w|^{-\frac{\beta}{2}-1}\},$$ 
where $\sup$ is w.r.t. points in the region $\delta^{\frac{4}{3}}\leq |w|_{|V_t}\leq 1$. Since $\beta>-2$ the above quantity attains its maximum when $|w| \approx \delta^{\frac{4}{3}}$ (i.e., at the gluing region). Thus $||E_{\delta,t}(0)||_{L^{\infty}_{\tilde{\omega}_{t,\delta}, \beta-2}(\delta^{\frac{4}{3}}\leq |w|_{|V_t}\leq 1)}=\mathcal{O}(\delta^{\frac{8-2\beta}{3}})$ (note that $\delta^4 \leq \delta^{\frac{8-2\beta}{3}}$ for $\delta$ small).  
\qed 
\end{dimo}

By Proposition \ref{I2}, $$||\mathcal{D}_{t,\delta}^{-1}[E_{\delta,t}(0)]||_{\mathcal{C}^{2,\gamma}_{\tilde{\omega}_{t,\delta}, \beta}}=\mathcal{O}(\delta^{\frac{8-2\beta}{3}}).$$
Now  observe that if $||\varphi||_{\mathcal{C}^{2,\gamma}_{\tilde{\omega}_{t,\delta}, \beta}}\leq \mathcal{O}(\delta^{\frac{8-2\beta}{3}})$ then un-weighted norms behave as $||\varphi||_{L^\infty} \leq \mathcal{O}(\delta^{\frac{8+\beta}{3}})$, $||\nabla \varphi||_{L^\infty_\delta} \leq \mathcal{O}(\delta^{\frac{5+\beta}{3}})$ and 
$||\nabla^2 \varphi||_{L^\infty_\delta} \leq \mathcal{O}(\delta^{\frac{2+\beta}{3}})$. In particular these norms go to zero as soon as $\beta > -2$. Thus the  preimage by the linearized operator of the initial error is small in the pointwise norm up to the second derivatives. This is important when we will prove GH convergence of the metrics.

The non linearities are given by the operator
$$\mathcal{R}_{t,\delta} [\varphi]:= \frac{\Dt \varphi \wedge \Dt \varphi}{ \omega_{t,\delta}^2}-e^{f_\delta}\left(\varphi-1+e^{-\varphi}\right).$$

If $||\varphi||_{L^\infty} << 1 $, since $e^{f_d}=\mathcal{O}(1)$ the non-linearities behave as  
$$\mathcal{R}_{t,\delta} [\varphi]=\frac{\Dt \varphi \wedge \Dt \varphi}{ \omega_{t,\delta}^2}-\varphi^2 +\mathcal{O}(\varphi^3),$$
which implies that
$$||\mathcal{R}_{t,\delta}(\varphi_1)-\mathcal{R}_{t,\delta}(\varphi_2)||_{\mathcal{C}^{0,\gamma}_{\tilde{\omega}_{t,\delta}, \beta-2}} \leq C \delta^{(\beta-2)}||\varphi_1-\varphi_2||_{\mathcal{C}^{2,\gamma}_{\tilde{\omega}_{t,\delta}, \beta}}(||\varphi_1||_{\mathcal{C}^{2,\gamma}_{\tilde{\omega}_{t,\delta}, \beta}}+||\varphi_2||_{\mathcal{C}^{2,\gamma}_{\tilde{\omega}_{t,\delta}, \beta}})$$

We are now ready to state and prove the main result:

\begin{prop}\label{M}   If $\delta$ (hence $t$) is sufficiently small and $\beta \in (-2,0)$ then the equation
$$E_{t,\delta}(\varphi_t)=0,$$
admits a (unique) solution  with $||\varphi_t||_{\mathcal{C}^{2,\gamma}_{\tilde{\omega}_{t,\delta}, \beta}} =\mathcal{O}(\delta^{\frac{8-2\beta}{3}})$. 

Moreover $||\nabla^2 \varphi_t||_{L^\infty_\delta} \leq \mathcal{O}(\delta^{\frac{2+\beta}{3}}) \rightarrow 0$ as $\delta \rightarrow 0$.
\end{prop}
\begin{dimo}
We want to apply Lemma \ref{In} to our operators $E_{t,\delta}$. Take $r_0=r_0(\delta)=C_1 \delta^{\frac{8-2\beta}{3}}$. Since $||\varphi||_{L^\infty} \leq   \mathcal{O}(\delta^{\frac{8+\beta}{3}})$, we can use the estimate of the non-linearities $\mathcal{R}_{t,\delta}$ which gives $L=L(\delta)=C_2 \delta^{(\beta-2)}$. In order to apply the Lemma we need that the initial error $E_{t,\delta}(0)$ is much smaller then $\frac{1}{L(\delta)}$, that is
$$\delta^{\frac{8-2\beta}{3}} << \delta^{2-\beta},$$
which is true for $\beta >-2$. Then we can take $r(\delta) \cong r_0(\delta)$ and apply the Lemma. The estimate on the pointwise second derivative follows by the previous observations.

\qed \\
\end{dimo}

Rephrasing the above Proposition, we have that the form
$$\omega_{t,\delta}=\tilde{\omega}_{t,\delta} +\Dt \varphi_t>0$$
is the K\"ahler form of a KE metric provided $\delta$ (hence $t$) is sufficiently small.

Finally we show that the KE metric constructed above converges in the GH topology to the singular metric on the central fiber.

\begin{prop} The KE Del Pezzo surface $(X_t,\omega_{t,\delta})$ converges in the Gromov-Hausdorff sense to the original KE Del Pezzo orbifold $(X_0,\omega_0)$.
\end{prop}

\begin{dimo}

It is evident from the construction that the pre-glued metric $\tilde{\omega}_{t,\delta}$ is  GH close to the orbifold metric $\omega_0$. That is $d_{GH}((X_t,\tilde{\omega}_{t,\delta}),(X_0,\omega_0)) \rightarrow 0$ as $\delta$ (hence $t$) goes to zero.

On the other hand it follows from the implicit function argument that
$$||\omega_{t,\delta}-\tilde{\omega}_{t,\delta}||_{\tilde{\omega}_{t,\delta}}=||\Dt\phi_t||_{L^\infty_{\tilde{\omega}_{t,\delta}}}\leq\mathcal{O}(\delta^\frac{2+\beta}{3}),$$
where the norms are computed in the standard point-wise norm w.r.t. the pre-glued metric. In particular, the identity map on $X_t$ is a $C \delta^{\frac{2+\beta}{6}}$-quasi isometry between the KE and the pre-glued metric. 

Then it follows by Proposition \ref{GHCC} that
$$d_{GH}((X_t,\omega_{t,\delta}),(X_0,\omega_0)) \leq $$  $$d_{GH}((X_t,\omega_{t,\delta}),(X_t,\tilde{\omega}_{t,\delta}))+ d_{GH}((X_t,\tilde{\omega}_{t,\delta}),(X_0,\omega_0))
\rightarrow  0,$$ as $\delta \rightarrow 0$.

\qed
\end{dimo}

\end{chapter}

\begin{chapter}{Compact Moduli spaces of Del Pezzo surfaces}
 
Smooth Del Pezzo surfaces of degree less than or equal to four come in continuous families. An important Theorem of G. Tian \cite{T90} states that all such smooth Del Pezzo surfaces admit  KE metrics. It is then natural to try to understand the  global geometry of the set of KE Del Pezzo surfaces of fixed degree and  its metric compactification.

 The case of Del Pezzo surfaces of low degree can be seen as the first interesting case where one can try to verify the Main Conjecture \ref{MC} stated in the first Chapter. The only situation where the Main Conjecture is known to be true is the case of degree four Del Pezzo surfaces, essentially due to the work of T. Mabuchi and S. Mukai \cite{MM90}. 

In this Chapter we review the case of degree four Del Pezzo surfaces recalling the main ideas in Mabuchi-Mukai  and extending a little their results by showing that the Main Conjecture completely holds in this case. Then we study the geometry of the GH compactification.

Next we investigate the case of degree three Del Pezzos (i.e., cubic surfaces). Even if we are not able to find a complete answer to the compactification problem, we show how to possibly solve it thanks to our results obtained in the previous Chapter on the deformations of nodal Del Pezzo surfaces under some additional hypothesis (which we believe should be always satisfied). 

Finally we briefly consider the case of degree one and two Del Pezzo surfaces by discussing an example, related to an old Conjecture of G. Tian, of smoothable KE Del Pezzo surfaces with log-terminal but not canonical singularities.    

A small note on the relations between the above moduli spaces and K-polystability is added at the end of the Chapter.

\section{Del Pezzo Quartics}

Let $X_4$ be a smooth Del Pezzo of degree $4$. It is classically well-known that the anticanonical divisor is very ample. The anticanonical embedding realizes a degree 4 Del Pezzo surface as a complete intersection $X_4=Q_1 \cap Q_2$ of two quadrics $Q_i$ in $\p^4$. 

The moduli problems of degree $4$ Del Pezzo surfaces fits naturally inside a classical GIT picture. Given two quadrics $Q_1$ and $Q_2$, one can consider the pencil $\lambda Q_1+ \mu Q_2$ (i.e., a (projective) line in the space of quadric polynomials on $\p^4$). Thus we may identify the space parameterizing  intersections of two quadrics with $Gr(2, H^0(\p^4, \mathcal{O}(2)))$. Taking the Pl\"ucker embedding, we consider the following GIT picture:
$$SL(5,\C) \curvearrowright Gr(2, H^0(\p^4, \mathcal{O}(2))) \hookrightarrow \p(\Lambda^2 H^0(\p^4, \mathcal{O}(2))).$$
(i.e. we are considering the natural GIT stability on the $2$-Hilbert point of $X_4$). Thus it is natural to consider 
$$\overline{\mathcal{M}_4}^{ALG}:= Gr(2, H^0(\p^4, \mathcal{O}(2))) // SL(5,\C),$$
as a possible candidate to be a good algebraic compactification of the space of smooth Del Pezzo degree $4$ (if all smooth degree $4$ Del Pezzo surfaces are indeed stable).

In order to study the above GIT picture, it is useful to consider the natural map of sets: 
$$\begin{array}{cccc} disc:& Gr(2, H^0(\p^4, \mathcal{O}(2))) & \longrightarrow & \{ \mbox{Binary quintics} \}/SL(2,\C) \\
    & [\lambda Q_1+ \mu Q_2] & \longrightarrow & det(\lambda Q_1+ \mu Q_2).
  \end{array}
$$

The central result is the following \cite{MM90}:

\begin{thm}[T. Mabuchi-S. Mukai] The map $\mbox{disc}$ descends to the GIT quotient and it induces an isomorphism between $\overline{\mathcal{M}_4}^{ALG}$ and the moduli space of binary quintics $\overline{\mathcal{M}_{bin, 5}}^{GIT}:= \p H^0(\p^1, \mathcal{O}(5)) // SL(2,\C).$

 Every polystable intersection of two quadrics can be ``diagonalized'', i.e., it can be expressed as
$$ X: \,\begin{cases} x_0^2+x_1^2+x_2^2+x_3^2+x_4^2=0\\
 \lambda_0 x_0^2+ \lambda_1x_1^2+\lambda_2x_2^2+\lambda_3x_3^2+\lambda_4x_4^2=0 \end{cases}, $$
with \emph{at most} pairs of equal $\lambda_i$. (Conversely, every intersection of quadrics that admits the above representation is polystable).

Moreover $X=Q_1 \cap Q_2$ is stable iff all $\lambda_i$ in the above representation are distinct iff $X$ is smooth. Note that all polystable intersections of two quadrics admit at most $A_1$ (nodal) singularities.
\end{thm}

It is known that all smooth Del Pezzo surfaces of degree $4$ admit a KE metric. This follows, for example, by computing the alpha invariant \cite{MM90} (and using the fact that all smooth Del Pezzo surfaces of degree $4$ can be simultaneously diagonalized, as proven in M. Reid PhD Thesis). 

Let $(X_{\lambda^i},\omega_i)$ be a sequence of smooth KE Del Pezzo quartics. Suppose that $X_{\lambda^i}\rightarrow X_\mu$ (in the Hilbert scheme), where $X_\mu$ is polystable intersection of quadrics (i.e.,  $\lambda_j^i \rightarrow \mu_j$, with no more that two equal $\mu_j$). Using the fact that we may assume that the KE metric is invariant w.r.t. the finite group of automorphisms generated by $x_i \mapsto -x_i$ and $x_j \mapsto x_j$ for $j\neq i$, it is easy to see that the (ri)-embedding of $X_i$ using $L^2$ KE-orthonormal section (Tian's gauge) are still in a simultaneously diagonalized form, i.e.,
$$ T(X_i): \,\begin{cases} N_0^ix_0^2+N_1^i x_1^2+ N_2^i x_2^2+ N_3^i x_3^2+ N_4^i x_4^2=0\\
 N_0^i\lambda_0^i x_0^2+ N_1^i\lambda_1^ix_1^2+N_2^i\lambda_2^ix_2^2+N_3^i\lambda_3^ix_3^2+N_4^i\lambda_4^ix_4^2=0 \end{cases}, $$
for positive constant $N_j^i$ with $N_j^i \geq N_{j+1}^i$.

Taking $i\rightarrow \infty$, it may happen that $N_j^i$ goes to zero, i.e., that the image of the limit anticanonical map (the flat limit) is contained in a degenerate (diagonalized) intersection of two quadrics:
$$W:\sum_{i=0}^{k} z_j^2= \sum_{j=0}^{k} \mu_j z_j^2=0,$$
where $z_j (N_j^{\infty})^{\frac{1}{2}}=x_j$ if $N_j^\infty \neq 0$ and $x_j=z_j$ otherwise and $k\leq 4$. The crucial argument in Mabuchi Mukai consists in showing that $k=4$ and that the corresponding limit anticanonical map is a biholomorphism. Since $\mu$ is generic, this proves that every polystable intersection of quadrics admit a KE metric.

On the other hand given any sequence of smooth degree $4$ Del Pezzo $(X_{\lambda^i})$, by compactness of the algebraic moduli space, we may assume that $\lambda_j^i \rightarrow \mu_j$, with no more that two equal $\mu_j$, (i.e., $X_{\lambda^i}$ converge to a polystable intersection). Then the GH limits are all polystable intersections of quadrics.

Recalling that the moduli space of binary quintics is isomorphic to $\p(1,2,3)$, we can state the main result:

\begin{prop}
$$\overline{\mathcal{M}}^{AG}_{4} \cong \overline{\mathcal{M}_{bin, 5}}^{GIT} \cong \p(1,2,3). $$
The set of smooth Del Pezzo quartics is parametrized by the complement of the ample  divisor of equation (in weighted homogeneous coordinates) $$ x^2=128y, $$
see  \cite{DI03} pag. 151 for the explicit choice of invariants. 
Moreover, the natural map $$D:\overline{\mathcal{M}}^{AG}_{4} \cong \p(1,2,3)\longrightarrow \overline{\mathcal{M}}^{GH}_{4},$$ given by associating to each polystable Del Pezzo quartic its KE metric structure, is continuous and generically $2:1$.
\end{prop}
\begin{dimo}
 This follows for example by \cite{DI03}, where invariants of binary quintics are computed. Note that the space of smooth Del Pezzo quartics is quasi-projective. The continuity of the map $D$ is obvious by construction. $D$ is generically $2:1$ since Del Pezzo quartics are not given by the product of two lower dimensional Del Pezzo surfaces according to Theorem \ref{S}.
\qed\\
\end{dimo}

As we have pointed out during the proof of the last Proposition,  two isometric  KE degree four del Pezzo surfaces must be either biholomorphic or complex conjugate. Therefore, it may be interesting to study the action of complex conjugation on the above algebraic moduli space. Not too surprisingly, the above  complex conjugation is induced by an antiholomorphic map on the moduli space:

\begin{prop}
 $$\overline{\mathcal{M}}^{GH}_{4} \cong_{top} \p(1,2,3) / \Z_2,$$
where the $\Z_2$-action is given by the antiholomorphic involution obtained by the ordinary complex conjugation on weighted projective space.

In particular the ramification locus of the map $D$ in the previous Proposition is the real surface $\R \p (1,2,3)\cong_{top} \R \p^2$.

\end{prop}

\begin{dimo}
 First of all let us recall that if $(X=Q_1 \cap Q_2, J)$ is a given intersection of two quadrics in $\p^4$, then its complex conjugate variety $\overline{X}$ given by changing $J$ to $-J$ is simply  $\overline{X}=\overline{Q_1} \cap \overline{Q_2}$, i.e., the complex conjugate of the two quadratic polynomials.

Moreover it is immediate to see that the  above action descend, set theoretically, to the quotient by the $SL(5,\C)$ action. Therefore, restricting our considerations to the set of $2$-Hilbert polystable intersections, we see that the action of complex conjugation is given by:

$$ X: \,\begin{cases} x_0^2+x_1^2+x_2^2+x_3^2+x_4^2=0\\
 \lambda_0 x_0^2+ \lambda_1x_1^2+\lambda_2x_2^2+\lambda_3x_3^2+\lambda_4x_4^2=0 \end{cases} $$
$$\downarrow$$   $$  \overline{X}:\,\begin{cases} x_0^2+x_1^2+x_2^2+x_3^2+x_4^2=0\\
 \overline{\lambda}_0 x_0^2+ \overline{\lambda}_1x_1^2+\overline{\lambda}_2x_2^2+\overline{\lambda}_3x_3^2+\overline{\lambda}_4x_4^2=0 \end{cases}   $$

To each polystable intersection of two quadrics is associated the binary quintics $q(z,w)= \prod_{i=0}^{4}(z-\lambda_i w) \in \C[z,w]$. Thus we see that if $X \mapsto q$, then $\overline{X} \mapsto \overline{q}$. In particular the complex conjugation coincides with the natural complex conjugation on $\p^5$ parametrising binary quintics.

It is a classical well-known result (compare \cite{DI03}) that the polynomials generating the algebra of invariants $I_i$ are \emph{real} and satisfy a real syzygy. Then the statement follows immediately by the previous consideration looking at the isomorphism
$$ \overline{\mathcal{M}}^{AG}_{4}\cong \p^5 //SL(2,\C) \rightarrow \p(1,2,3)$$
given, with a little abuse of notation, by $q \mapsto (I_4(q):I_8(q):I_{12}(q))$.

\qed
\end{dimo}

It is maybe interesting to note that an intersection of quadrics $X=Q_1 \cap Q_2$,  with $Q_1$ and $Q_2$ simultaneously diagonalized as in the above proofs, can be biholomorphic to its complex conjugate $\overline{X}$ even if the coefficients of $Q_2$ are not real: in fact, if the coefficients of $Q_2$ appear in complex conjugate pairs, e.g., $(\lambda_0,\lambda_1=\overline{\lambda}_0, \mbox{others real})$, then  linear maps of the type $x_0 \leftrightarrow x_1$ give  biholomorhisms between $X$ and $\overline{X}$. 

Thanks to the previous results, it is easy to compute some topological invariants of the GH compactification.

\begin{cor}
$$H^{k}(\overline{\mathcal{M}}^{GH}_{4};\Q)=\begin{cases} \Q& k=0,4;\\
0 & \mbox{otherwise}. \end{cases}.
$$
In particular the Euler characteristic $e(\overline{\mathcal{M}}^{GH}_{4})=2$.
 
\end{cor}

\begin{dimo} Consider the ramified covering map:
$$\begin{array}{cccc}
   p:& \p^2 &\longrightarrow & \p(1,2,3) \\
  & (x:y:z) & \longmapsto & (x:y^2:z^3) \end{array}
$$
It is well known \cite{K73} that the induced map on rational cohomology 
$$p^{*}: H^{*}\left(\p(1,2,3);\Q \right) \longrightarrow H^{*} \left(\p^2;\Q \right)$$ is an isomorphism, which is also manifestly equivariant w.r.t. the $\Z_2$-action induced by  complex conjugations on $\p^2$ and $\p(1,2,3)$.

Now we recall (compare \cite{S83}) that the transfer homomorphism in cohomology induces the following isomorphism:
$$\pi^{*}: H^{*}\left(\p(1,2,3)/\Z_2;\Q\right) \longrightarrow H^{*}\left(\p(1,2,3);\Q\right)^{\Z_2},$$
where $H^{*}\left(\p(1,2,3);\Q\right)^{\Z_2},$ denotes the $\Z_2$-invariant part of the cohomology. 

Thus the map

 $$p^{*} \pi^{*}:H^{*}\left( \overline{\mathcal{M}}^{GH}_{4};\Q \right)\longrightarrow H^{*} \left(\p^2;\Q \right)^{\Z_2},$$
is an isomorphism. It is straightforward to see that the $\Z_2$-action is trivial on the $0$ and $4$ cohomology and  changes the hyperplane class $[\omega]$ to $-[\omega]$ inside $H^2$. Then the claim follows trivially.\qed
\end{dimo}

We conclude this section by describing some special points in the moduli space and the corresponding intersection of quadrics (the proof is just a trivial computation of the value of the invariants as in \cite{DI03}, so we omit it).

\begin{prop} Under the identification $\overline{\mathcal{M}}^{AG}_{4} \cong \p(1,2,3), $ we have that
\begin{itemize}
 \item $p_1=[1:0:0]$ is a \emph{smooth} point corresponding to the \emph{smooth} intersection of quadrics given by: $Q_1=Id$ and $Q_2=\mbox{diag}(1,\mu, \mu^2,\mu^3,\mu^4)$ with $\mu^5=1$ primitive;
\item $p_2=[0:1:0]$ is a \emph{singular} point of type $A_1$ corresponding to the \emph{smooth} intersection of quadrics given by: $Q_1=Id$ and $Q_2=\mbox{diag}(0,1, \mu^1,\mu^2,\mu^3)$ with $\mu^4=1$ primitive;
\item $p_3=[0:0:1]$ is a \emph{singular} point of type $A_2$ corresponding to the \emph{singular} intersection of quadrics given by: $Q_1=Id$ and $Q_2=\mbox{diag}(0,0,1, \mu^1,\mu^2)$ with $\mu^3=1$ primitive;
\item $p_4=[16:2:6]$ is a \emph{smooth} point corresponding to the toric \emph{singular} intersection of quadrics $x_0x_1=x_2^2=x_3x_4$.
\end{itemize}
The ramification locus $\R \p (1,2,3)$ goes through all the above special points and intersects the divisor parameterizing singular intersections of quardrics $x^2=128y$ non transversely along an $S^1$. Thus the quotient $\p(1,2,3)/ \Z_2 \cong \overline{\mathcal{M}}^{GH}_{4}$ fails to be a topological manifold only 
at two points (the images of $p_2$ and $p_3$).
\end{prop}

Finally it is interesting to note that there are singular points in the moduli corresponding to a smooth variety and viceversa smooth points in the moduli corresponding to singular varieties.

\section{Del Pezzo Cubics}

We begin by recalling some known results on cubic surfaces and their moduli space.

First of all let us note that  two (orbifold) cubic  surfaces $X,Y \subseteq \p^3$ are abstractly biholomorphic if and only if there is a matrix in $SL(4,\C)$ such that $X=g.Y$ (where the action is the natural one induced on the space of cubic polynomials $\p^{19} = \p(H^0(\p^3,\mathcal{O}(3)))$). This simply follows by the very ampleness of the anticanonical divisor.

Then, denoting with $\mathcal{M}_3$  the set of biholomorphism classes of \textit{smooth} cubics, we have
$$\mathcal{M}_3 \cong \left(\{ \mbox{Smooth cubics polynomials}\} \subseteq \p^{19}\right) / SL(4,\C).$$
In order to check that the above space satisfies the Hausdorff properties, it  is important to study the GIT quotient for the  $SL(4,\C)$ action on the space of cubics (note that the associated notion of stability coincides with Chow-stability). The following results, which follows by an easy but long application of the Hilbert-Mumford criterion and which was already contained in Hilbert's Doctoral Dissertation of 1885, describes completely the GIT quotient.

\begin{prop}[D. Hilbert]\label{Hilb} A cubic surface $X$ is
\begin{itemize}
  \item  Chow stable if and only if it has at most $A_1$ singularities.
  \item Chow semistable if and only if has at most $A_1$ or $A_2$ singularities.
\end{itemize}
Moreover there exists only one Chow strictly polystable orbit corresponding to the singular cubic defined by $T:=\{xyz=t^3\}$.
\end{prop}

Thus 
$$\mathcal{M}_3  \subseteq {\overline{\mathcal{M}_3}}^{ALG} := \p^{19} // SL(4,\C),$$
where ${\overline{\mathcal{M}_3}}^{ALG}$ is a connected projective compact algebraic variety. Moreover the boundary $$\partial{\overline{\mathcal{M}_3}}^{ALG}:={\overline{\mathcal{M}_3}}^{ALG} \setminus \mathcal{M}_3$$ consists of two strata 
$$\partial{\overline{\mathcal{M}_3}}^{ALG}= \mathcal{SS}_3 \cup \{T\},$$
where $\mathcal{SS}_3$ denotes the set of singular stable cubics, and $T$ the (orbit) of the toric cubic $xyz=t^3$ (with $3A_2$ singularities). Moreover we can stratify $SS_3$ by the number of $A_1$ singularities
$$\mathcal{SS}_3= \mathcal{SS}_3^1 \cup \mathcal{SS}_3^2 \cup \mathcal{SS}_3^3 \cup \{C\},$$
where $\mathcal{SS}_3^j$ denotes the set of cubics with $jA_1$ singularities and $C$ the (orbit) of the Cayley's cubic $xyz+yzt+ztx+txy=0$ with $4A_1$ singularities.

Next we focus on the existence problem of  KE metrics on cubic surfaces. Denote with $\mathcal{U}^{KE}$ the set defined by
$$\mathcal{U}^{KE}:=\{[X] \in  {\overline{\mathcal{M}_3}}^{ALG} \, | \, \mbox{ the polystable rep. admits a KE orbi. metric }\}.$$

\begin{thm}[G. Tian \cite{T90}, B. Wang \cite{W10}]
$$\mathcal{M}_3 \cup \mathcal{SS}_3^1 \cup \{C\} \cup \{T\} \subseteq \mathcal{U}^{KE}.$$ 	 
Moreover, Chow unstable cubics cannot admit KE metrics. (This last point follows by Tian's computation of the generalized Futaki invariant \cite{DT92}).
\end{thm}

Note that $C$ and $T$ admit obvious KE orbifold metrics since they are quotients of a KE Del Pezzo of degree $6$ and of $\p^2$ respectively.

Some remarks on the above results are needed. The proofs of these results \textit{don't imply} that the singular KE cubics in $\mathcal{SS}_3^1 \cup \{C\} \cup \{T\}$ are GH limits of smooth cubics. Moreover it is not a priori obvious that Chow semistable but not polystable cubics do not admit KE metrics. 

Using our result Theorem \ref{MT} we can extend the knowledge on the existence and degenerations of KE metric on cubic surfaces:

\begin{prop} $$\mathcal{SS}_3^j \cap \mathcal{U}^{KE} \neq \emptyset $$ for $j=2,3$, i.e., there are some KE cubics with $2$ or $3$ nodal singularities.  Moreover all KE cubics with only nodal singularities appear in the ``boundary'' of the KE compactified moduli space $ \overline{\mathcal{M}_3}^{GH}$.
\end{prop}

With the above results on algebraic and metric compactifications in mind, we are now able to state the precise statement of the Main Conjecture \ref{MC}  in the case of cubic surfaces:

\begin{conj}[Degenerations of KE cubic surfaces]\label{C} There exists a $2:1$ ramified covering mapping
$$D:{\overline{\mathcal{M}_3}}^{ALG} \longrightarrow \overline{\mathcal{M}_3}^{GH},$$
where ${\overline{\mathcal{M}_3}}^{ALG} := \p^{19} // SL(4,\C)$. More precisely, every (Chow)-polystable cubic admits a KE (orbifold) metric and the continuous map $D$ is given by associating to each polystable Del Pezzo cubic its induced KE metric structure.
In particular GH limits of cubics have at most $A_1$ or exactly $3A_2$ singularities.

Beside Chow-polystable cubics, no other (normal) cubic surface admits a KE metric. 
\end{conj}

The rest of the section is dedicated to prove the above Conjecture assuming our Deformation Conjecture \ref{KEDC} and to the study of the geometry of ${\overline{\mathcal{M}_3}}^{ALG}$. We remark that the argument given below can be used to reprove the existence of KE metrics on smooth cubic surfaces once we know that the set of Chow-polystable cubics is not empty. We may point out that similar techniques can be used to reprove the case of Del Pezzo quartics.

For the reader's convenience, we divide the argument in four subsections:

\begin{itemize}
 \item GH degenerations of cubics in Tian's gauge;
 \item Deformations of KE orbifold cubics;
 \item Proof of the Conjecture under some assumptions;
 \item The geometry of the compactified moduli space.
\end{itemize}

\subsection*{GH degenerations of cubics in Tian's gauge}

We begin by recalling some basic facts on Tian's gauge (compare \cite{T90}) and limit anticanonical maps discussed in the first Chapter.

Let $(X_i, \omega_i)$ be a sequence of smooth KE cubics. For each $i$ choose an $L^2$-KE orthonormal basis of $H^{0}(K_{X_i}^{-1})$ and denote with $T_i$ the associated holomorphic embedding of $X_i$ into the fixed $\p^3$. We say that $T_i(X_i)$ is a cubic in Tian's gauge. Being  defined only up to the compact group $SU(4)$ and by the compactness of the Hilbert scheme (in this case simply $\p^{19}$), we can assume that there is a subsequence  of the $T_i(X_i)$ converging to a possible very singular (non-reduced, reducible and with bad singularities) cubic $W$, which we call the \textit{flat limit}.

On the other hand, we know that, again up to subsequences, $X_i$ converge to a KE orbifold $G$ in the Cheeger-Gromov sense, i.e. GH and smoothly up to diffeomorphisms on compact sets away from the singularities. We call this limit $G$ the\textit{ GH limit}. Moreover it is possible to show (Tian) that an  $L^2$-KE orthonormal basis of $H^{0}(K_{X_i}^{-1})$ converges smoothly to an orbifold basis of $H^{0}(K_{G}^{-1})$ which is also $L^2$-KE orthonormal w.r.t. the limit orbifold KE metric $\omega_{\infty}$.

In conclusion, we can  w.l.o.g assume that $(X_i,\omega_i)$ GH converges to an orbifold $G$ and to the flat limit $W$, once embedded using Tian's gauge. Moreover, by construction, the biholomorphic embedding $T_i$ converges to a \textit{limit anticanonical meromorphic map} $$T_\infty:G \dashrightarrow W$$ given by $4$ anticanonical sections.  We can also assume that
\begin{itemize}
	\item The indeterminacy locus, i.e. the set of points in $G$ where the function fails to be  holomorphic (or can be extended to one), which in general is a just a subset of the base locus (a priori the linear system can have some fixed part), consists of finitely many isolated points.
	\item The image  $T_\infty(G) $ is of dimension $\geq 1$ (by the Open map Theorem).
\end{itemize}

Using the estimate of Theorem \ref{NDP}, it is easy to prove the following result. 

\begin{prop} \label{HC} The map $T_\infty$ is a biholomorphism between $G$ and $W$. In particular $G \cong W$ is an irreducible, Chow semistable, orbifold cubic which admits a KE metric.
\end{prop}

\begin{dimo}
By general theory (compare Chapter two), it is known that a GH limit of a sequence of smooth KE cubic surfaces must be a Del Pezzo of degree three with $\Q$-Gorenstein smoothable singularities. Combining the estimate of Theorem \ref{NDP} with the Koll\'ar, Shepherd-Barron  classification of $\Q$-Gorenstein smoothable singularities (see Theorem \ref{TS}), it follows that $G$ must have at most canonical singularities of type $A_1$ or $A_2$. In particular the anticanonical divisor is very ample, i.e. the limit anticanonical map is actually an embedding (\cite{D80}). Finally, $G\cong W\subseteq \p^3$ must be Chow semistable by Proposition \ref{Hilb}. \qed 
\end{dimo}

 \subsection*{Deformations of KE orbifold cubics}

 By our Theorem on nodal Del Pezzo, we know that generic deformations of KE cubic surfaces with $A_1$ singularities admit KE metrics. If we assume our deformation Conjecture \ref{KEDC}, then we know that  the set of KE polystable cubics with at most nodal singularities is open in the analytic topology of $\overline{\mathcal{M}_3}^{ALG}$.

 Now let $T:=\{xyz=t^3\}$ be the unique KE Chow polystable cubic in $\overline{\mathcal{M}_3}^{ALG}$ with non trivial $Aut_0(T)$. We would like to apply our Conjecture \ref{KEDC} to this particular example. 
Consider the $A_2$ singularity of $T$ at the point $[0,0,1,0]$ of $\p^3_{x,y,z,t}$. All the cubics obtained by (partial)-smoothings of this singularity are given the  equations:
 $$xyz=t(t-a_3z)(t-b_3z).$$
We can think the pair $(a_3,b_3)$ as coordinates for $\mathcal{E}xt_{0}^{1}(\Omega^1_T,\mathcal{O}_T)$ at $p$. Thus we have the following:
 \begin{itemize}
 \item If $a_3=0$ or $b_3=0$ or $a_3=b_3$ (but not all zero) then the $A_2$ singularity is deformed to an $A_1$ singularity.
 \item If $a_3\neq 0$ and $b_3 \neq 0$ and $a_3 \neq b_3$ then the $A_2$ singularity is smoothed out.
 \end{itemize}

A similar argument applies to the other two $A_2$ singularities. In particular, by Proposition \ref{Manetti}, the versal family of  $\Q$-Gorenstein deformations of $T$ is given by
$$xyz=t[(t-a_1x)(t-b_1x)+(t-a_2y)(t-b_2y)+(t-a_3z)(t-b_3z)].$$
It is then easy to see how  $Aut_0(T)$ acts on the space of versal deformations:
 \begin{lem} Let $T:=\{xyz=t^3\}$. Then
 \begin{itemize}
 \item $Aut_0(T) \cong \left( \begin{array}{ccc}
\lambda_1 & 0 & 0 \\
0 & \lambda_2 & 0 \\
0 & 0 & \lambda_3 \end{array} \right)$ with $\lambda_1\lambda_2\lambda_3=1$, $\lambda_i \in \C^*$; 

It acts on $T$ by $x\mapsto\lambda_1 x$, $y\mapsto\lambda_2 y$, $z\mapsto\lambda_3 z$ and $t \mapsto t$.
 \item $\T^1_T\cong H^0(T,\mathcal{E}xt^1(\Omega^1_T,\mathcal{O}_T)) \cong \C^2_{a_1,b_1} \oplus \C^2_{a_2,b_2} \oplus \C^2_{a_3,b_3}$;
 \item The action $Aut_0(T) \curvearrowright \T^1_T \cong \C^6$  is given by  $\mbox{Diag}\left(\lambda_1,\lambda_1, \lambda_2,\lambda_2,\lambda_3,\lambda_3\right).$
 \end{itemize}
 \end{lem}
 
\begin{dimo}
 The above properties are immediate once observed that $T$ is isomorphic to $\p^2 / \Z_3$ via the map:
$$[z_0:z_1:z_2] \longmapsto [z_0^3:z_1^3:z_2^3:z_0z_1z_2].$$
\qed
\end{dimo}

Then, we can study the corresponding GIT stability notion:
 
 \begin{lem}\label{T} Let $Aut_0(T)\cong \mbox{diag}\left(\lambda_1,\lambda_2,\lambda_3\right)$ with $\lambda_1\lambda_2,\lambda_3=1$ acts on $\T^1_T$ as explained in the previous Lemma and let $v\neq 0$ be a vector in $\T^1_T \cong \C^6$. Then the followings are equivalent
 \begin{enumerate}
 \item $v$ is polystable;
 \item No one of the pairs $(a_i,b_i)$ is equal to zero;
 \item The corresponding deformed cubics have no more $A_2$ singularities. 
 \end{enumerate} 
 \end{lem}

 \begin{dimo} We have already discuss the implications $(2)\Leftrightarrow(3)$. 
 
 $(1)\Rightarrow(2)$. Suppose w.l.o.g that $(a_1,b_1)=0$. Then it is sufficient to consider the action of the 1-ps $\sigma_t:= \mbox{diag}(\frac{1}{t^2}, t, t)$ to  destabilize $v$, i.e. $\mbox{lim}_{t\rightarrow 0} \sigma_t v=0$.
 
 $(2) \Rightarrow (1)$. First of all note that by the nature of the action, the only possibility for an orbit of not being closed is that zero is in its closure. Thus what we need to prove is that under assumption $(2)$ is impossible to converge to zero. W.l.o.g we can assume that $a_1 \neq 0$. If $0$ is in the closure of the orbit of $v$ then we need $|\lambda_1| \rightarrow 0$. But, since $\lambda_1\lambda_2\lambda_3=1$, this would implies that at least one of $|\lambda_i| \rightarrow \infty$, say $\lambda_2$. But condition $(2)$ is telling us that at least one between $a_2$ and $b_2$ is different from zero. If $a_2 \neq 0$ then we would get $a_2 \lambda_2 \rightarrow \infty$ which contradicts our assumption.

 \qed
 \end{dimo}
 
 Finally we can relate the above notion of stability with Chow stability of cubics. By combining Theorem \ref{Hilb} with Lemma \ref{T}, we immediately infer the following:

 \begin{lem} Under the previous notations:
 \begin{itemize}
 \item $v \in \T^1_T$ is polystable $\Leftrightarrow$ the corresponding cubic surfaces is Chow polystable.
 \item $v \in \T^1_T$ is not-polystable  $\Leftrightarrow$ the corresponding cubic surfaces is Chow semistable but not polystable.
 \end{itemize}
 \end{lem}
 
In conclusion, we have:
 
 \begin{prop} \label{Open} If Conjecture \ref{KEDC} holds then
 $\mathcal{U}^{KE}$ is open in the analytic topology of $\overline{\mathcal{M}_3}^{Chow}$.
 \end{prop}

 \subsection*{Proof of the conjecture under some assumptions}
 
We are now going to prove Conjecture \ref{C} assuming that the general picture on deformations of KE metric holds. Note that this assumption is quite natural: the only point where the study of the deformation theory is more complicated, is  at the toric cubics $T$.

First of all we prove that $\mathcal{U}^{KE}$ is closed in $\overline{\mathcal{M}_3}^{Chow}$. Before proceeding, recall that, since all Chow semistable but not polystable cubics can degenerate to the toric cubic $T$,  GH limits of cubics are indeed Chow polystable.
 
 \begin{prop}\label{Closed} If Conjecture \ref{KEDC} holds then
 $\mathcal{U}^{KE}$ is closed in the analytic topology of $\overline{\mathcal{M}_3}^{Chow}$.
 \end{prop}

 \begin{dimo} First of all notice that, being $\overline{\mathcal{M}_3}^{Chow}$ metrizable, we can prove closure using sequences. We will show that $\mathcal{U}^{KE}$ coincides with the closure in the analytic topology of the locus parametrizing smooth KE cubics.
 
 Thus, let us start by defining $\mathcal{U}_{Smooth}^{KE} \subseteq \mathcal{U}^{KE}$ to be the set of smooth KE cubics. We claim that $\mathcal{U}_{Smooth}^{KE}$ is pre-compact in $\mathcal{U}^{KE}$ w.r.t. the (induced) analytic topology of the Chow quotient.
 
 Let $(p_i) \subseteq \mathcal{U}_{Smooth}^{KE}$ be a sequence. All we need to prove is that we can find a subsequence converging to a point $p_0 \in \mathcal{U}^{KE}$. W.l.o.g we can assume we have $p_i=[T(X_i)]$ i.e. $p_i$ is represented by a smooth cubic in Tian's Gauge. Then we know (Proposition \ref{HC}) that up to subsequence $[T(X_{i_j})] \rightarrow [T_{\infty}(G)]$ where $G$ is the GH limit (hence KE) and $T_\infty$ the biholomorphic limit anticanonical map. As we have remarked $T(G)$ must be Chow polystable. Thus $p_{i_j} \rightarrow [T_{\infty}(G)]=:p_0 \in \mathcal{U}^{KE}$, i.e. $$\overline{\mathcal{U}_{Smooth}^{KE}} \subseteq \mathcal{U}^{KE}$$ (here bar denotes the closure in the Analytic topology of the Chow quotient). Note that we are tacitly  using that the GIT quotient is Hausdorff (more precisely that we have uniqueness of the limit).   
 
 To finish the proof it remains to prove that $$\overline{\mathcal{U}_{Smooth}^{KE}} \supseteq \mathcal{U}^{KE}.$$
 Thus let $p_0 \in \mathcal{U}^{KE}$ and we claim that we can find a sequence of $p_i \in \mathcal{U}_{Smooth}^{KE}$ converging to $p_0$. By assumption we know that $p_0=[X_0]$ where $X_0$ is a Chow polystable cubic admitting a KE metric. Applying Conjecture \ref{KEDC},  we can find a sequence of smooth KE $X_i$ GH converging to $X_0$. Now setting the $X_i$ in Tian's Gauge $T(X_i)$, eventually taking a subsequence we can assume $[T(X_{i_j})] \rightarrow [T_\infty(G)]$, where $G$ denotes the GH limit of the $X_i$ which, by the uniqueness of the GH limit, must be isomorphic to $X_0$ and, since the limit anticanonical map $T_\infty$ is a biholomorphism, to $T(G)$. Thus defining $p_j:=[T(X_{i_j})] \in  \mathcal{U}_{Smooth}^{KE}$ we have that $p_j \rightarrow [T(G)]=[T(X_0)]=[X_0]=p_0$, as we desired. 
 
 \qed
 \end{dimo}
 
 Now we are  able to prove the Conjecture \ref{C}:
 
 \begin{thm} If Conjecture \ref{KEDC}  then Conjecture \ref{C} holds.
 
 \end{thm}
 
 \begin{dimo} First of all note that, since the Hilbert scheme $\p^{19}$ is connected, ${\overline{\mathcal{M}_3}}^{Chow}$ is also connected (a GIT quotient of a connected variety is connected). Thus, since $\mathcal{U}^{KE} \neq \emptyset$, $$\mathcal{U}^{KE}={\overline{\mathcal{M}_3}}^{Chow}.$$
The fact that the natural map $D$ is continuous follows by by the deformation theory of KE orbifolds that we are assuming. Note that $D$ is also  automatically closed and proper, being a continuous map between two Hausdorff spaces. The fact that the map is generically $2:1$ follows by our Theorem \ref{S}, since cubic surfaces are not product of two lower dimensional Fano manifolds.

\qed
 \end{dimo}

\subsection*{ The geometry of the compactified moduli space}

The ring of invariants of cubic surfaces in $\p^3$ is classically very well-known (compare for example \cite{DI12} pag. 534):
$$\mathcal{R}= \C[I_8,I_{16},I_{24}, I_{32}, I_{40}, I_{100}] / <I^{2}_{100}-F(I_8, I_{16}, I_{24}, I_{32},I_{40} >,$$
where $I_j$ is a (real) homogeneus polynomial in $Sym^j(\C^{20})$. Thus
$$\overline{\mathcal{M}_3}^{ALG} = Proj\, \mathcal{R} \cong \p(1,2,3,4,5).$$
Moreover the locus of smooth cubics is given by the complement of the ample divisor in $\p(1,2,3,4,5)=Proj \, \C [x,y,z,t,u]$ (with obvious degrees) with equation $(x^2-64y)^2-2^{11}(8t+xz)=0$. Observe that the equation which defines the divisor is not canonical, i.e., depends on the choice of invariants (which we have taken to be as in \cite{DI12}).  For example, the point corresponding to the Cayley cubic is $[1:-1:-1:-1:2]$.

Then it is immediate to observe that $\overline{\mathcal{M}_3}^{ALG}$ is singular at the points $[0:0:1:0:0]$ (a quotient by $\Z_3$), $[0:0:0:0:1]$ (a quotient by $\Z_5$) and along the rational line $[0:\bigstar:0:\bigstar:0]=\p(2,4) \cong \p(1,2) \cong_{top} S^2$.

It can be shown that the set theoretically defined involution on the moduli space of cubics given by taking a cubic to its complex conjugate reduces to the obvious complex conjugation map on the moduli space (hence it is anti-holomorphic). The fixed set of the involution is $\R \p(1,2,3,4,5)$. The $\Z_2$-quotient of $\overline{\mathcal{M}_3}^{ALG}$   by this involution  (i.e., the expected GH compactification) is a smooth orbifold. Its topological Euler characteristic is equal to three.

\section{Del Pezzo of degree two and one}

It is a classical  fact (compare \cite{IP99}) that every Del Pezzo surface of degree $2$ can be realized as a degree $2$ branched cover of $\p^2$ ramified along a smooth quartic surface (the anticanonical linear system provides such  a covering).

Therefore a good candidate for the moduli (and its compactification) of degree 2 Del Pezzo  surfaces is given by the moduli of quartics in $\p^2$, i.e. by non-hyperelliptic curves of genus three:
$$\mathcal{M}_2 \subseteq \overline{\mathcal{M}}_{quar}^{GIT}:= \p(\mbox{Sym}^4 \C^3) // SL(3,\C).$$

Comparing \cite{M77} we see that:

\begin{prop} Let $Q$ be a quartic surface in $\p^2$. Then
\begin{itemize}
 \item $Q$ is stable iff $Q$ is smooth or it has at most nodes (i.e., locally $z^2=w^2$) or cusps (i.e., $z^2=w^3$) as singularities;
 \item $Q$ is strictly polystable iff $Q=(x^2+y^2+z^2)^2$ (up to $SL(3)$) or a union of two tangent conics where at least one is smooth.  
\end{itemize} 
\end{prop}

Ghigi and Kollar proved in \cite{GK07} that every Del Pezzo surface of degree $1$ with $A_1$ or $A_2$ singularities admits a KE orbifold metric. Taking double covers of $\p^2$ ramified along singular quartics, we can rephrase their result as follow:

\begin{prop}
 Let $X_2\rightarrow  \p^2$ be a degree $2$ Del Pezzo whose ramification locus is a GIT-stable quartic $Q\in \p^2$. Then $X_2$ admits a (unique) KE orbifold metric.
\end{prop}

By our result on deformations of nodal KE Del Pezzo surfaces, we know that many of these singular Del Pezzo surfaces appear as GH limits. It is natural to believe that in reality all such singular Del Pezzo appear as GH degenerations, i.e., $$\left(\p(\mbox{Sym}^4 \C^3)^{s}/SL(3)\right)/\Z_2 \subseteq \overline{\mathcal{M}}^{GH}_2,$$
where the $\Z_2$-action is induced by complex conjugation.

However, we can easily observe that (a $\Z_2$-quotient of) the compact  space $\overline{\mathcal{M}}_{quar}^{GIT}$ has to be different from the GH compactification:
\begin{prop}
 (A  $\Z_2$-quotient of) $\overline{\mathcal{M}}_{quar}^{GIT}$ is not a ``coarse'' moduli space for smooth KE degree $2$ Del Pezzo  surfaces and their degenerations.
\end{prop}

\begin{dimo} If it were the case, we would have a KE metric $\omega$ on the space $X$ given by taking the double cover of $\p^2$ branched on the not-reduced quartic $Q=(x^2+y^2+z^2)^2$. Moreover $(X,\omega)$ would be a GH limit of smooth Del Pezzo surfaces. But this is clearly impossible, since the space $X$ has two irreducible components while GH limits of Del Pezzo surfaces must be irreducible. 
 \qed
\end{dimo}
  
We can heuristically describe from the point of view of differential geometry what is happening  around the special point $[ (x^2+y^2+z^2)^2] \in \overline{\mathcal{M}}_{quar}^{GIT}$ in this way: let $([Q_i])$ and $([P_i])$ two distinct sequences of smooth quartics both converging to the special point $[ (x^2+y^2+z^2)^2]$, and let $(X_{Q_i},\omega_i)$ and $(X_{P_i}, \eta_i)$ the corresponding degree $2$ KE Del Pezzo surfaces. Then generically we  have that
$d_{GH}(X_{Q_i},X_{P_i})\geq C >0$. In particular, up to subsequences, $X_{Q_i} \rightarrow X_1$ and $X_{P_i} \rightarrow X_2$ in the GH topology with $X_1 \neq X_2$, i.e. if we approach the special points by different directions we have different metric limits.

 What the previous considerations suggest is that the GH compactification of degree two Del Pezzo surfaces should be a sort of blow up of $\overline{\mathcal{M}}_{quar}^{GIT}$ at the special point.

In the remaining part of the section, we describe what we expect is really happening around the non-reduced quartic. First of all, let us recall that smooth quartics in $\p^2$ are never hyperelliptic (recall that a curve is called hyperelliptic if it admits a branched double cover to $\p^1$). The space where hyperelliptic curves of genus three naturally embed, is the weighted projective space $\p(1,1,4)$. In this case the equations are
$$z^2=\sum_{i=0}^{8} a_i x^iy^{8-i}.$$
Now let $X$ be a Del Pezzo of degree $2$. Then $X$ defines a pair $(\p^2,Q)$, where $X \rightarrow \p^2$ is a $2:1$ map with ramification divisor a quartic $Q$. In order to find a nicer compactification of degree $2$ Del Pezzo, we may think to consider degenerations of pairs $(\p^2,Q)$. It is known \cite{HP10} that $\p^2$ admits $\Q$-Gorenstein degenerations to $\p[1,1,4]$. This suggests that we can consider degenerations of degree $2$ Del Pezzo surfaces to double cover of $\p[1,1,4]$  ramified along hyperelliptic genus $3$ curves (possibly degenerate).
Observe that these degenerate Del Pezzo surfaces have always two $\C^2 / \Z_4$ singularities (non-canonical). 

It is natural to ask if one can find KE metrics on these strictly log-terminal Del Pezzo surfaces. It turns out that we have already described an example (compare Proposition \ref{12}):  it is easy to see that the double cover $X_2$ of $\p(1,1,4)$, ramified along the degenerate hyperelliptic curve $z^2=x^4y^4$, is a \emph{toric} degree $2$ del Pezzo surface with $Sing(X_2)=2A_3+2\C^2 /\Z_4$. Thus

\begin{lem}
The $\Q$-Gorenstein smoothable KE toric degree 2 Del Pezzo $Y_2$ in Proposition \ref{12}  is isomorphic to the surface $X_2$ obtained as double cover of the weighted projective plane $\p(1,1,4)$ branched along the degenerate genus $3$ hyperelliptic curve of equation $z^2=x^4y^4$. 
\end{lem}

The degenerate hyperelliptic curve $z^2=x^{4}y^4$ is somewhat special. A hyperelliptic curve of genus $3$ defines naturally $8$ points on $\p^1$ (ramifications). Then one can study the action of $SL(2,\C)$ on the space parametrizing unordered sets of $8$ points in $\p^1$, i.e., binary octics. Then by classical GIT, we know that $8$ points are stable if no more than three points are equal and that there is only one strictly polystable orbit corresponding to the cycle given by $4$ points at 0 and $4$ points at infinity in $\p^1$. This last most degenerate case corresponds exactly to the hyperelliptic curve $z^2=x^4y^4$.

Thus it is natural to think of the  KE Del Pezzo surface $X_2=Y_2$ as the most degenerate ``stable'' degree $2$ Del Pezzo surface. Moreover by Theorem \ref{NDP},  the singularities which appears in the GH compactification must be exactly of the type $\C^2 /\Z_4$, beside the canonical ones.

With the above observations in mind, we can now formulate a more precise version of the Main Conjecture \ref{MC} in the case of Del Pezzo of degree $2$: 
\begin{conj}
There is a compact projective algebraic variety $\mathcal{X}_2= \overline{\mathcal{M}_2}^{ALG}$ together with a $2:1$ ramified continuous maps
$D:\mathcal{X}_2 \longrightarrow \overline{\mathcal{M}}^{GH}_2,$
such that
\begin{itemize}
 \item There exists a birational morphism $$\pi: \mathcal{X}_2 \longrightarrow \overline{\mathcal{M}}_{quar}^{GIT}:= \p(\mbox{Sym}^4 \C^3) // SL(3,\C),$$
which contracts a divisor $E\cong \overline{\mathcal{M}}^{GIT}_{Octics}=\p(Sym^{8}(\C^2))//SL(2,\C)$ to the special point $[(x^2+y^2+z^2)^2]$ and is an isomorphism away from $E$.
\item Each point $p\in E$ represents a singular degree $2$ Del Pezzo surface obtained as a ramified degree two covering of $\p(1,1,4)$ ramified along a possibly degenerate hyperelliptic curve $C$ of genus three given as a degree two cover of $\p^1$ branching on a polystable binary octic. In particular all such Del Pezzo surfaces contain a pair of $\C^2/\Z_4$ singularities.
\item Every Del Pezzo surface naturally associated to a point of $\mathcal{X}_2$ admits an (orbifold) KE metric. The map $D$ is then the canonical map that associates to each KE Del Pezzo surface the induced metric structure;
\end{itemize}

In particular, GH degenerations are parametrized by a reducible divisor $F=E \cup G \subseteq \mathcal{X}_2$. The singularities which can appear in the GH limits, are at most of type $A_1,A_2,A_3$ (canonical) or $\C^2 /\Z_4$ (log-terminal).
\end{conj}

Observe that $\mathcal{X}_2$ is not equal to the Deligne-Mumford compactification of the moduli of genus three curves. On the other hand it seems quite important to understand if the above $\mathcal{X}_2$ is itself given by a GIT quotient. A natural guess could be the following. Del Pezzo surfaces with canonical singularities embed (using the bianticanonical map) in $\p^6$. It is easy to check that also $Y_2=X_2$ above has the same property. Then one can consider the GIT quotient of the Chow variety parametrizing, in particular, the cycles in $\p^6$  associated to degree $2$ del Pezzo. It seems interesting to see if the above $\mathcal{X}$ can be realized as (a component of) the GIT quotient variety parameterizing Chow polystable points.  

We finish this chapter with  a small remark on degree $1$ Del Pezzo surfaces. It is known (compare \cite{IP99}) that the bianticanonical embedding realizes any degree $1$ Del Pezzos as double cover over a quadric cone in $\p^3$. The ramification locus is given by the intersection of the quadric cone with a cubic surface (i.e., a genus $4$ curve). In the degree one case we have not yet identified a possible candidate for the GH compactification. However, we should recall that in \cite{CK10} the authors proved that there are KE orbifold metrics on many singular degree one Del Pezzo surfaces with  $A_1$, $A_2$, $A_3$, $A_4$, $A_5$ and $A_6$ singularities. Moreover we should recall our example of a smoothable Del Pezzo surface of degree $1$ with genuine non-canonical singularities. It seems to be interesting to realize this example ``geometrically'' (we know only that it embeds in $\p^6$ by the three anticanonical linear system): as before, it may be used to obtain a precise picture of what the GH compactification should be.
It is also interesting to see if the moduli space of degree $1$ Del Pezzo surfaces is related to the DM compactification of genus $4$ curves (observe that, by dimension counting, the moduli space of degree $1$ Del Pezzo surfaces is one dimension smaller that the DM moduli space). It is also reasonable to believe that the GH compactification agrees with (a component of ) the Chow quotient of degree $1$ Del Pezzos in their third anticanonical embeddings (in $\p^6$).

\section{Some concluding remarks}

By the above discussion it is evident what the GH compactification picture is (or should be) in the case of Del Pezzo surfaces. Moreover   the relations between GH limits and ``stable'' Del Pezzo surfaces is particularly interesting. A way of rephrasing the results of this Chapter in the spirit of the YTD Conjecture \ref{YTD} is the following:

Let $X$ be a degree $d$ (with $d\in\{1?,2?,3(?),4\}$) $\Q$-Gorenstein smoothable Del Pezzo surface (in a family $\mathcal{H}_d$). Then

$$ X \mbox{ is } KE \,\Longleftrightarrow \, X \mbox{ is ``$GIT_d$-polystable'' },$$
Moreover, up to finite covering, $$\overline{\mathcal{M}}^{GH}_d\cong_{top} \overline{\mathcal{M}}^{ALG}_d:= \mathcal{H}_d // G$$ (at least for degree $3$ and $4$. A ``glueing'' of GIT quotients may be required in deg 2 or 1 cases).

Here $GIT_d$-stability denotes a particular GIT  notion of stability, ``ad hoc'' chosen for each $d$ (e.g., $2$-Hilbert for degree $4$, Chow for degree  $3$). In this two cases, $\mathcal{H}_4$ is the Grassmanian parameterizing intersections of two quadrics and $\mathcal{H}_3=\p^{19}$ (the Hilbert scheme parameterizing cubic surfaces).

Recall that  ``$\Longrightarrow$'' direction follows by proving that GH limits are always inside the family $\mathcal{H}_d$ and that $GIT_d$ unstable varieties can not admit a KE metric (hence they must be $GIT_d$ polystable). The direction ``$\Longleftarrow$'' is a bit indirect. Essentially it is  a consequence of a topological argument which implies that \emph{all} $GIT_d$-polystable Del Pezzo surfaces are indeed GH limits (of smooth) KE Del Pezzo surfaces. Hence by GH degeneration a $GIT_d$-polystable Del Pezzo surface must be KE too.

It would be particularly important to understand if the above situations can be unified using a similar stability notion in all cases. A possible classical candidate is Chow stability (the one used to study cubic surfaces). It may be the case that the GH compactification is related to Chow quotient for all $d$ Del Pezzo. For degree $1$ and $2$ the candidate seems to be Chow stability in $\p^6$.

On the other hand, it is interesting to see the relation with K-polystability. By the very recent paper of Berman \cite{B12} it follows that GH limits of KE Del Pezzo surfaces are indeed K-polystable. Hence the moduli spaces  $\overline{\mathcal{M}}^{ALG}_d$ are \emph{compact and projective} varieties ``parameterizing'' (a priori) a subclass of $\Q$-Gorenstein smoothable K-polystable Del Pezzo of degree $d$. Thus a natural (algebraic geometric) question is: are there any other $\Q$-Gorenstein smoothable K-polystable Del Pezzo surfaces of degree $d$ beside the ones appearing in the GH ($GIT_d$) compactification? It is seems that if any ``reasonable'' (i.e., ``Hausdorff'') algebraic (moduli) structure exists on the \emph{set} of isomorphism classes of $\Q$-Gorenstein smoothable degree $d$ Del Pezzo surfaces, then this ``space'' must coincide with the $GIT_d$ compactification considered above. In this direction, recall that the set of \emph{smooth} K-polystable Del Pezzo surfaces has a natural structure of an \emph{irreducible} quasi-projective variety (it follows by \cite{T90}, \cite{D05b} together with the identification of this set with the complement of a divisor in the above $\overline{\mathcal{M}}^{ALG}_d$).

\end{chapter}

\begin{chapter}{Degenerations $\&$ Moduli of higher dimensional Fano manifolds}
 
In this last Chapter we discuss some possible generalizations of the problems discussed in the previous part of the Thesis. In the first section we introduce the notion of cKE metric (conical KE) at rate $\mu$ on $\Q$-Fano varieties with isolated singularities. Of course, not all $\Q$-Fano varieties can admit KE metrics. In particular, we prove a  criterion for the non-existence for these cKE metrics. The criterion depends on some properties of the singularities.
Then we briefly discuss the deformation theory of cKE $\Q$-Fano varieties, focusing on the case of nodal singularities.

Finally, we show an example of a (possible) GH compactification of KE Fano manifolds in dimension three (intersections of two quadrics in $\p^5$).

\section{ KE metrics on singular Fano varieties}

As we briefly discussed in the first Chapter, the existence of KE metrics on singular $\Q$-Fano varieties is a subtle problem. In the recent paper \cite{BBEGZ11}, the authors consider ``weak'' solutions to the complex Monge-Amp\`ere equation in the sense of pluripotential theory. They provided some sufficient conditions for the existence of such  metrics. However, nothing is known about the asymptotic behavior of the metric near the singular set. This last fact is particularly unsatisfactory when one want to study the local structure of the moduli space at these singular KE Fano manifolds (compare the discussions on Del Pezzo surfaces in the previous Chapters).

Since it is expected that the KE metric  is asymptotic near the singularities to  CY cones (in the case of isolated singularities), we consider an ``ad hoc'' definition of singular KE metrics with prescribed asymptotic behavior at the singularities. However, we should remark that there are no known examples of these kinds of singular spaces.     

\begin{defi}
Let $X^n$ be a $\Q$-Fano variety whose singularity set consists of isolated points $\{p_1, \dots , p_n\}$. We say that $X^n$ admits a conical KE metric (cKE) (at rate $(\mu_i)$, with $\mu_i \in (0,2]$),  if the following properties hold:

\begin{itemize}
 \item On $X^n\setminus Sing{X^n}$ there exists a smooth KE metric $\omega$ satisfying $Ric (\omega)= \omega$;
 \item The metric completion of $X^n\setminus Sing{X^n}$ is naturally identified with $X^n$ itself;
 \item Around each singular point we can find biholomorphic maps $$\psi_i:U_i \subseteq X^n \rightarrow V_i \subseteq C(L_i),$$
where $(C(L_i)=\R^{+} \times L_i,g_i, \eta_i, J)$ is a Calabi-Yau cone, i.e., it admits a Ricci flat K\"ahler metric of the form
$$g_i=dr^2+r^2g_{L_i},$$
with $(L_i,g_{L_i})$ a smooth $(2n-1)$ real manifold (the link), necessarily (Sasaki)-Einstein.

The metric $\omega$ can be written around the singularities as
 $${\psi_i^{-1}}^{*} \omega= \eta_i+ \D h_i$$ 
for a continuous function $h_i \in C^{0}(V_i) \cap C^{\infty}(V_i\setminus \{0\})$ satisfying:
$$||\nabla_{\eta_i}^j h_i ||_{\eta_i}= \mathcal{O}(r^{2+\mu_i-j})$$
where $r$ denotes the distance from the tip w.r.t. the CY metric.

\end{itemize}

\end{defi}

Observe that since $\mu_i > 0$ the metric is genuinely asymptotic to a CY metric cone. Moreover note that when $n=2$ every KE Del Pezzo orbifold satisfies the above definition with $\mu_i= 2$ (it is sufficient to take (invariant) normal coordinates). Observe that cKE $\Q$-Fano with rate $\mu_i > 2$ cannot exist (this is the reason why in the definition we impose $\mu_i \leq 2$).

The basic example of these kind of singularities is the node $$V_0:=\{z_1^2+\dots+z_n^2=0\} \subseteq \C^n.$$ It is well-known that this singularity admits a CY conical metric. As we will show later, the CY metric has the explicit form 
$$\eta_{0,Stn}=\D(|z|^{2(1-\frac{1}{n})})_{|V_0},$$
where $|z|^2=\sum_i |z_i|^2$. Moreover, a standard computation in Sasakian geometry  shows that, setting $r^2:= \lambda |z|^{2(1-\frac{1}{n})}$  for a suitable constant $\lambda$, the metric is a cone metric, i.e.,
$$g_{0,Stn}= d r^2 + r^2 g,$$
where $g$ is an (Einstein) metric on the link $L=S^{n-1}$-bundle over $S^n$. Note the relation between the intrinsic distance function from the tip $r$ and the ``extrinsic'' distance $|z|$: $r \sim |z|^{(1-\frac{1}{n})}$.

In order to understand GH  compactifications of KE Fano manifolds in higher dimension, it is important to consider not only isolated singularities. For example in the case of threefolds one needs to consider the possibility that the singular locus consists of a complex curve (compare $\cite{DS12}$). If the singular locus is smooth then the KE metrics to consider should be of orbifold type. If the singular locus is singular, then one expects that the KE metric (of a GH limit) should be ``generically'' orbifold and, at the singular point of the curve, it should be asymptotic to CY cones with complex link equal to singular KE Del Pezzo surfaces (with genuine complex codimension $2$ singularities).  

 More generally, one should consider $n$-dimensional Fano varieties where the KE metric should be generically  asymptotic at the singularities to ``bundles'' of CY cones. Of course, very little is known on the properties of these ``stratified'' spaces.

In the next section we show an elementary obstruction to the existence of cKE metrics on $\Q$-Fano varieties.

\section{A criterion for the non-existence of cKE metrics on $\Q$-Fano varieties with isolated singularities}

In this section we extend   Proposition \ref{NDP} to higher dimensional  varieties. That is, we show that a $\Q$-Fano variety $X$ with singularities modeled on CY cones cannot admit a  cKE metric (of any rate) if the singularities are too bad compared to the degree of $X$. The main ingredients of the proof are the generalizations to the conical setting of the classical Bishop-Gromov volume comparison and Myers' Theorem. 

To start, observe the following fact:
\begin{lem}\label{GC} Let $(X^n,\omega)$ be a cKE Fano variety. Then $ X^n \setminus \mbox{Sing}(X^n) $ is geodesically convex, i.e., for all $p,q \in X^n \setminus \mbox{Sing}(X^n)$ there is a minimizing geodesic all contained in the smooth locus. 

\end{lem}
\begin{dimo}
Suppose that $X^n \setminus \mbox{Sing}(X^n)$ is not geodesically convex. Then we have a minimizing curve connecting $p$ to $q$ going through a singularity $s\in \mbox{Sing}(X^n)$. Taking the metric tangent cone at $s$ (which exists and is equal to the CY cone metric by definition), we find a minimizing curve in the CY cone going through the tip of the cone. However this is not possible, since singular CY cones minus the tip are geodesically convex.
\qed \\
\end{dimo}

The next Proposition is the generalization to the conical case of the Bishop-Gromov's volume comparison.

\begin{prop} Let $(X^n,\omega)$ be a cKE Fano variety (with normalization $\mbox{Ric}(\omega)=(2n-1)\omega$). Then for all $p \in X^n$, the function
$$r \mapsto \frac{Vol(B(p,r))}{Vol(\overline{B}(r))},$$
is non-increasing with limit as $r \rightarrow 0$ equal to $\frac{Vol(L)}{Vol(S^{2n-1}(1))}$. Here $\overline{B}(r)$ denotes the ball in the $2n$-sphere of radius $1$, and $L$ the real link in the model CY cone at $p$. 
\end{prop}

\begin{dimo}

It is known that the Bishop-Gromov volume comparison holds for geodesically convex domains (compare Remark 4.1 in \cite{CGT82}). Thus Lemma \ref{GC}, together with $$\mbox{Meas}(\mbox{Sing}(X^n))=0,$$ implies that the Proposition is true for points in the smooth locus. 

Let $p \in \mbox{Sing}(X^n)$ and take $(p_i) \in X^n \setminus \mbox{Sing}(X^n)$ converging to $p$. Since $$B(p_i,r) \rightarrow B(p,r)$$ in the Gromov-Hausdorff topology, the volume comparison estimate holds also for $p$ in the singular locus.

It remains to show that  the limit as $r \rightarrow 0$ is equal to $\frac{Vol(L)}{Vol(S^{2n-1}(1))}$. Since the metric is asymptotic to the CY metric on the cone, we find that
$$
Vol(B(p,r))=\frac{r^{2n}}{2n} Vol(L)+\mathcal{O}(r^{2n+1})
$$
Similarly
$$Vol(\overline{B}(r))=\frac{r^{2n}}{2n} Vol(S^{2n-1}(1))+\mathcal{O}(r^{2n+1}).$$
 The Proposition follows immediately.

\qed \\
\end{dimo}

The following Proposition generalizes  Myers' Theorem:
\begin{prop} Let $(X^n,\omega)$ be a cKE Fano variety (with normalization $\mbox{Ric}(\omega)=(2n-1)\omega$). Then
$$\mbox{diam}(X^n) \leq \pi.$$
\end{prop}

\begin{dimo}
Since $\mbox{diam}(X^n)=\mbox{diam}(X^n \setminus \mbox{Sing}(X^n))$, it is sufficient to estimate the distance between two points in the smooth locus.

 Let $p,q\in X^n \setminus \mbox{Sing}(X^n)$. By Lemma \ref{GC}, we have a smooth minimal geodesic $\gamma$ all contained in the smooth locus connecting $p$ to $q$. Then the usual argument in the proof of Myers' Theorem (compare for example \cite{GL87}) gives 
$$d(p,q)= \mbox{length}(\gamma) \leq \pi.$$
The Proposition follows by taking the sup of the distance of pairs $(p,q)$. 
\qed \\
\end{dimo}

\begin{rmk}
 The above ``Bishop-Gromov and Myers Theorems'' hold more generally under the weaker assumption $\mbox{Ric}(\omega)\geq C > 0$. 
\end{rmk}

Now we are ready to state and prove the main result
\begin{thm} \label{NF}Let $X$ be a $\Q$-Fano variety  with an isolated singularity modeled on CY cones with real link $L_p$.  If $X$ admits a cKE metric, then

\begin{equation}\label{OBS}
 \frac{\mbox{Vol}(S^{2n-1}(1))}{\mbox{Vol}(L_p)} \mbox{deg}(X) \leq n! \frac{(2n-1)^n}{(2\pi)^n}\mbox{Vol}(S^{2n}(1))
\end{equation}

\end{thm}

\begin{dimo}
 The proof is identical to the proof of Proposition \ref{NDP}, keeping track of Ricci tensor normalizations.
\qed \\ 
\end{dimo}

Let us note that $\tfrac{\mbox{Vol} (S^{2n-1}(1))}{\mbox{Vol}(L_p)}\geq 1$ (by Myers' Theorem), and it is equal to one only if $p$ is a smooth point. 

Before showing that the inequality \ref{OBS} can be actually used to prove that some singular Fano varieties cannot admit cKE metrics asymptotic near the singularity to a specific CY cone metric, we make some useful remarks.

\begin{itemize}
\item A computation (see next Corollary) of the RHS of \ref{OBS} gives $\mbox{RHS}=2,12,100$ and $\frac{5488}{5}$ for $n=1,2,3$ and $4$ respectively. Thus the RHS manifests roughly an ``exponential behavior''.

\item  For an orbifold singularity $\C^n / \Gamma$, with $\Gamma \subseteq U(n)$ finite, the quantity $\frac{\mbox{Vol}(S^{2n-1}(1))}{\mbox{Vol}(L_p)}$ is simply given by $|\Gamma|$. 

More generally (compare \cite{GMSY07} and \cite{S11}), if the CY cone is modeled on a (quasi)-regular Sasakian cone with $Z_p$ as quotient of the real link $L_p$  by the Reeb field, then
\begin{equation}\label{Bis}\frac{\mbox{Vol}(L_p)}{\mbox{Vol}(S^{2n-1}(1))}=\frac{\mbox{ind}(Z_p)}{n^n}\,\int_{Z_p}c_1^{orb}(Z_p)^{n-1},
\end{equation}
where $\mbox{ind} (Z_p)$ denotes the index, i.e., the maximal integer divisor of the anticanonical class in the Picard group. 
For example the standard CY metric on the $n$-dimensional double point, i.e., $\{z_1^2+\dots+z^{n+1}=0\}$, satisfies
$$ \frac{\mbox{Vol}(L_p)}{\mbox{Vol}(S^{2n-1}(1))} = 2\left(1-\frac{1}{n} \right)^n \rightarrow \frac{2}{e}.$$

\item Observe that combining the formula \ref{Bis} with  $\tfrac{\mbox{Vol} (S^{2n-1}(1))}{\mbox{Vol}(L_p)}\geq 1$,  the authors of \cite{GMSY07} found what they called ``Bishop obstruction'': a necessary condition for a Fano orbifold $X^n$ (suppose, for simplicity, with codimension two quotient singularities) to be the quotient by the Reeb vector field of a link in a Sasaki CY cone -thus $X^n$ will be in particular KE-  is that:
$$ deg(X^n) \leq \frac{(n+1)^{n+1}}{\mbox{ind}(X^n)}.$$
This obstruction is weaker than the restriction on the degree of our Theorem \ref{NF}. In \cite{GMSY07} the authors are mainly interested in CY cones with the tip as the only singular point. For this reason, they apply the Bishop-Gromov inequality to the \emph{smooth} (Einstein) link of a CY cone.  On the other hand, our obstruction  has been found applying the Bishop-Gromov volume comparison \emph{directly} on the KE Fano orbifold, which a-priori should be related  to CY cones singular along complex codimension two subsets.

Moreover, we should recall that in the case of KE Fano manifolds with a $\C^*$-action with a finite number of fixed points, the ``Bishop obstruction'' has been recently improved by R. Berman and B. Berndtsson to $ deg(X^n) \leq (n+1)^{n}$ \cite{BB12}.

\item Finally recall that the degree $deg(X)$ of a $\Q$-Fano variety is unbounded in the class of log-terminal Fano varieties. However it is conjectured to be bounded as long as the discrepancies are bigger than $-1+\epsilon$ (Borisov-Alexeev-Borisov Conjecture, proven in dimension $2$ by V. Alexeev in \cite{A94}). 
\end{itemize}

By combining the above observations with the computation of the volume of the even dimensional spheres, we find the following rephrasing of Theorem \ref{NF}:

\begin{cor} Let $X^n$ be a $\Q$-Fano variety  with an isolated singularity at $p$  modeled on a (quasi)-regular Sasakian CY cone and let $Z_p$ be the (orbifold) Fano variety given by the quotient 
 of the real link of the cone by the Reeb  foliation, then
$$ \frac{n^n}{\mbox{ind}(Z_p) \,c_1^{orb}(Z_p)^{n-1}}\,\mbox{deg}(X) \leq 2\frac{n!(2n-1)^n}{(2n-1)(2n-3)(2n-5)...1}.$$
\end{cor}

Now we give an example of a family of Fano 3-folds violating the constraint \ref{OBS}.

\begin{prop} Let $X_d$ be the toric $\Q$-Fano 3-fold whose Fano polytope $P_d$ is defined as 
$$P_d=\mbox{Convex} \, (\{e_1,e_2,e_3,-e_1-e_2-de_3\}).$$
Then $X_d$ does not admit a cKE (actually orbifold) metric for $d>1$. 
\end{prop}

\begin{dimo}
It is easy to see that $X_d$ is a $\Q$-Fano with only one singular point corresponding to 
$$\mbox{Cone}(e_1,e_2,-e_1-e_2-de_3\}).$$
Computing the algebra of the dual cone, we find that the singularity is given by
 $$\mbox{Spec}\left(\C[z,zy,\dots,zy^d,zx,zxy,\dots,zxy^{d-1},zx^2,\dots,\dots,zx^d]\right).$$
Changing the variable to $z=a^d,y=\frac{b}{a},x=\frac{c}{a}$, we immediately observe that the above singularity is simply the cone over the degree $d$ Veronese embedding of $\p^2$, i.e., $$\mbox{Spec} \left(\C[a,b,c]^{{\Z_d}}\right)$$

The degree of a toric Fano 3-fold is well-known to be equal to $\mbox{deg(X)}=3! \mbox{Vol}(P^{\checkmark})$, where $P^{\checkmark}$ denotes the dual polytope of P.  In our case $P_d^{\checkmark}$ is given by the convex hull of the following vectors
$$v_1=(-1,-1,-1);v_2=(-1,-1,\frac{d}{3});v_3=(-1,d+2,-1); v_4=(d+2,-1,-1).$$
It follows that $\mbox{deg}(X_d)=\frac{(d+3)^{3}}{d}$.

Hence the LHS of \ref{OBS} is equal to $(d+3)^{3}\geq 125 $, for $d>1$. Since the RHS of \ref{OBS} is equal to $100$, we find that $X_d$ cannot admit a cKE (orbifold) metric.

\qed \\
\end{dimo}

The above example can be generalized in all dimensions simply considering the Fano polytope $P_d^n=\mbox{Convex}(\{e_1,\dots,e_n,-e_1-\dots-e_{n-1}-de_n\})$, giving log-terminal Fano varieties with no cKE metric (note that the LHS of \ref{OBS} is equal to $(d+n)^{n}$).
Moreover it is important to note that these singular Fano varieties are not smoothable, the cones over the Veronese embeddings being rigid singularities
 in dimension bigger than two. 

Finally we remark that the non-existence of  KE orbifold metrics on the above Fano varieties follows also from toric geometry, since the baricenter of the dual polytope is not the origin \cite{SZ11}. From a purely combinatorial point of view, it would be interesting to see if the above obstruction to the existence of a KE metric gives  properties that ``balanced'' polytopes must satisfy which are unknown (compare the discussion on the Ehrhart Conjecture in \cite{BB12}).

Moreover we remark that it seems that the above obstruction holds, in particular, for GH limits of smooth KE Fano manifolds (\cite{DS12}). It would be nice to see if this obstruction can be used to study GH compactifications as we did in the case of Del Pezzo surfaces. However ``experimental computations'',  that we performed on some classes of Fano $3$-folds, suggest that this restriction is too weak to be particularly useful.

\section{Some observations on deformations of cKE $\Q$-Fano varieties}

Let $(X_0, \omega_0)$ be a cKE metric on a $\Q$-Fano variety as defined in the first section of this Chapter and  let $\mathcal{X} \rightarrow \Delta_t \subseteq \C$ be a (partial)-smoothing. The fundamental question is the following: when does $X_t$ admit a (cKE) metric? As we discussed in the previous Chapters, the automorphism group plays a central role in understanding the deformation theory of (singular)  KE Fano varieties.

The next Proposition is a result relating the automorphism group of a cKE variety $(X,\omega)$ to eigenfunctions of the Laplacian. We say that $\phi \in C^k_{\beta, \omega}(X)$, if $\phi \in C^k_{\beta, \omega}(X \setminus Sing(X),\R)$ and near each singularity $p_i\in Sing(X)$ we have $|\nabla_{\omega}^j \phi|_{\omega} = \mathcal{O}(r^{\beta-j})$ for $r<<1$ and $j\leq k$, where $r$ is the distance from the tip of the cone w.r.t. the CY  metric $\eta_{i}$.

Assume that $X$ has only discrete automorphism group. This implies that $H^0\left(X, \mathcal{H}om(\Omega_X^1, \mathcal{O}_X)\right)=0$ (here $\Omega_X^1$ denotes the sheaf of K\"ahler differentials). If the local cohomology group $H^1_{Sing(X)}\left(X, \mathcal{H}om(\Omega^1_X, \mathcal{O}_X)\right)=0$, then also $H^0\left(X\setminus Sing(X), TX^{1,0}\right)=0$, i.e. there are no non-trivial holomorphic vector fields on the smooth part.  
Now we are ready to state the Proposition.

\begin{prop} \label{ker} Let $(X^n,\omega)$ be a cKE $\Q$-Fano variety with discrete automorphism group and such that $$H^1_{Sing(X)}\left(X, \mathcal{H}om(\Omega_X^1, \mathcal{O}_X)\right)=0.$$  Suppose that there is a  function $\phi \in C^2_{\beta, \omega}(X,\R)$ for $\beta>2-n$ satisfying
$$\Delta_\omega \phi = \phi \; \mbox{ on $X \setminus Sing(X) $}.$$
Then $\phi\equiv0$. ($\Delta_\omega$ denotes the $\bar{\partial}$-Laplacian)
\end{prop}

\begin{dimo}

Consider the positive quantity
$$g^{i\bar{j}}g^{k\bar{l}} \phi_{;ik}\phi_{;\bar{j}\bar{l}}= | (\nabla^2\phi)^{(2,0)}|^2 \geq 0,$$
where $g_{i\bar{j}}$ is the expression of the K\"ahler metric in a coordinate system.
Working in normal coordinates, one can check that  the following Bochner-type formula is true \cite{T00}:

$$| (\nabla^2\phi)^{(2,0)}|^2 = div (\nabla^2 \phi\ast \nabla \phi)+ (\Delta \varphi)^2-Ric(\partial \phi, \bar{\partial} \phi),$$
where $\ast$ denotes a ``product like expression'' (between first and second derivatives). 

We claim that if $\phi$ is a eigenfunction of eigenvalue $1$ on a cKE with discrete automorphism, then $| (\nabla^2\phi)^{(2,0)}|^2\equiv 0$ on $X \setminus Sing(X)$. For simplicity we assume that $Sing(X)$ consists of only one node, say $p$.

Integrating the Bochner-identity on $X \setminus Sing(X)$:

$$\int_{X \setminus Sing(X)}| (\nabla^2\phi)^{(2,0)}|^2 \omega^n = $$
 $$\lim_{r \rightarrow 0} \left(\int_{X \setminus B_{\eta}(p,r)}
 div (\nabla^2 \phi\ast \nabla \phi)+ (\Delta \varphi)^2-Ric(\partial \phi, \bar{\partial} \phi))\right)=$$
 $$\underbrace{\lim_{r \rightarrow 0} \left(\int_{\partial B_{\eta}(p,r)}<\nabla^2 \phi \ast  \nabla \phi, \eta >\right)}_{A}+\underbrace{\int_{X \setminus Sing(X)} (\Delta \varphi)^2-Ric(\partial \phi, \bar{\partial} \phi))}_{B}.
$$

Here we are using the fact that, near the singular point $p$, $r \approx r_{\omega}$ (i.e. the ``real'' distance is equivalent to the CY  distance). 

\begin{itemize}

\item Estimate of $A$.

Since $\omega^n=(1+\mathcal{O}(r^{\mu})) \eta^n$ and the volume form of the metric $\eta$ restricted to the real link at distance $r$ is $r^{2n-1}dV_1$,  we find
$$ \rvert \int_{\partial B_{\eta}(p,r)}<\nabla^2 \phi \ast  \nabla \phi, \eta >\rvert \, \leq C
r^{\beta-2}r^{\beta-1}r^{2n-1} \leq C r^{2(\beta+n)-4}.$$
Thus $\beta>2-n$ implies $A=0$.

\item Estimate of $B$.

Since $\Delta\phi=\phi$ and $Ric(\omega)=\omega$, we have
$$\int_{X \setminus Sing(X)} (\Delta \varphi)^2-Ric(\partial \phi, \bar{\partial} \phi)= \int_{X \setminus Sing(X)}  \phi\Delta \phi  - \int_{X \setminus Sing(X)}|\bar{\partial}\phi|^2.$$
Moreover
$$\int_{X \setminus Sing(X)} \phi \Delta  \phi = \lim_{r\rightarrow 0} \left(\int_{\partial B_{\eta}(p,r)} <\bar{\partial} \phi \ast \phi, \eta>\right) + \int_{X \setminus Sing(X)}|\bar{\partial}\phi|^2,$$

and, as before,

$$\rvert\int_{\partial B_{\eta},(p,r)} <\bar{\partial} \phi \ast \phi, \eta> \rvert \leq C r^{\beta-1} r^{\beta} r^{2n-1} \leq C r^{2(\beta+n)-2},$$
implies that, if $\beta>1-n$ then $B=0$.

\end{itemize}

Combining the two previous estimates, we find for $\beta>2-n$ that
$$\int_{X \setminus Sing(X)}| (\nabla^2\phi)^{(2,0)}|^2 \omega^n =0.$$

Thus the tensor $(\nabla^2\phi)^{(2,0)}\equiv 0$ on $X \setminus Sing(X)$. This implies that, since the metric is parallel,
$$\nabla_{\bar{k}}(g^{i\bar{j}} \phi_{\bar{,j}})=0.$$
In particular we find that the vector field $\left(\bar{\partial} \varphi\right)^\sharp:= g^{i\bar{j}}\phi_{\bar{,j}}$ is holomorphic, and hence it must be the zero field by our assumption on the automorphisms. This implies that $\bar{\partial}\phi=0$, and, being $\phi$ real-valued, also $\partial \phi=0$. Hence $\phi$ must be constant. Finally, since $\Delta \phi= \phi$, this constant must be zero.

\qed

\end{dimo}

The above Proposition applies in the case of nodal Fano varieties (or more generally for any variety with isolated complete intersection singularities) provided $n\geq3$ (compare for example \cite{F91}).

Another  crucial ingredient for the study of deformations of (singular) KE metrics is given by the existence of good metric models (Asymptotically Conical CY) on smoothings of the singularities. The simplest possible singularity is given by the node. In this case, M. Stenzel constructed an AC CY metric on the (unique) smoothing \cite{S93}. Due to the symmetries of the node, the construction of a CY metric is quite easy (and essentially explicit).
We recall his argument. 
Let $Q_0^n$ be the $n$-dimensional complex cone
$$Q_0^n:= \{ \underline{z} \in \C^{n+1} \, | \, z_1^2+ \dots + z_{n+1}^2=0 \},$$
and let $Q^n$ be its smoothing, i.e.
$$Q^n:= \{ \underline{z} \in \C^{n+1} \, | \, z_1^2+ \dots + z_{n+1}^2=1 \}.$$
 We look for metrics of the form
$$\D f(s),$$ where $s=|z_1|^2+\dots+ |z_{n+1}|^2$. Consider the  holomorphic $(n,0)$-form
$$\Omega:=\frac{(-1)^{n}}{2z_{n+1}}dz_1 \wedge \dots \wedge dz_n.$$
Thus we seek solutions of the equation
$$(\D f(s))^n=C \Omega \wedge \overline{\Omega},$$
on $Q_0^n$ and $Q^n$, where $C \in \R$.
The equation can be rewritten as
$$i^n((f^{'}(s))^n (\partial \bar{\partial} s)^n +n {f^{''}(s)}({f^{'}(s)})^{n-1} \partial s \wedge \bar{\partial} s \wedge (\partial \bar{\partial} s)^{n-1})=C \Omega \wedge \overline{\Omega},$$
An elementary computation shows that the equation reduces to the following ODEs:
On $Q_0^n$,
\begin{equation}\label{ODEeq1}
s(f^{'}(s))^n+ s^2{f^{''}(s)}({f^{'}(s)})^{n-1}=c \geq 0,
\end{equation}
for $s>0$. On $Q^n$
\begin{equation}\label{ODEeq2}
s(f^{'}(s))^n+ (s^2-1){f^{''}(s)}({f^{'}(s)})^{n-1}=c > 0,
\end{equation}
for $ s\geq 1$.
On $Q_0$ the equation can be completely solved giving as solution the following  complete  at infinity Ricci-flat  metric
$$\eta_{0,Stn}=\D(r^{2(1-\frac{1}{n})}),$$
where $r^2=s$, i.e. $r$ is the distance function with respect to the std metric in $\C^{n+1}$. 
On $Q^n (= T^*S^n)$, it is in general harder to find an exact solution of $\eqref{ODEeq2}$. However, for $n=2$ one can check that
$$\eta_{Stn}:= \D \sqrt{r^2+1},$$
where $r^2=s$ is a complete CY metric (the Eguchi-Hanson metric, with K\"ahler form given by the complex structure of the smoothing).

In general, defining $g(s):= (f^{'}(s))^n$, the equation $\eqref{ODEeq2}$ becomes
$$s g(s) + \frac{1}{n} (s^2-1) g^{'}(s)=c.$$
Thus  $g(s)=\mathcal{O}(\frac{1}{s})$, which implies that the leading term at infinity of the Stenzel potential is $\sim r^{2(1-\frac{1}{n})}$. To perform a gluing construction, it is essential to understand the asymptotics of the CY conical metric. In this case we want an $x>0$ such that
$$\eta_{Stn}= \D(r^{2(1-\frac{1}{n})}+ \mathcal O(r^{-x})).$$

Recent progress in the asymptotics of AC CY metrics on smoothings of CY cones are contained in the paper \cite{CH12}.

Combining the precise value of the Stenzel asymptotics with the decay rate of a  cKE Fano, it is not unrealistic to expect that our Theorem \ref{MT} generalizes to higher dimensional setting, but we have not carried out all the details yet. 
Thus we expect the following generalization of our Theorem \ref{MT}:

\begin{conj}
Let $\mathcal{X} \rightarrow \Delta_t \subseteq \C$ be a (partial)-smoothing of a nodal cKE Fano variety $(X_0,\omega_0)$, for some specific rates of convergence $\mu_i$. Assume that the automorphism group is discrete. Then $X_t$ admits a (cKE) metric $\omega_t$ for $t$ sufficiently small. These metrics $\omega_t$ are GH close to $\omega_0$.
\end{conj}

\begin{rmk}
 The deformation theory of $\Q$-Fano varieties is in general more complicated than deformation theory of log-terminal Del Pezzo surfaces described in Chapter two. In particular, it is no longer true  that we have always no local-to-global obstructions. Some results on the deformation theory of Fano $3$-folds are contained in \cite{F86} and \cite{N97}.

\end{rmk}

\section{Moduli of intersections of two quadrics}

In this final Chapter we discuss what the GH compactification of the space of intersections of two quadrics in $\p^5$ should be. Similar arguments  also hold for intersections of two quadrics in $\p^n$.

First of all, we recall that all smooth intersections of two quadrics in $\p^5$ ($\p^n$) admit a KE metric (\cite{N90} and \cite{AGP06}). Note that the moduli of smooth intersections of two quadrics is ``complete'' (i.e., all smooth Fano manifolds in the same family are indeed given by intersections of two quadrics). With this in mind, it becomes particularly interesting to study the GIT stability of such intersections.

The (classical) GIT picture is identical to the one discussed in the case of Del Pezzo quartics (and the proof generalizes verbatim, compare for example \cite{AL00}).  Given two quadrics $Q_1$ and $Q_2$, one  considers the pencil $\lambda Q_1+ \mu Q_2$. Thus we identify the space of intersections of two quadrics with $Gr(2, H^0(\p^5, \mathcal{O}(2)))$. Taking the Pl\"ucker embedding, we consider the following GIT picture:
$$SL(6,\C) \curvearrowright Gr(2, H^0(\p^5, \mathcal{O}(2))) \hookrightarrow \p(\Lambda^2 H^0(\p^5, \mathcal{O}(2))).$$
Define
$$\overline{\mathcal{M}_{IntQuar}}^{ALG}:= Gr(2, H^0(\p^5, \mathcal{O}(2))) // SL(6,\C),$$
as a possible candidate to be a good algebraic compactification of intersections of two quartics. Consider the natural map of sets 
$$\begin{array}{cccc} disc:& Gr(2, H^0(\p^5, \mathcal{O}(2))) & \longrightarrow & \{ \mbox{Binary sextics} \}/SL(2,\C) \\
    & [\lambda Q_1+ \mu Q_2] & \longrightarrow & det(\lambda Q_1+ \mu Q_2).
  \end{array}
$$
Then
\begin{thm} The map $\mbox{disc}$ descends to the GIT quotient and it induces an isomorphism between $\overline{\mathcal{M}_{IntQuart}}^{ALG}$ and the moduli space of binary sextics $\overline{\mathcal{M}_{bin, 6}}^{GIT}:= \p H^0(\p^1, \mathcal{O}(6)) // SL(2,\C).$

 Every polystable intersection of two quadrics can be ``diagonalized'', i.e., it can be expressed as
$$ X: \,\begin{cases} x_0^2+x_1^2+x_2^2+x_3^2+x_4^2+x_5^2=0\\
 \lambda_0 x_0^2+ \lambda_1x_1^2+\lambda_2x_2^2+\lambda_3x_3^2+\lambda_4x_4^2+\lambda_5 x_5^2=0 \end{cases}, $$
with \emph{at most} pairs of equal $\lambda_i$ or with $\lambda_0=\lambda_1=\lambda_2$ and  $\lambda_3=\lambda_4=\lambda_5$ (this (unique) variety corresponds to the (unique) polystable but not stable binary sestic $x^3y^3=0$ via the map $disc$ ).  Conversely, every intersection of quadrics that admits the above representation is polystable.

Moreover $X=Q_1 \cap Q_2$ is stable iff all $\lambda_i$ in the above representation are distinct iff $X$ is smooth. By an obvious change of base, the polystable intersection of quadrics with $\lambda_0=\lambda_1=\lambda_2$ and  $\lambda_3=\lambda_4=\lambda_5$ is given by the equations: 
$$X_0 := \begin{cases} x_0x_1=x_2^2\\
x_3x_4=x_5^2 \end{cases}.$$
Note that, beside $X_0$, all polystable quadric admit at most  nodal singularities. $X_0$ is singular along two disjoint rational lines (see Proposition \ref{pKE}). 
\end{thm}
For the sake of completeness, note that the principal invariants of $X:=Q_1 \cap Q_2\subseteq \p^5$ are easily computed to be: $deg(X)= 32$, $\rho(X)=1$ (Picard rank), $h^{1,2}(X)=2$.
The GIT-invariants of binary sextics are classically well-known. It follows by computing the algebra of invariants (compare for example \cite{DI03} pag. 153) that 
$$\overline{\mathcal{M}_{IntQuar}}^{ALG} \cong \p(1,2,3,5).$$

In particular it has only $3$ isolated quotient singularities. The smooth point $[1:0:0:0]$ corresponds to the polystable intersection with non isolated singularities $X_0$. With the choice of invariants as in \cite{DI03}, the point $[0:0:0:1]$ corresponds to a smooth intersection of two quadrics. The other two singular points in the moduli space correspond to singular intersections. The set of smooth intersections of two quadrics is given by the complement of the divisor $t=0$ (where $t$ denotes the coordinate of weight $5$ in $\p(1,2,3,5)$). 

Another interesting point in $\overline{\mathcal{M}_{IntQuar}}^{ALG}$ is $[-16:4:8:0]$. This point corresponds to the toric polystable intersection of quadrics of equation $$x_0x_1=x_2x_3=x_4x_5.$$
It is the unique polystable intersection of two quadrics with $6$-nodal singularities. Its moment polytope is given by the regular octahedron in $\R^3$.  

Finally we briefly discuss the problem of existence of singular metrics on polystable intersections of two quadrics. First of all note that the strictly polystable intersection with non isolated singularities admits an orbifold KE metric (this fits nicely within the general degeneration theory of KE Fano manifolds \cite{DS12}).

\begin{prop}\label{pKE}
The polystable intersection of quadrics
$$Q_1 \cap Q_2 := \begin{cases} x_0x_1=x_2^2\\
x_3x_4=x_5^2 \end{cases}$$
admits an orbifold KE metric with singularities locally modeled on $\C^2/{\Z_2}\times \C$. 

Moreover $ Sing(  Q_1 \cap Q_2 )\cong \p^1 \sqcup \p^1$ is a totally geodesic submanifold with induced metric given by the Fubini-Study metric on $\p^1$. 

\end{prop}

\begin{dimo}
By considering the $2:1$ covering map
$$\begin{array}{cccc}
   p:& \p^3 &\longrightarrow & Q_1 \cap Q_2 \subseteq \p^5 \\
  & (a:b:c:d) & \longmapsto & (a^2:b^2:ab:c^2:d^2:cd), \end{array}
$$
we  see that $Q_1 \cap Q_2 \cong \p^3/\Z_2$ where the $\Z_2$-action is generated by $$[a:b:c:d] \longrightarrow [a:b:-c:-d].$$
This $\Z_2$-action is clearly an isometry w.r.t. the Fubini-Study metric on $\p^3$. Therefore the FS metric descends to an orbifold metric on the quotient, i.e., on the given intersection of quadrics.  

Finally, it is immediate to see that the singular locus of $Q_1 \cap Q_2$, given by $x_0=x_1=x_2=0$ and $x_3=x_4=x_5=0$, has as preimage under $p$ the two disjoint lines $L_1:=[a:b:0:0]$ and $L_2:=[0:0:c:d]$ inside $\p^3$.
Then the claim follows by recalling that lines in $\p^3$ are totally geodesic and that the restriction of the ambient FS metric coincides with the round metric on $\p^1$.

\qed
\end{dimo}

For intersections of quadrics with nodal singularities the situation is more subtle. A particularly nice point in the algebraic compactification is given by the toric intersection $x_0x_1=x_2x_3=x_4x_5$. Its moment polytope has ``barycenter zero'' (i.e., the Futaki invariant vanishes). Thus it is expected that such a variety admits a KE metric (actually cKE). In \cite{BB11} the authors claim (but the proof is not appeared yet) that the vanishing of the Futaki invariant in the $\Q$-Fano toric case is sufficient to have a ``weak solution'' of the KE equation. It is interesting to see if it is possible to compute the asymptotic behavior of the KE metric. In particular, it would be important to understand if the KE metric is asymptotic to the standard toric  Guillemin metric.

For the other polystable intersections the situation is more complicated. However if \cite{AGP06} generalizes to $\Q$-Fano varieties then one can try to construct singular KE metrics by considering, for example, double covers from the (singular) polystable intersections to smooth (KE) quadric threefolds in $\p^4$. More precisely, assume that $Q_1=1d$ and $Q_2=\underline{\lambda}$ with $\lambda_0=\lambda_1=0$ and distinct other $\lambda$s. Then define the (KE) quadrics $C_i:=\lambda_{i}Q_0- Q_1\subseteq \p^4$ for $i=2,3,4,5$. There are four $2:1$ maps $\pi_i: Q_0 \cap Q_1 \rightarrow C_i$, given by forgetting the $x_i$. Then the argument  in \cite{GK07} would provide the existence of a ``singular'' (weak solution in the sense of pluripotential theory) KE metric on $Q_0\cap Q_1$. A statement which seems related to the above picture can be found in $\cite{BBEGZ11}$ (compare Theorem $8.3$).

The above observations suggest that $\overline{\mathcal{M}_{IntQuar}}^{ALG}$ is indeed equal to the GH compactification $\overline{\mathcal{M}_{IntQuar}}^{GH}$. However we should remark that it does not seem immediate to generalize the Mabuchi-Mukai argument (or an argument close to the one given for the case of cubic surfaces) to this higher dimensional situation.

\end{chapter}

%Appendix
\appendix

\begin{chapter}{Fano Toric Varieties}
 
In this appendix we collect a few well-known facts about Fano Toric Varieties which we use in the Thesis. For a more detailed account see \cite{F93} or \cite{D03}.

Let $N \cong \Z^n$ and $M=Hom_{\Z}(N,\Z)$. A \emph{Fan} $\Delta$ is a finite set of cones $\sigma \subseteq N^{\R}:=N\otimes_{\Z}\R$ such that:
\begin{itemize}
 \item Every face of a cone in $\Delta$ is itself an element of $\Delta$.
 \item The intersection of two cones is a face.
\end{itemize}

For each $\sigma \in \Delta$, one can consider the affine toric variety
$$X_{\sigma}:= Spec \left( \C[\sigma^\checkmark \cap M] \right),$$
where $\sigma^\checkmark:=\{ y \in M^{\R}\; | \; <y,x> \geq 0 \; \forall x \in \sigma \}.$

Given $\Delta$, we define an $n$-dimensional toric variety $X_\Delta$ to be the normal, irreducible variety obtained by glueing $X_{\sigma}$ to $X_{\tau}$ along $X_{\sigma\cap\tau}$ for $\sigma,\tau \in \Delta$.
Note that $ ({\C^{*}})^n=X_{\{0\}} \subseteq X_\Delta$.

Let $P$ the Polytope defined as the convex hull of the set of the minimal generators in $N$ of the rays of the fan $\Delta$.
\begin{prop}
Let $X_\Delta$ be a toric variety defined as above, then 
\begin{itemize}
 \item $X_\Delta$ is smooth iff every cone in $\Delta$ is generated by basis elements of $N$.
 \item $X_\Delta$ is projective iff $\Delta$ is generated by the faces of the polytope $P$ containing the origin in its interior.
\item  $X_\Delta$ is $\Q$-factorial iff $P$ is simplicial.
\end{itemize}
In this last case the Picard rank is equal to $\rho(X_{\Delta})= \sharp \{\mbox{vertices of } P \} -n$.

\end{prop}

For projective toric varieties we use the notation $X_P$ instead of $X_\Delta$.

 A characterization of toric Fano varieties is the following:

\begin{prop}
 A projective toric variety is Fano iff $P$ has vertices which are primitive vectors in $N$ ($P$ is called  the Fano polytope).
\end{prop}

Let $P^\checkmark \in M^{\R}$ be the dual polytope (in general non integral) to $P$  defined by
$$P^\checkmark:=\{y \in M^{\R}\; | \; <y,x> \geq -1 \; \forall x \in \sigma \}.$$

Some properties of anticanonical divisor of the toric Fano variety $X_P$ are encoded in the dual polytope $P^\checkmark$:

\begin{prop}
 Let $P$ be a Fano polytope of dimension $n$ with $-K_X$ $\Q$-Cartier. Then
\begin{itemize}
 \item $deg(X_P):= (-K_{X_P})^n= n! Vol(P^\checkmark)$;
  \item $-kK_{X_P}$ is Cartier iff $kP$ is integral;
 \item $dim\, H^{0}(X_P,-kK_{X_P})= \sharp \{ \mbox{lattice points in}\, kP^\checkmark\}$.
\end{itemize}

\end{prop}

If $P$ is a simplicial Fano polytope then $X_P$ has at most orbifold (quotient) singularities. For example this is always true for Del Pezzo surfaces. More precisely every singularity of a toric Del Pezzo surface is, up to unimodular transformation, given by the cone in $N^{\R}$ with rays generated by $v_0=(0,1)$ and $v_1=(m,-k)$ with $(m,k)=1$. It is easy to check that the singularity is of type $\C^{2}/ \Z_m$ where the action is generated by $(\mu,\mu^k)$ with $\mu$ being primitive root of unity. 

For toric Fano orbifolds it is possible to detect the existence of an (orbifold) KE metric simply by looking at the dual polytope $P^\checkmark$ \cite{SZ11}: 

\begin{thm}[B. Wang-X. Zhu-Y. Shi] Let $X_P$ be a toric Fano orbifold. Then the following are equivalent:
\begin{itemize}
 \item $X_P$ admits a KE orbifold metric;
 \item $Fut(X_P)=0$;
 \item $Bar(P^\checkmark)=0$, i.e., $\int_{P^\checkmark}y_i \,dy=0$ for all $i$.
\end{itemize}
Moreover if $Bar(P^\checkmark)\neq0$, $X_P$ admits a  K\"ahler-Ricci soliton metric.
\end{thm}
Note that the proof for orbifolds is a verbatim generalization of the proof in the smooth case.

\end{chapter}

\begin{chapter}{Singularities dictionary}

In this Appendix we recall the basic definitions of some important classes of singularities which we considered in the Thesis. For a more systematic description compare \cite{KM98}, \cite{IP99} or \cite{HK10}.

Let $X$ be a $\Q$-Gorenstein normal variety (i.e., $K_X$ is a $\Q$-Cartier divisor). Let  $f:Y\rightarrow X$ be a resolution of singularities with exceptional divisor  $E=\bigcup E_i$ (where $E_i$ are the exceptional irreducible prime divisors ). Write
$$ K_Y=  f^{\ast}K_X+\sum a_i E_i, $$
where $a_i \in\Q$ (the discrepancies at $E_i$) , $f_\ast K_Y= K_X$ and equality is understood in the sense of $\Q$-divisors. Then we say that

\begin{itemize}
 \item $X$ has terminal singularities if $a_i > 0$;
 \item  $X$ has canonical singularities if $a_i \geq 0$;
 \item $X$ has log-terminal singularities if  $a_i >-1$;
 \item $X$ has log-canonical singularities if $a_i \geq -1$;
\end{itemize}
for all resolutions $f:Y\rightarrow X$.

It is customary to denote by $\mbox{Discrep}(p)$ (for a resolution of singularities a $p\in X$) the inf of all $a_i$ and all resolutions $f:Y \rightarrow X$. We also denote by $\mbox{Discrep}(X)$ the minimum of $\mbox{Discrep}(p)$ for all $p \in X$.

A characterization of log-terminal singularities in complex analytic geometry is the following: let $\Omega$ be a nowhere vanishing section of $K_X^k$ around $p\in Sing(X)$, then $p$ is a log-terminal singularity iff $\int_{U^{reg}} \Omega \wedge \overline{\Omega} < +\infty$, where $U$ denotes an open neighborhood of $p$ in $X$ and the integration is with respect to the Lebesgue measure.  
Moreover, recall that all log-terminal singularities are rational (i.e., $R^if_{\ast} \mathcal{O}_Y=0$ if $i >0$).

In dimension two terminal singularities are smooth points and canonical singularities are the classical RDP hypersurface singularities (also known as Du Val singularities or ADE singularities). They are obtained as  quotients of $\C^2$ by  finite subgroups $\Gamma$ of $SU(2)$ acting freely away from the origin. According to the classification of finite subgroups of $SU(2)$ there are only three families of singularities:
the cyclic quotients denoted by $A_k$ (order $k+1$), whose equations are given by $xy=z^{k+1}$; quotients by the binary dihedral group denoted by $D_k$ (order $4(k-2)$), with equations $x^2+y^2z+z^{k-1}=0$; the three exceptional cases $E_6$ (order $24$), with equation $x^2+y^3+z^4=0$, $E_7$ (order $48$) with equation $x^2+y^3+yz^3=0$ and $E_8$ (order $120$), with equation $x^2+y^3+z^5=0$.

More generally a two dimensional singularity is log-terminal iff it is a quotient of $\C^2$ by  a finite subgroup $\Gamma$ of $U(2)$ acting freely away from the origin.

The node $\{z_1^2+\dots+z_{n+1}^2=0\}\subseteq \C^{n+1}$ is an example of terminal (and non-orbifold) singularity (for $n\geq 3$). Three dimensional terminal singularities are also classified (M. Reid). In particular recall that they must be isolated. 

An important class of singularities are the ones which admit $\Q$-Gorenstein smoothings. In dimension two they are classified (see Chapter $2$).

\end{chapter} 
%Biblio

\bibliographystyle{plain}

\addcontentsline{toc}{chapter}{Bibliography}
\bibliography{CSThesis} % References file

\end{document}